\documentclass[journal]{new-aiaa}
\usepackage[utf8]{inputenc}

\usepackage{bm}
\usepackage{amsmath}
\usepackage{color,soul}
\usepackage{todonotes}
\usepackage{colortbl}
\definecolor{myblue}{rgb}{0, 0.23, 0.64}
\definecolor{WVUblue}{rgb}{0, 0.16, 0.33}

\usepackage{hyperref}
\hypersetup{
    colorlinks=true,
    linkcolor=myblue,
    filecolor=magenta,
    citecolor=myblue,
    urlcolor=myblue,
}
\usepackage{siunitx}
\usepackage{graphicx}
\usepackage[version=4]{mhchem}
\usepackage[group-separator={,}]{siunitx}
\usepackage{longtable,tabularx}
\usepackage{gensymb}
\usepackage{mathtools}
\usepackage{float}
\usepackage[ruled, vlined, linesnumbered]{algorithm2e}
\SetKwComment{Comment}{$\triangleright$\ }{}
\usepackage{subcaption}
\usepackage{enumitem}
\usepackage{multirow}
\usepackage{booktabs}
\usepackage{lscape}
\usepackage{mdframed}
\usepackage{adjustbox}
\usepackage{dsfont}
\usepackage{xspace}
\usepackage[flushleft]{threeparttable}
\newcommand{\MCRP}{\hyperlink{MCRP}{\textsf{MCRP}}\xspace}
\newcommand{\MP}{\hyperlink{MP}{\textsf{MP$(s)$}}\xspace}
\newcommand{\RHP}{\hyperlink{RHP}{\textsf{RHP$(s,L)$}}\xspace}

\let\svthefootnote\thefootnote
\newcommand\freefootnote[1]{%
  \let\thefootnote\relax%
  \footnotetext{#1}%
  \let\thefootnote\svthefootnote%
}

\begin{document}

\freefootnote{This paper is a substantially revised version of the Paper AAS 22-825, presented at the AAS/AIAA Astrodynamics Specialist Conference, Charlotte, NC, August 7-11, 2022. It offers new results and a better description of the materials.}

\title{Deterministic Multistage Constellation Reconfiguration Using Integer Programming and Sequential Decision-Making Methods}

\author{Hang Woon Lee\footnote{Assistant Professor, Department of Mechanical, Materials and Aerospace Engineering; hangwoon.lee@mail.wvu.edu. Member AIAA (Corresponding Author).}, David O. Williams Rogers\footnote{Ph.D. Student, Department of Mechanical, Materials and Aerospace Engineering, Student Member AIAA.}, and Brycen D. Pearl\footnote{Ph.D. Student, Department of Mechanical, Materials and Aerospace Engineering, Student Member AIAA.}}
\affil{West Virginia University, Morgantown, WV, 26506}
\author{Hao Chen\footnote{Assistant Professor, School of Systems and Enterprises, Member AIAA.}}
\affil{Stevens Institute of Technology, Hoboken, NJ, 07030}
\author{Koki Ho\footnote{Associate Professor, Daniel Guggenheim School of Aerospace Engineering, Senior Member AIAA.}}
\affil{Georgia Institute of Technology, Atlanta, GA, 30332}

\maketitle{}

\begin{abstract}
In this paper, we address the problem of reconfiguring Earth observation satellite constellation systems through multiple stages. The Multi-stage Constellation Reconfiguration Problem (MCRP) aims to maximize the total observation rewards obtained by covering a set of targets of interest through the active manipulation of the orbits and relative phasing of constituent satellites. In this paper, we consider deterministic problem settings in which the targets of interest are known \textit{a priori}. We propose a novel integer linear programming formulation for MCRP, capable of obtaining provably optimal solutions. To overcome computational intractability due to the combinatorial explosion in solving large-scale instances, we introduce two computationally efficient sequential decision-making methods based on the principles of a myopic policy and a rolling horizon procedure. The computational experiments demonstrate that the devised sequential decision-making approaches yield high-quality solutions with improved computational efficiency over the baseline MCRP. Finally, a case study using Hurricane Harvey data showcases the advantages of multi-stage constellation reconfiguration over single-stage and no-reconfiguration scenarios.
\end{abstract}

\section*{Nomenclature}
{\renewcommand\arraystretch{1.0}
\noindent\begin{longtable*}{@{}l @{\quad=\quad} l@{}}
    $N$ & number of stages \\
    $\mathcal{C}$ & constellation configuration \\
    $\mathcal{S}$ & set of stage indices (index $s$) \\
    $\mathcal{K}$ & set of satellite indices (index $k$) \\
    $\mathcal{J}$ & set of orbital slots (indices $i,j$) \\
    $\mathcal{P}$ & set of target points (index $p$) \\
    $\mathcal{T}$ & mission planning horizon (index $t$) \\
    $c$ & transfer cost \\
    $\pi$ & observation reward \\
    $r$ & coverage requirement threshold \\
    $V$ & visibility matrix \\
    $x$ & satellite transfer variable \\
    $y$ & coverage state variable \\
    $z$ & observation reward \\
\end{longtable*}}

\section{Introduction}

Earth observation (EO) is the act of gathering information about planetary phenomena to understand their underlying dynamics and impacts. The use of satellite systems enables remote sensing, providing observations at a global scale and high temporal frequency, as well as enabling a range of measurements from synthetic-aperture radar, radio frequency, and/or optical sensors. For instance, EO satellites have been instrumental in monitoring ecosystem changes through the observation of oceanic water circulation and salinity \cite{Chen2021Ocean, Irrgang2019Ocean}, land degradation \cite{Jong2011LandDeg}, natural disaster monitoring \cite{Guo2010Disaster}, sustainability of coastal regions \cite{Politi2020Coasts}, and emissions from biomass burning \cite{Duncan2003Biomass}. Furthermore, EO satellites have provided vital information on changes in societal development by monitoring agricultural droughts \cite{Srivastava2018Drought}, geological mapping \cite{VANDERMEER2012Geology}, disease spread and general public health \cite{Sogno2020Disease}, and land cover \cite{HANSEN2012LandCover}. The wide range of use cases provided by EO satellites proves paramount in the monitoring of current and future planetary phenomena.

Monolithic satellite systems are spatially and temporally limited by their orbital characteristics, leading to the popularity of distributed satellite systems comprised of multiple satellites. Satellite constellations, such as the Afternoon Train constellation operated for EO \cite{Kelly2009}, contain many satellites working toward a common goal, allowing for cooperation to capture the governing dynamics of targets that a monolithic satellite system would lack. For example, a multitude of satellites was used to gather data on hurricanes Katrina and Rita, such as sea surface temperature data gathered by NOAA polar orbiters NOAA-16, 17, and 18, sea surface height derived from altimeter measurements of Envisat, and chlorophyll data from NASA's Aqua \cite{Gierach2007}. Moreover, the Disaster Monitoring Constellation, operated for disaster response by a cooperative international team, contains a variety of satellites for varying purposes \cite{Stephens2003DMC}. The additional observational throughput provided by cooperative satellite systems allows for more informed analysis and decision-making related to natural disasters and other planetary phenomena. However, these traditional satellite constellations lack the ability to respond to highly dynamic mission environments and objectives, thereby limiting satellite performance when applied to fast-paced processes, resulting in long revisit times or a low quantity of data obtained.

Constellation reconfiguration, defined as the process of transforming a given configuration into another through the orbital maneuvers of satellites, is a state-of-the-art concept in satellite operation that improves upon the limitations imposed by fixed configuration systems \cite{deweck2008optimal, he2020reconfigurable, Zuo2022Surrogate, lee2023regional}. As a result of orbital maneuverability, constellation reconfiguration provides a high level of flexibility and responsiveness to satellite systems \cite{chen2015reconfiguration, paek2019optimization}, allowing for new operational or scientific tasks to be performed in an optimized constellation configuration. Applications of constellation reconfiguration include its use in telecommunications systems for staged deployment directed toward the minimization of transfer cost or overall cost \cite{deweck2008optimal, Appel2014Optimization, Anderson2022megaconstellation}. Additional investigations include reconfiguration in response to lost assets with the objective of compensating for lost performance \cite{ferringer2009many, Zuo2022Surrogate}.

Existing literature on constellation reconfiguration in application to EO demonstrates the value of maneuverable satellites. Extensive work has covered reconfiguration between two selected modes for constellation operation, one selected for global observation and one providing more frequent regional observations, with the objectives to maximize observations of a region of interest, minimize revisit time, and minimize reconfiguration time \cite{paek2019optimization, he2020reconfigurable}. In the case of Ref.~\cite{paek2019optimization}, the formulation of the reconfiguration process is thoroughly described in the application to a latitude of interest in the regional observation mode, and a band of latitudes in the global observation mode. Meanwhile, Ref.~\cite{he2020reconfigurable} applies the reconfiguration formulation to unpredictable disaster locations selected from a given set of possible locations, as well as to a forest fire. Additional research has been provided on altitude change maneuvers for response to mobile target tracking such as cyclones \cite{Morgan2023} or earthquake impact zones \cite{McGrath2019General}. The main objectives provided by altitude change maneuvers include the reduction of reconfiguration cost and overall revisit time. One further extension of altitude change reconfiguration incorporates single-stage maneuvers and phasing maneuvers with the objective of maximizing data gathered on, and rapid response to, coasts impacted recently by tsunamis \cite{Jiaxin2021Evolution}. Furthermore, single-stage reconfiguration has been investigated for application to natural disaster impact monitoring with emphasis on high levels of maneuverability to maximize observations \cite{Lee2020binary, Lee2021lagrangian, lee2023regional}. It should be noted that `stage,' in this context, refers to an opportunity for reconfiguration, differing from the definition of `stage' in Refs.~\cite{deweck2008optimal, Appel2014Optimization, Anderson2022megaconstellation} relating to stages of constellation deployment. Such applications of reconfiguration in EO provide significant results as opposed to fixed constellation configuration systems.

The extensive research conducted in previous literature has proved paramount to the overall investigation of constellation reconfiguration and the value it provides to EO, however, there are certain gaps that can be addressed and expanded upon. In the case of reconfiguration between global and regional modes \cite{paek2019optimization, he2020reconfigurable, McGrath2019General} the constellation is inherently limited to the two specific predetermined constellation modes, resulting in a lack of both individual satellite flexibility and overall constellation variability. Similarly, limiting satellite maneuvers to a change in altitude \cite{McGrath2019General, Morgan2023, Jiaxin2021Evolution} removes overall flexibility. Limitations to satellite and constellation flexibility, in the application of available maneuver types and satellite individuality respectively, restrict reconfiguration to an extremely rigid set of conditions. A more expanded set of conditions allowing wider flexibility provides a larger degree of freedom in the reconfiguration process, possibly benefiting objectives over the course of the mission lifetime. Finally, allowing only single-stage reconfiguration \cite{Jiaxin2021Evolution, Lee2020binary, Lee2021lagrangian, lee2023regional} limits observational throughput, especially over long time horizons that may benefit greatly from multiple reconfiguration opportunities.

The contributions of this paper are as follows: In response to the gaps identified in the literature, we propose a novel problem referred to as the \textit{Multi-stage Constellation Reconfiguration Problem} (MCRP). This problem extends our prior work on a single-stage reconfiguration problem \cite{lee2023regional} to provide an increased degree of freedom, enabling the uncovering of aspects of the design space that may remain unnoticed in zero or single reconfiguration stages. To model the sequence of orbital maneuvers by satellites over time, we apply the concept of TEGs. Based on this graph-theoretic modeling, we propose a novel integer linear programming (ILP) formulation for the MCRP, which can be solved using conventional mixed-integer linear programming (MILP) methods to obtain provably optimal solutions. Additionally, we devise two sequential decision-making methods to address the computational intractability found in solving large-scale problems. We empirically demonstrate that the proposed methods provide high-quality solutions and are computationally efficient. The versatility of the proposed MCRP framework is further demonstrated through a case study of tracking a real-world historical storm system, Hurricane Harvey. This paper extends a previous version of this research \cite{Lee2022maximizing} by providing new results and a more detailed description of the materials.

The remainder of this paper is organized as follows: In Sec.~\ref{sec:mcrp}, we provide a formal description of MCRP and propose a novel ILP formulation of it. Section~\ref{sec:solution} discusses two solution methods aimed at addressing large-scale MCRP instances. Then, in Sec.~\ref{sec:computational_experiments}, we conduct a comparative analysis to validate the applicability of the solution methods and present a case study of a real-world scenario, using Hurricane Harvey's historical trajectory, to demonstrate the value of multi-stage constellation reconfiguration as a means to increase the system's observational throughput. Lastly, in Sec.~\ref{sec:conclusions}, we suggest several interesting future work directions to enhance the applicability of the proposed work and conclude this paper.


\section{Multi-Stage Constellation Reconfiguration Problem} \label{sec:mcrp}
In this section, we describe and propose a mathematical optimization formulation of MCRP.

\subsection{Problem Description} \label{sec:description}
The mission planning horizon, denoted by $\mathcal{T}$, comprises discrete time steps $t=1,2,\ldots,T$. This horizon is finite and is evenly segmented into $N$ intervals, referred to as \textit{stages}. At each of these stages, a satellite constellation---a group of satellites---may undergo a reconfiguration process. We define reconfiguration as the process of transitioning a constellation system from one configuration to another through maneuverable satellites executing orbital transfers.

The objective of MCRP is to determine the optimal sequence of orbital maneuvers for satellites, aimed at maximizing the observation rewards obtained by covering targets of interest. These targets, which may be static or dynamic and are situated on the ground, in the air, or in space, are chosen in alignment with specific mission objectives or client requirements. This paper focuses on the deterministic variant of MCRP, where targets and their attributes, such as observation rewards and coverage thresholds, are known \textit{a priori} for the entire mission planning horizon.

The MCRP is subject to various physical and operational constraints, including the observation reward mechanism that governs the process of computing visible time windows (VTWs) and obtaining rewards, as well as budget constraints that dictate the feasibility of orbital transfer maneuvers. These constraints will be formally discussed later in this section. Furthermore, the MCRP must account for potential heterogeneity within the constellation, where each satellite can be equipped with unique sensors, propulsion subsystem specifications, propellant states, and orbital characteristics.

We define $\mathcal{K}$ as the set of satellites, indexed by $k = 1, 2, \ldots, K$. The target set is denoted by $\mathcal{P}$, where each target is indexed by $p$, and the total number of targets is $P$. Associated with each target are time-dependent observation rewards $\pi_{tp}$ and coverage thresholds $r_{tp}$. The set of stages is represented as $\mathcal{S} = \{0, 1, \ldots, N\}$, with each stage indexed by $s$ and the total number of reconfiguration stages being $N$ (note that $|\mathcal{S}| = N + 1$). In $\mathcal{S}$, $s = 0$ signifies the initial condition of the MCRP. Each stage $s \in \mathcal{S}$ is associated with a start time $t_s \in \mathcal{T}$ and a stage planning horizon $\mathcal{T}_s = \{t : t_s \leq t < t_{s+1}, t \in \mathcal{T}\}$. The mission planning horizon can also be represented as $\mathcal{T}=\{\mathcal{T}_1,\ldots,\mathcal{T}_N\}$.

In Sec.~\ref{sec:graph}, we introduce a graph-theoretic model of the reconfiguration process. Additionally, we introduce the observation reward mechanism, detailed in Sec.~\ref{sec:mechanism}.

\subsubsection{Modeling the Sequence of Orbital Maneuvers Using Time Expanded Graphs} \label{sec:graph}
To model the sequences of orbital transfers over time, we leverage the concept of a TEG. In a TEG, a vertex represents a specific state at a certain time, and an edge represents a transition from one state to another. The edges in TEGs are directed along the flow of time.

Applying this concept to our problem, each vertex in the graph corresponds to an \textit{orbital slot}, with its states---position and velocity vectors---defined at a given time. Each directed edge represents an \textit{orbital transfer} from one orbital slot to another over time. Throughout the paper, we refer to a \textit{path} as the sequence of vertices connected by directed edges, and to a \textit{reconfiguration process} as the set of these paths.

\begin{figure}[htb]
	\centering
		\includegraphics[width=0.45\linewidth]{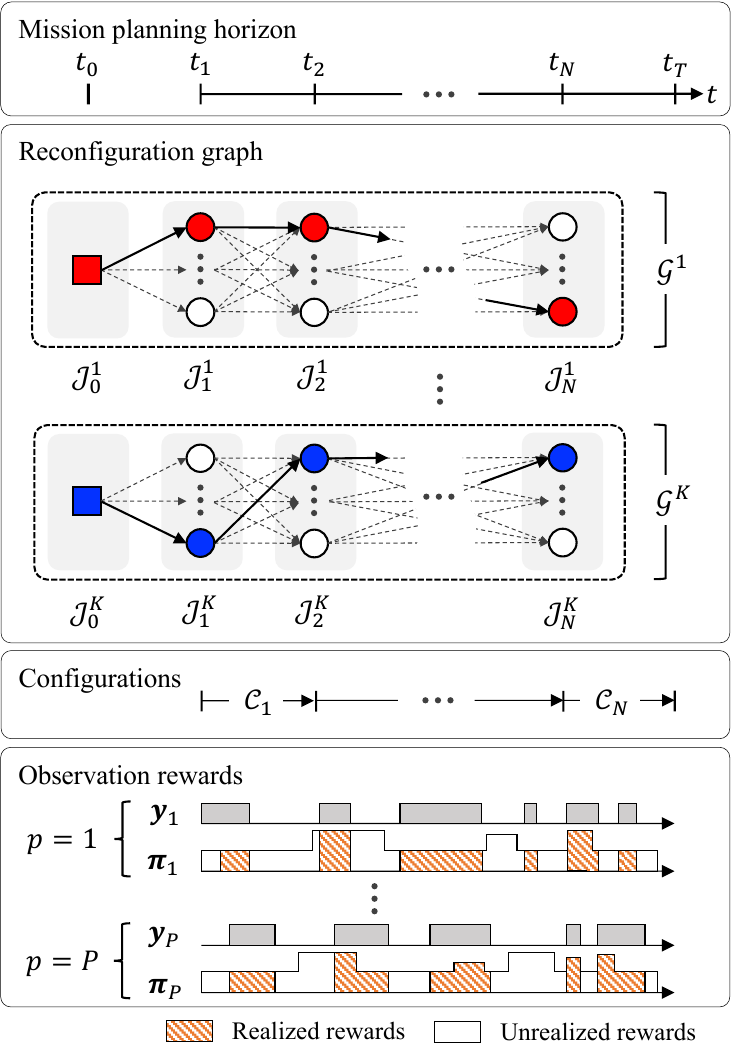}
		\caption{Mission planning horizon, reconfiguration graph, constellation configurations, and observation rewards.}
		\label{fig:graph}
\end{figure}

A \textit{reconfiguration graph} consists of $K$ TEGs, as depicted in Fig.~\ref{fig:graph}. We assume that each satellite $k$ is associated with its own TEG, denoted as $\mathcal{G}^k=(\mathcal{J}^k,\mathcal{E}^k)$, where $\mathcal{J}^k$ and $\mathcal{E}^k$ are the vertex set and edge set for satellite $k$, respectively. Let $\mathcal{J}_s^k$ denote the set of vertices for stage $s$ with the cardinality $J_s^k$; the vertex set can be further decomposed such that $\mathcal{J}^k=\{\mathcal{J}_0^k,\mathcal{J}_1^k,\ldots,\mathcal{J}_N^k\}$. Note that $\mathcal{J}_0^k$ is a singleton set, containing only one element: the orbital slot to which satellite $k$ is initialized. Similarly, let $\mathcal{E}_s^k$ denote the set of directed edges from orbital slots in $\mathcal{J}_{s-1}^k$ to orbital slots in $\mathcal{J}_{s}^k$; the edge set can be further decomposed such that $\mathcal{E}^k=\{\mathcal{E}_1^k,\mathcal{E}_2^k,\ldots,\mathcal{E}_N^k\}$. Decomposing the vertex and edge sets based on stages allows us to take into account generic mission scenarios in which satellites encounter different sets of orbital slots at every stage. In this paper, we assume that all orbital slots in $\mathcal{J}_s^k$ at $t_s$ are created by propagating orbital slots in $\mathcal{J}_{s-1}^k$ at $t_{s-1}$ forward in time and that all satellites have the identical number of orbital slots such that $J_s^k=J$ for all $s=1,\ldots,N$ and $k=1,\ldots,K$.

At each stage $s$, satellite $k$ performs an active orbital maneuver to transfer from the orbital slot $i \in \mathcal{J}_{s-1}^k$ of the previous stage to a new orbital slot $j \in \mathcal{J}_s^k$. We model the utilization of the transfer edge $(i,j)$ by satellite $k$ between stages $s-1$ and $s$ by defining the following binary decision variable:
\begin{equation}
    \label{eq:x}
    x_{ij}^{sk} = \begin{cases}
        1, & \text{if satellite $k$ transfers from orbital slot $i\in\mathcal{J}_{s-1}^k$ to orbital slot $j\in\mathcal{J}_s^k$} \\
        0, & \text{otherwise}
    \end{cases}
\end{equation}
where we refer to this as the \textit{satellite transfer variable}. Note that, in Eq.~\eqref{eq:x}, the stage index $s$ is defined with respect to the destination stage.

To model the cost incurred due to orbital transfers, we associate each edge $(i,j)\in\mathcal{E}_s^k$ with a non-negative transfer cost $c_{ij}^{sk}\ge0$. This cost represents the propellant or the $\Delta V$---the change in velocity---required to transfer satellite $k$ from orbital slot $i \in \mathcal{J}_{s-1}^k$ to orbital slot $j \in \mathcal{J}_s^k$. It can be computed by solving an orbital boundary value problem, given $i$ and $j$ as the initial and final conditions, respectively, along with the specifications of satellite $k$. In this paper, we assume that each vertex set $\mathcal{J}_s^k$ for $s=1,\ldots,N$ contains orbital slots that are propagated from the orbital slots of the previous stage to allow the possibility of remaining in the same orbit without any orbital maneuver. Consequently, there exists an edge $(i,j)\in\mathcal{E}_s^k$ for all $s\in\mathcal{S}\setminus\{0\}$ and for all $k\in\mathcal{K}$, with $c_{ij}^{sk}=0$. An edge is considered feasible if its transfer cost is less than the satellite's maximum available resource, $c_{\max}^{k}$.

\begin{figure}[htbp]
	\centering
		\includegraphics[width=0.31\linewidth]{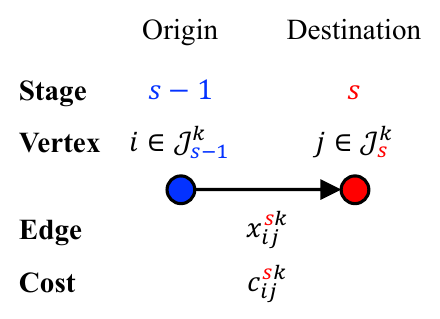}
		\caption{Vertex-edge relationship.}
		\label{fig:vertexedge}
\end{figure}

Figure~\ref{fig:vertexedge} visualizes the relationship between two adjacent vertices and a directed edge that connects them.

For each pair of stage $s$ and satellite $k$, only one $x_{ij}^{sk}$ should take the value of one. We enforce this condition using a set of constraints when formulating the MCRP later in Sec.~\ref{sec:formulation}.

\subsubsection{Constellation Configuration, Coverage, and Observation Reward Mechanisms} \label{sec:mechanism}
Upon reconfiguration, at each stage, a new constellation \textit{configuration} is formed. A configuration for stage $s$ can be represented as a set comprising orbital slots newly occupied by all satellites, or mathematically as $\mathcal{C}_s := \{j:x_{ij}^{sk}=1, k\in\mathcal{K}, i\in\mathcal{J}_{s-1}^k, j\in\mathcal{J}_{s}^k\}$. This configuration remains valid for the time interval $[t_s,t_{s+1})$, as shown in Fig.~\ref{fig:graph}.

Each orbital slot is associated with a Boolean, time-dependent visibility profile for each target. The visibility profile, $V_{tjp}^{sk}$, is set to one (true) if target $p$ is visible from orbital slot $j\in\mathcal{J}_s^k$ at time step $t\in\mathcal{T}_s$, and zero (false) otherwise. This is obtained by propagating satellite $k$ placed in orbital slot $j\in\mathcal{J}_s^k$ forward in time and checking whether the satellite can view target $p$ at time step $t\in\mathcal{T}_s$.

Depending on how satellites are relatively phased in a configuration, it is possible for multiple satellites to simultaneously view target $p$. We consider target $p$ to be covered at time step $t\in\mathcal{T}_s$ if at least $r_{tp}^s$ satellites are viewing it simultaneously, where $r_{tp}^s \in \mathbb{N}$ represents the minimum number of satellites required for coverage, also known as the \textit{coverage threshold}. Here, $\mathbb{N}$ is the set of natural numbers. Although a common scenario requires single-fold continuous coverage, implying $r_{tp}^s = 1, \forall s \in \mathcal{S}\setminus\{0\}, \forall t \in \mathcal{T}_s$, allowing the threshold to take a natural number accommodates more general cases, such as multi-fold time-varying coverage.

Using the definitions of the satellite transfer variable $x_{ij}^{sk}$, the visibility profile $V_{tjp}^{sk}$, and the coverage threshold $r_{tp}^s$, we can mathematically express the following indicator variable, representing the Boolean \textit{coverage state} of a target:
\begin{equation}
        \label{eq:y}
        y_{tp}^s = \begin{cases}
        1, & \text{if } \displaystyle \sum_{k\in\mathcal{K}} \sum_{i\in\mathcal{J}_{s-1}^k} \sum_{j\in\mathcal{J}_s^k} V_{tjp}^{sk} x_{ij}^{sk} \ge r_{tp}^s\\
        0, & \text{otherwise}
    \end{cases}
\end{equation}
where we refer to this as the \textit{coverage state variable}. In Eq.~\eqref{eq:y}, the satellite transfer variable $x_{ij}^{sk}$ allows for the summation of the visibility profiles of orbital slots that are actually occupied by satellites. Using Eq.~\eqref{eq:y}, the VTW of target $p$ is represented as $\bm{y}_p=(y_{tp}^s\in\{0,1\}:s\in\mathcal{S}\setminus\{0\},t\in\mathcal{T}_s)$.

\begin{figure}[htbp]
	\centering
		\includegraphics[width=0.8\linewidth]{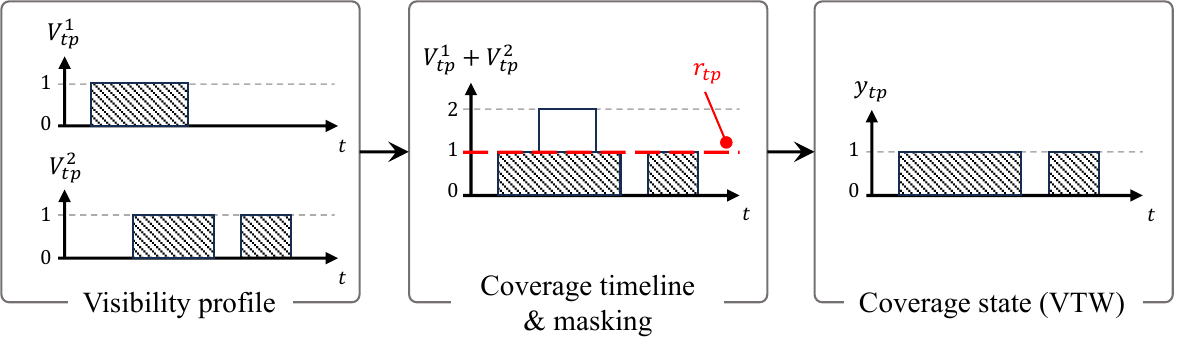}
		\caption{Illustration of the coverage mechanism.}
		\label{fig:masking}
\end{figure}

To illustrate the coverage mechanism, consider a simple case involving two satellites, $k=1,2$, each occupying one orbital slot, and a single target, $p$, as depicted in Fig.~\ref{fig:masking}. Each orbital slot has a distinct, time-dependent visibility profile of the target, $V_{tp}^{1}$ and $V_{tp}^{2}$, as shown in the left figure. Note that the indices $s$ and $j$ are omitted for ease of illustration. By aggregating these visibility profiles, a coverage timeline is generated that encodes the number of satellites in view of the target at each time step, as shown in the center figure. It can be seen that there are time intervals when both satellites are simultaneously viewing the target, where $V_{tp}^{1} + V_{tp}^{2} = 2$. The coverage requirement for the target, in this case, is single-fold continuous coverage, represented by a red dashed line. Applying coverage masking, the Boolean coverage state for the target, $y_{tp}$, is retrieved, indicating whether the coverage requirement is met, as shown in the right figure: one (true) if the target is covered, and zero (false) otherwise.

The observation reward, denoted as $\pi_{tp}^s$, is realized only when target $p$ is covered at time step $t\in\mathcal{T}_s$, specifically when $y_{tp}^s=1$, as depicted in Fig.~\ref{fig:graph}. To maximize the attainment of observation rewards over the entire mission planning horizon, expressed as $\sum_{s\in\mathcal{S}\setminus{0}} \sum_{t\in\mathcal{T}_s} \sum_{p\in\mathcal{P}} \pi_{tp}^s y_{tp}^s$, it is crucial to devise an optimal reconfiguration process. Optimal reconfiguration processes would be those that optimize the design of constellation configurations to ensure the alignment of VTWs coincides with periods yielding the highest observation rewards. However, solving the MCRP presents a significant computational challenge due to the complexity of finding a feasible and optimal solution. It demands consideration of various factors, including the coverage mechanism described above, multiple targets, and their respective time-dependent rewards. Additionally, it is subject to propellant budget constraints that limit the use of certain costly transfer edges, which could render configurations suboptimal and consequently result in suboptimal VTW-reward alignments. These complexities underscore the need to frame the MCRP as an optimization problem.

\subsection{Mathematical Formulation} \label{sec:formulation}
MCRP, in a general sense, can be classified as a multi-period deterministic decision-making problem. Although its reconfiguration processes are modeled on TEGs, MCRP, which involves maximizing observation rewards, cannot be fully represented within a graph-theoretic framework due to the interdependency between satellites and their interactions on targets via the coverage mechanism discussed in Sec.~\ref{sec:mechanism}. Therefore, conventional graph-based algorithms such as Dijkstra's algorithm and the A$^\ast$ algorithm, typically used for shortest path problems (based on the negation of the problem objective), are inadequate for solving MCRP.

To overcome the aforementioned mathematical programming challenges, we propose formulating MCRP as an ILP optimization problem. In the rest of this section, we define sets, parameters, and decision variables in Table~\ref{tab:symbol} and then introduce the ILP formulation of the MCRP.

\begin{table}[htbp]
    \caption{Definitions of sets, indices, parameters, and decision variables.}
    \label{tab:symbol}
    \centering
    \begin{tabular}{ll}
    \hline
    \hline
    Symbol & Description \\
    \hline
    \multicolumn{2}{l}{Sets and indices} \\
    $\mathcal{S}$ & Set of stage indices (index $s$; cardinality $N+1$) \\
    $\mathcal{K}$ & Set of satellite indices (index $k$; cardinality $K$) \\
    $\mathcal{J}_s^k$ & Set of orbital slot indices of stage $s$ for satellite $k$ (indices $i,j$; cardinality $J_s^k$) \\
    $\mathcal{P}$ & Set of target point indices (index $p$; cardinality $P$) \\
    $\mathcal{T}$ & Mission planning horizon (index $t$; cardinality $T$) \\
    $\mathcal{T}_s$ & Planning horizon for stage $s$ (index $t$) \\
    \multicolumn{2}{l}{Parameters} \\
    $c_{ij}^{sk}\ge0$ & Cost of transferring satellite $k$ from orbital slot $i\in\mathcal{J}_{s-1}^k$ to orbital slot $j\in\mathcal{J}_s^k$ \\
    $c_{\max}^k\ge0$ & Maximum resource available for satellite $k$ \\
    $\pi_{tp}^s\ge0$ & Observation reward for target point $p$ at time step $t\in\mathcal{T}_s$ \\
    $r_{tp}^s\in\mathbb{N}$ & Minimum number of satellites required to receive the reward of target point $p$ at time step $t\in\mathcal{T}_s$ \\
    $V_{tjp}^{sk}$ & $\begin{cases}
      1, & \text{if satellite $k$, placed in orbital slot $j\in\mathcal{J}_s^k$, is visible from target point $p$ at time step $t\in\mathcal{T}_s$} \\
      0, & \text{otherwise}
    \end{cases}$ \\
    \multicolumn{2}{l}{Decision variables} \\
$x_{ij}^{sk}$ & $\begin{cases}
1, &\text{if satellite $k$ transfers from orbital slot $i\in\mathcal{J}_{s-1}^k$ to orbital slot $j\in\mathcal{J}_s^k$} \\
0, &\text{otherwise}
\end{cases}$ \\
    $y_{tp}^s$ & $\begin{cases}
1, &\text{if target point $p$ is covered at time step $t\in\mathcal{T}_s$} \\
0, &\text{otherwise}
\end{cases}$ \\
    \hline
    \hline
    \end{tabular}
\end{table}

Building on the problem description in Sec.~\ref{sec:description} and the notations used therein (summarized in Table~\ref{tab:symbol}), the \MCRP can be formulated as an integer linear program as follows:\hypertarget{MCRP}{}
\begin{subequations}
\begin{alignat}{2}
z=\text{maximize} \quad & \sum_{s\in\mathcal{S}\setminus\{0\}} \sum_{t\in\mathcal{T}_s} \sum_{p\in\mathcal{P}} \pi_{tp}^s y_{tp}^s \label{mcrp:obj} \\
\text{subject to} \quad & \sum_{j\in\mathcal{J}_{1}^k} x_{ij}^{1k} = 1, \quad && \forall k\in\mathcal{K}, i\in\mathcal{J}_0^k \label{mcrp:init} \\
& \sum_{j\in\mathcal{J}_{s+1}^k}x_{ij}^{s+1,k}-\sum_{q\in\mathcal{J}_{s-1}^k}x_{qi}^{sk} = 0, && \forall s\in\mathcal{S}\setminus\{0,N\}, \forall k\in\mathcal{K}, \forall i\in\mathcal{J}_{s}^k \label{mcrp:flow} \\
& \sum_{k\in\mathcal{K}} \sum_{i\in\mathcal{J}_{s-1}^k} \sum_{j\in\mathcal{J}_{s}^k} V_{tjp}^{sk} x_{ij}^{sk} \ge r_{tp}^s y_{tp}^s, \quad && \forall s \in\mathcal{S}\setminus\{0\}, \forall t\in\mathcal{T}_s, \forall p\in\mathcal{P} \label{mcrp:cov} \\
& \sum_{s\in\mathcal{S}\setminus\{0\}} \sum_{i\in\mathcal{J}_{s-1}^k} \sum_{j\in\mathcal{J}_{s}^k} c_{ij}^{sk} x_{ij}^{sk} \le c_{\max}^{k}, && \forall k \in\mathcal{K}^\prime \label{mcrp:res} \\
& x_{ij}^{sk} \in \{0,1\}, && \forall s\in\mathcal{S}\setminus\{0\}, \forall k\in\mathcal{K}, \forall i\in\mathcal{J}_{s-1}^k, \forall j\in\mathcal{J}_{s}^k \label{mcrp:x} \\
& y_{tp}^s \in \{0,1\}, && \forall s\in\mathcal{S}\setminus\{0\}, \forall t\in\mathcal{T}_s, \forall p\in\mathcal{P} \label{mcrp:y}
\end{alignat}
\end{subequations}

The \MCRP aims to maximize the objective function given by Eq.~\eqref{mcrp:obj}, which represents the total observation reward obtained by covering a set of target points of interest over the entire mission planning horizon $\mathcal{T}$. We denote the optimal objective value of \MCRP by $z$. Constraints~\eqref{mcrp:init} and \eqref{mcrp:flow} represent the \textit{path contiguity constraints} that guarantee the existence of a connected path in each satellite's TEG. Specifically, constraints~\eqref{mcrp:init} ensure that all satellites are assigned new orbital slots for the first stage, while constraints~\eqref{mcrp:flow} balance the outflow (the first term) and inflow (the second term) for the vertices of intermediate stages $s=1,\ldots,N-1$, indicating that each path must be contiguous. Constraints~\eqref{mcrp:cov} are the \textit{configuration-coverage linking constraints}, ensuring that target point $p$ is covered at time step $t$ only if at least $r_{tp}$ satellites are in view. Typically, $r_{tp}=1$ for all $t$ in $\mathcal{T}$ (\textit{i.e.}, at least one satellite is needed for coverage). Constraints~\eqref{mcrp:res} are the \textit{budget constraints}, ensuring that the sum of all $\Delta V$ (or propellant) used by satellite $k$ does not exceed $c_{\max}^{k}$. The set $\mathcal{K}^\prime \subseteq \mathcal{K}$ denotes the subset of satellites subject to these resource availability constraints, although it is likely that all satellites are subject to such constraints. Finally, constraints~\eqref{mcrp:x} and \eqref{mcrp:y} define the domains of the decision variables as in Eqs.~\eqref{eq:x} and \eqref{eq:y}, respectively.

The \MCRP formulation is general enough to handle features such as heterogeneity in satellite hardware specifications and orbital characteristics, as well as asymmetry in satellite distribution. That is, the MCRP does not impose any prescribed constraints on the set of satellites $\mathcal{K}$ and the set of orbital slots $\mathcal{J}$, and thus can be defined according to users' needs. The consideration of heterogeneity can be particularly useful in modeling federated satellite-based emergency mapping missions, such as the International Charter: Space and Major Disasters, which facilitates coordination between space-based EO assets with different specifications for increased availability and diverse data products for disaster monitoring and supporting relief activities \cite{denis2016,voigt2016global}. Asymmetric satellite distributions can lead to highly efficient constellation pattern sets applicable to EO systems, as they can result in fewer satellites being required to provide the same coverage performance as their symmetric counterparts \cite{lee2020satellite,david2023designing}, which is particularly advantageous for monitoring regional targets.

Several remarks can be made about the presented \MCRP formulation. Constraints~\eqref{mcrp:res} couple the stages, while constraints~\eqref{mcrp:cov} couple the satellites. First, without constraints~\eqref{mcrp:res}, each individual stage can be solved independently and in parallel; the optimal solution for stage $s$ does not impact those of its subsequent stages $s+1,\ldots,N$. Second, without constraints~\eqref{mcrp:cov}, each satellite can occupy orbital slots that maximize its sum of observation rewards over the entire time horizon, without necessarily considering the loss of observation rewards due to coverage overlaps with other satellites, as illustrated in Fig.~\ref{fig:masking}. These two sets of constraints inherently make the \MCRP challenging to solve, as all decisions over all stages and all satellites need to be concurrently derived to obtain globally optimal solutions.

Each stage of reconfiguration entails optimizing two elements: (i) the design of a maximum-reward destination configuration and (ii) the assignment of satellites from one configuration to another. \MCRP generalizes the single-stage design-transfer problem, as discussed in Ref.~\cite{lee2023regional}, by extending it along the temporal dimension, allowing for multiple reconfiguration opportunities. The \MCRP framework is versatile enough to be applied to scheduling either a short segment of the entire mission horizon or the entire mission horizon itself, depending on the users' requirements. Multiple, short-segmented \textsf{MCRP}s can be consecutively applied, and their solutions can be concatenated to generate a solution for a longer planning horizon. However, this approach does not ensure optimality across the entire mission.

\section{Sequential Decision-Making Methods for Large-Scale MCRP Instances} \label{sec:solution}
In Sec.~\ref{sec:mcrp}, the ILP formulation of MCRP is proposed. This approach enables the utilization of generic MILP methods, such as the branch-and-bound algorithm, to solve the problem and obtain provably optimal solutions, which are well-suited for small-scale instances. Several commercial software packages, such as the Gurobi Optimizer, CPLEX, and MATLAB's \texttt{intlinprog}, are available for ease of implementation.

Despite these advantages, \MCRP, however, is a combinatorial optimization problem that suffers from an exponentially expanding solution space when attempting to optimize all stage decisions simultaneously; the total number of potentially feasible plans increases exponentially with a linear increase in $J$, $N$, and $K$. There are at most $J^{NK}$ reconfiguration plans to consider (obtained without considering the budget constraints). For example, an instance $I$ of \MCRP with three reconfiguration stages, five satellites, and fifty candidate orbital slots per satellite can have up to \num{3.05e25} potentially feasible reconfiguration plans. Both exhaustive enumeration of these plans to check their feasibility and optimality, and the use of MILP methods, are computationally prohibitive for large-scale instances.

To address the computational intractability challenge posed by the exponentially expanding combinatorial solution space when attempting to optimize all stages simultaneously for large-scale instances of \MCRP, we propose two sequential decision-making solution methods. These methods are based on the principles of a myopic policy and the rolling horizon procedure \cite{powell2009}. The main idea behind these methods is to partition the original problem, \MCRP, into smaller subproblems by stages, making them more manageable in size. These subproblems are then solved sequentially, stored, and later aggregated and processed to form feasible solutions to \MCRP.

\subsection{Myopic Policy}
The proposed \textit{Myopic Policy} (MP) method partitions an $N$-stage \MCRP into $N$ coupled single-stage subproblems. Each subproblem aims to optimize the assignment of satellites such that the total observation rewards obtained are maximized for a given stage, without considering future stages. Each subproblem is sequentially solved, with the results from the prior stage being passed down to the subsequent subproblem as fixed parameters. This approach has a significant benefit: although $N$ single-stage subproblems must be solved in series to obtain a feasible solution to \MCRP, each subproblem has only $J^K$ potentially feasible reconfiguration plans to consider, which is significantly smaller in scale than that of \MCRP. For example, if we consider the same instance $I$ of \MCRP discussed earlier, we can partition it into three subproblems, each with up to \num{3.12e8} potentially feasible plans. These smaller subproblems can be efficiently solved using commercial software packages. Moreover, additional algorithmic improvements can be applied to subproblems to enhance performance in terms of solution quality and time complexity.

\begin{figure}[htbp]
	\centering
		\includegraphics[width=0.25\linewidth]{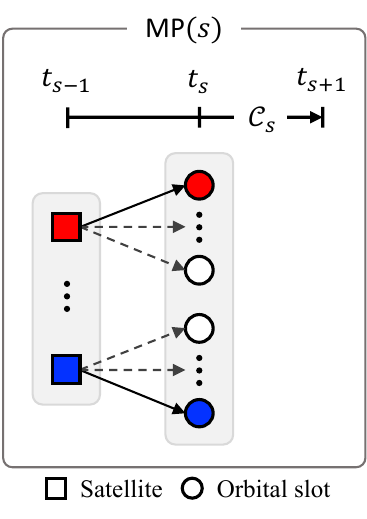}
		\caption{MP subproblem for stage $s$.}
		\label{fig:mp}
\end{figure}

A subproblem \MP is parameterized with the stage index $s=1,\ldots,N$ and the associated origin orbital slot index $i$. The destination orbital slots from the previous stage, $s-1$, become the origin orbital slots for the stage $s$ subproblem, thus parameterizing $i\in\tilde{\mathcal{J}}_{s-1}^k$. Here, $\tilde{\mathcal{J}}_{s-1}^k$ is defined as follows:
\begin{equation}
    \label{eq:priorvertex}
    \tilde{\mathcal{J}}_{s-1}^k := \left\{j : \tilde{x}_{ij}^{s-1,k}=1,i\in\mathcal{J}_{s-2}^k,j\in\mathcal{J}_{s-1}^k\right\}
\end{equation}
where the tilde symbol $(\tilde{\cdot})$ is used to denote fixed parameters and variables from previous stages. Also, note that $\tilde{\mathcal{J}}_{s-1}^k$ is a singleton set, containing the destination orbital slot from the previous stage as its only element; $\tilde{\mathcal{J}}_0^k$ corresponds to the initial conditions of the MCRP. Figure~\ref{fig:mp} illustrates the scope of \MP.

The mathematical formulation of \MP is as follows:\hypertarget{MP}{}
\begin{subequations}
\begin{alignat}{2}
z_s = \text{maximize} \quad & \sum_{t\in\mathcal{T}_s}\sum_{p\in\mathcal{P}} \pi_{tp}y_{tp} \label{eq:mh1} \\
\text{subject to} \quad & \sum_{j\in\mathcal{J}_{s}^k} x_{ij}^{sk} = 1, && \forall k\in\mathcal{K}, i\in\tilde{\mathcal{J}}_{s-1}^k \label{eq:mh2} \\
& \sum_{k\in\mathcal{K}} \sum_{j\in\mathcal{J}_{s}^k} V_{tjp}^{sk} x_{ij}^{sk} \ge r_{tp}^s y_{tp}^s, \quad && \forall t\in\mathcal{T}_s,\forall p\in\mathcal{P} \label{eq:mh4} \\
& \sum_{j\in\mathcal{J}_{s}^k} c_{ij}^{sk} x_{ij}^{sk} \le c_{\max}^{sk}, && \forall k \in\mathcal{K}^\prime \label{eq:mh5} \\
& x_{ij}^{sk} \in \{0,1\}, && \forall k\in\mathcal{K}, \forall j\in\mathcal{J}_{s}^k \label{eq:mh6} \\
& y_{tp}^{s} \in \{0,1\}, && \forall t\in\mathcal{T}_s, \forall p\in\mathcal{P} \label{eq:mh7}
\end{alignat}
\end{subequations}

The objective of \MP is to maximize the objective function~\eqref{eq:mh1}, which represents the total observation reward for stage $s$. We denote the optimal objective value of \MP by $z_s$. Constraints~\eqref{eq:mh2} ensure that all satellites are assigned to new orbital slots. Constraints~\eqref{eq:mh4} are the configuration-coverage linking constraints for stage $s$. Constraints~\eqref{eq:mh5} represent the budget constraints. Unlike constraints~\eqref{mcrp:res}, constraints~\eqref{eq:mh5} are bounded by $c_{\max}^{sk}$, which reflects the remaining budget (\textit{i.e.}, delta-v) for satellite $k$ at stage $s$. This can be computed by considering the resource consumption from the previous stages, $1,\ldots,s-1$, as follows:
\begin{equation}
    \label{eq:csmax}
     c_{\max}^{sk} = c_{\max}^{k}-\sum_{q=1}^{s-1}\sum_{i\in\mathcal{J}_{q-1}^k}\sum_{j\in\mathcal{J}_{q}^k} c_{ij}^{qk} \tilde{x}_{ij}^{qk}
\end{equation}
Constraints~\eqref{eq:mh6} and \eqref{eq:mh7} define the domains of the decision variables. Note that both $x_{ij}^{sk}$ and $y_{tp}^s$ are essentially two-dimensional variables because $s$ and $i$ are now fixed; $k$ and $j$ are the only independent indices.

Algorithm~\ref{alg:mp} outlines the overall solution procedure of MP, which solves \MP sequentially, progressing stage-by-stage in a greedy manner. In each stage, \MP yields the optimal objective value $z_s$, the optimal assignment solution $\bm{x}_s^\ast=(x_{ij}^{sk\ast}=\{0,1\}: k\in\mathcal{K}, i\in\tilde{\mathcal{J}}_{s-1}^k, j\in\mathcal{J}_{s}^k)$, and the optimal coverage state solution $\bm{y}_s^\ast=(y_{tp}^{s\ast}\in\{0,1\}:t\in\mathcal{T}_s,p\in\mathcal{P})$. It is important to note that the optimality here pertains to \MP, not to \MCRP. Upon enumerating through all stages from $s=1$ to $s=N$, the algorithm outputs the overall solution objective value $z_{\text{mp}}$, which is the aggregate of all $z_s$ values, and a feasible solution $(\bm{x}^\ast,\bm{y}^\ast)$, which is the set of $(\bm{x}_s^\ast,\bm{y}_s^\ast)$ for $s=1,\ldots,N$, for the original problem, \MCRP.

\begin{algorithm}[htbp]
\DontPrintSemicolon
\caption{Myopic policy}
\label{alg:mp}
\For{$s=1,\ldots,N$}{
\If{$s\ge2$}
{
Update $\tilde{\mathcal{J}}_{s-1}^k$ as in Eq.~\eqref{eq:priorvertex} \;
}
Solve \MP and store: $z_s$ and $(\bm{x}_s^\ast,\bm{y}_s^\ast)$ \;
}
\Return $z_{\text{mp}}=\sum_{s=1}^N z_s$ and $(\bm{x}^\ast,\bm{y}^\ast) = \left((\bm{x}_s^\ast,\bm{y}_s^\ast):s\in\mathcal{S}\setminus\{0\}\right)$ \;
\end{algorithm}

The formulation of \MP is mathematically equivalent to the Regional Constellation Reconfiguration Problem with Individual Resource Constraints (RCRP-IRC), as described in Ref.~\cite{lee2023regional}. This problem embeds an assignment problem with a budget constraint and a maximal covering location problem. As previously noted, any dedicated algorithm can be applied to \MP to improve both solution quality and computational efficiency. In Ref.~\cite{lee2023regional}, a Lagrangian relaxation-based solution method is used to approach large-scale RCRP-IRCs.

\subsection{Rolling Horizon Policy}
The \textit{Rolling Horizon Policy} (RHP), also known as the receding horizon procedure or model predictive control in operations research and control literature, utilizes the impact of current-stage decisions on future stages to make informed decisions in the current stage \cite{Sethi1991theory}. It comprises \textit{control} stages and \textit{lookahead} stages. The control stages are those for which solutions are derived and kept to make up a solution for the original problem, while the lookahead stages aid in decision-making at the control stage by providing deterministic forecasts of future stages. In this paper, we consider a 1-stage control and an $L$-stage deterministic lookahead RHP. Consequently, the RHP partitions an $N$-stage \MCRP into $N-L$ coupled subproblems, each resembling a smaller-scale \MCRP. Due to the lookahead, each subproblem is larger than a subproblem in \MP, but smaller than the original \MCRP, with at most $J^{(L+1)K}$ potentially feasible plans. Considering the same instance $I$ of \MCRP and implementing a 1-stage lookahead RHP, each subproblem can have up to \num{9.77e16} potentially feasible plans.

\begin{figure}[htbp]
	\centering
		\includegraphics[width=0.5\linewidth]{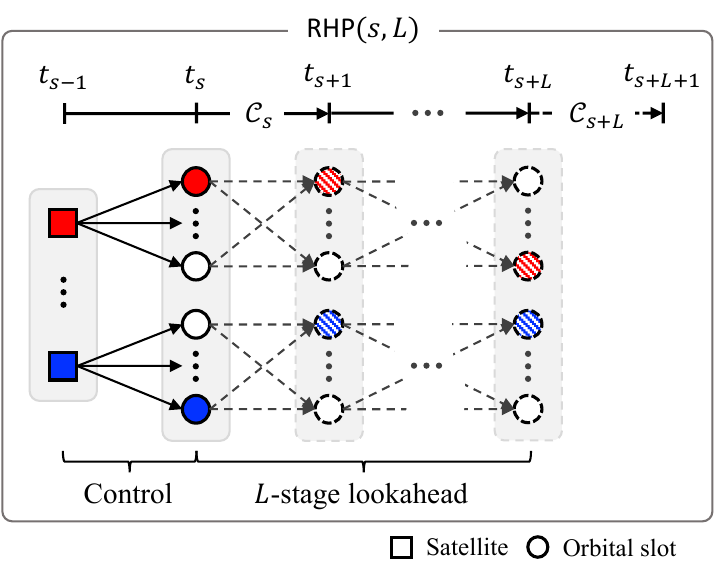}
		\caption{RHP subproblem for stage $s$ with the $L$-stage lookahead.}
		\label{fig:rhp}
\end{figure}

We denote by \RHP a subproblem parameterized with the stage index $s$ and the $L$-stage deterministic lookahead. As previously discussed in MP, associated with the stage parameter $s$ is the origin orbital slot $i \in \tilde{\mathcal{J}}_{s-1}^k$. At stage $s$, \RHP concurrently optimizes stages $s$ through $s+L$, with stage $s$ as the control and stages $s+1$ to $s+L$ as the lookahead. Figure~\ref{fig:rhp} illustrates the scope of \RHP.

The mathematical formulation of \RHP is as follows:\hypertarget{RHP}{}
\begin{subequations}
\begin{alignat}{2}
\sum_{\ell=s}^{s+L} z_{\ell} = \text{maximize} \quad & \sum_{\ell=s}^{s+L} \sum_{t\in\mathcal{T}_{\ell}} \sum_{p\in\mathcal{P}} \pi_{tp}^{\ell} y_{tp}^{\ell} \label{eq:rhp1} \\
\text{subject to} \quad & \sum_{j\in\mathcal{J}_{s}^k} x_{ij}^{sk} = 1, && \forall k\in\mathcal{K}, i\in\tilde{\mathcal{J}}_{s-1}^k \label{eq:rhp2} \\
& \sum_{j\in\mathcal{J}_{\ell+1}^k}x_{ij}^{\ell+1,k}-\sum_{q\in{\mathcal{J}}_{\ell-1}^k}x_{qi}^{\ell k} = 0, && \forall \ell\in\{s,\ldots,s+L-1\}, \forall k\in\mathcal{K}, \forall i\in\mathcal{J}_{\ell}^k \label{eq:rhp3} \\
& \sum_{k\in\mathcal{K}} \sum_{i\in{\mathcal{J}}_{\ell-1}^k} \sum_{j\in\mathcal{J}_{\ell}^k} V_{tjp}^{\ell k} x_{ij}^{\ell k} \ge r_{tp}^\ell y_{tp}^\ell, \quad && \forall \ell\in\{s,\ldots,s+L\}, \forall t\in\mathcal{T}_{\ell}, \forall p\in\mathcal{P} \label{eq:rhp4} \\
& \sum_{\ell=s}^{s+L} \sum_{i\in{\mathcal{J}}_{\ell-1}^k} \sum_{j\in\mathcal{J}_{\ell}^k} c_{ij}^{\ell k} x_{ij}^{\ell k} \le c_{\max}^{sk}, && \forall k \in\mathcal{K}^\prime \label{eq:rhp5} \\
& x_{ij}^{\ell k} \in \{0,1\}, && \forall \ell\in\{s,\ldots,s+L\}, \forall k\in\mathcal{K}, \forall i\in\mathcal{J}_{\ell-1}^k, \forall j\in\mathcal{J}_{\ell}^k \label{eq:rhp6} \\
& y_{tp}^\ell \in \{0,1\}, && \forall \ell\in\{s,\ldots,s+L\}, \forall t\in\mathcal{T}_{\ell}, \forall p\in\mathcal{P} \label{eq:rhp7}
\end{alignat}
\end{subequations}

The objective of \RHP is to maximize the objective function~\eqref{eq:rhp1}, which represents the sum of observation rewards from stages $s$ to $s+L$. Here, $\ell$ is used to index the stages. The optimal objective value of \RHP is equal to $z_{s}+\cdots+z_{s+L}$. Constraints~\eqref{eq:rhp2} and \eqref{eq:rhp3} represent the path contiguity constraints, as previously discussed for \MCRP. Constraints~\eqref{eq:rhp4} are the configuration-coverage linking constraints for all stages from $s$ to $s+L$. Constraints~\eqref{eq:rhp5} represent the budget constraints; $c_{\max}^{sk}$ is defined in Eq.~\eqref{eq:csmax}. Constraints~\eqref{eq:rhp6} and \eqref{eq:rhp7} define the domains of the decision variables.

\begin{figure}[htbp]
	\centering
		\includegraphics[width=0.725\linewidth]{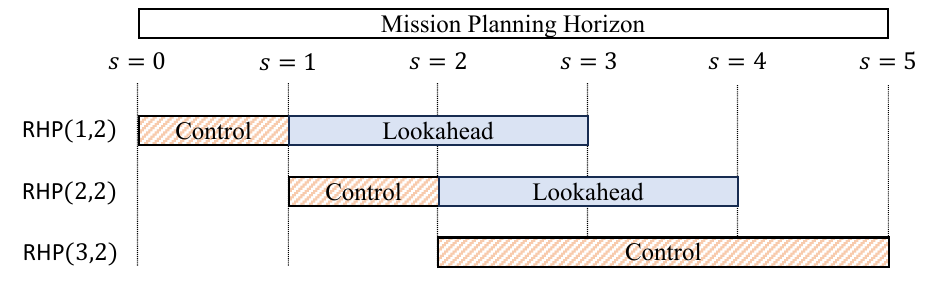}
		\caption{Example RHP with $N=5$ and $L=2$.}
		\label{fig:rhp_eg}
\end{figure}

Algorithm~\ref{alg:rhp} provides an overview of the RHP. The \RHP solutions are the satellite transfer and coverage state variables of the control stage $s$, which are represented as $(\bm{x}_s^\ast,\bm{y}_s^\ast)$. As mentioned, the lookahead stage solutions, $(\bm{x}_{s+1}^\ast,\bm{y}_{s+1}^\ast),\ldots,(\bm{x}_{s+L}^\ast,\bm{y}_{s+L}^\ast)$, are not utilized. Similarly, only the optimum of the control stage, $z_s$, is stored, which is computed as
\begin{equation}
    \label{eq:zs}
    z_s = \sum_{t\in\mathcal{T}_s}\sum_{p\in\mathcal{P}} \pi_{tp}^s y_{tp}^{s\ast}
\end{equation}
If $s=N-L$, the entire remaining mission planning horizon can be deterministically solved, making all remaining stages the control. A set of all control solutions is a feasible solution to \MCRP, and the optimum obtained by the RHP algorithm is the sum of all $z_s$ from $s=1$ to $s=N$. To illustrate the algorithm, we present a simple RHP case in Fig.~\ref{fig:rhp_eg}. With parameters $N=5$ and $L=2$, the algorithm sequentially solves three subproblems: \textsf{RHP$(1,2)$}, \textsf{RHP$(2,2)$}, and \textsf{RHP$(3,2)$}. For the first two subproblems, the controls are at stages 1 and 2, respectively. In the final subproblem, the entire planning horizon---stages 3, 4, and 5---is the control, and thus it is deterministically solved.

\begin{algorithm}[htbp] 
\DontPrintSemicolon
\caption{Rolling horizon policy}
\label{alg:rhp}
\For{$s=1,\ldots,N-L-1$}{
\If{$s\ge2$}{
Update $\tilde{\mathcal{J}}_{s-1}^k$ as in Eq.~\eqref{eq:priorvertex} \;
}
Solve \RHP and store: $z_s$ [using Eq.~\eqref{eq:zs}] and $(\bm{x}_s^\ast,\bm{y}_s^\ast)$ \;
}
$s\gets N-L$ \;
Solve \RHP and store: $\{z_s,\ldots,z_{N}\}$ and $\left((\bm{x}_s^\ast,\bm{y}_s^\ast),\ldots, (\bm{x}_{N}^\ast,\bm{y}_{N}^\ast)\right)$ \;
\Return $z_{\text{rhp}} = \sum_{s\in\mathcal{S}\setminus\{0\}}z_s$ and $(\bm{x}^\ast,\bm{y}^\ast) = \left((\bm{x}_s^\ast,\bm{y}_s^\ast):s\in\mathcal{S}\setminus\{0\}\right)$ \;
\end{algorithm}

\subsection{Upper Bound of MCRP for Solution Gap Analysis} \label{sec:upper_bound}
The primary motivation for exploring the solution methods---MP and RHP---is to circumvent the computational intractability problem that arises when solving large-scale \MCRP instances. Although these methods could provide a faster means of computing feasible solutions for \MCRP compared to addressing the entire mission planning horizon at once, they do not guarantee the optimality of the solutions obtained.

To evaluate the quality of a generic algorithmic objective value $z_{\text{alg}}$, obtained from MP or RHP, in the absence of the optimal objective value $z$, we establish an upper bound (UB), $\hat{z}$, for \MCRP. This upper bound can be used to calculate the duality gap (DG), which represents the difference between $\hat{z}$ and $z_{\text{alg}}$. By ensuring that the optimal solution is always bounded between the upper bound and the algorithmic solution, we can use this information to gauge the quality of the algorithmic solution relative to the unknown optimal solution. A common approach to finding an upper bound for this purpose is to solve for the linear programming (LP) relaxation bound of \MCRP. However, obtaining the LP relaxation bound can still be computationally challenging in large-scale instances. In the following discussion, we propose a more computationally efficient upper bound metric that can be calculated using a given set of parameters.

To derive an upper bound $\hat{z}$ for \MCRP, we begin by relaxing, that is, removing, the budget constraints [\textit{i.e.}, constraints~\eqref{mcrp:res}]. As discussed in Sec.~\ref{sec:formulation}, the relaxation of the budget constraints decouples stages and therefore decomposes the problem into $N$ stage-subproblems. Moreover, this relaxation guarantees that any upper bound metric derived hereafter will always satisfy $\hat{z} \ge z$, thereby proving that the optimum $z$ is bounded between $z_{\text{alg}}$ and $\hat{z}$.

Each stage-subproblem, parameterized by $s\in\mathcal{S}\setminus\{0\}$, is equivalent to \MP without constraints~\eqref{eq:mh5}. These subproblems can be independently solved in parallel, and the results can be aggregated later to obtain $\hat{z}$. However, it is important to note that each subproblem remains a complex combinatorial optimization problem as it embeds the Maximum Coverage Problem, which is shown to be NP-hard \cite{lee2023regional,megiddo1983}. To efficiently compute an upper bound metric, we propose the following approach.

We start by examining constraints~\eqref{eq:mh4} for a given subproblem with the stage parameter $s$. Aggregating all $t$ in $\mathcal{T}_s$ and $p$ in $\mathcal{P}$, we obtain the following inequality:
\begin{equation}
    \label{eq:inequality}
    \sum_{t\in\mathcal{T}_s} \sum_{p\in\mathcal{P}} \frac{\pi_{tp}^s}{r_{tp}^s} \sum_{k\in\mathcal{K}} \sum_{i\in\mathcal{J}_{s-1}^k} \sum_{j\in\mathcal{J}_{s}^k}  V_{tjp}^{sk} x_{ij}^{sk} \ge \sum_{t\in\mathcal{T}_s} \sum_{p\in\mathcal{P}} \pi_{tp}^s y_{tp}^s
\end{equation}
In this inequality~\eqref{eq:inequality}, we observe that the right-hand side---which represents the sum of the observation rewards obtained at stage $s$ and is of interest---is bounded by the left-hand side. Computing the left-hand side as it stands would yield a valid upper bound to the right-hand side. However, this approach may be impractical due to the potentially large gap between the left-hand side and the right-hand side. We can obtain a tighter upper bound by casting it as a maximization problem as follows:
\begin{equation}
    \label{eq:ubopt}
    \hat{z}_s=\max_{x_{ij}^{sk}\in\mathcal{X}_s^k,k\in\mathcal{K}} \left\{ \sum_{t\in\mathcal{T}_s} \sum_{p\in\mathcal{P}} \frac{\pi_{tp}^s}{r_{tp}^s} \sum_{k\in\mathcal{K}} \sum_{i\in\mathcal{J}_{s-1}^k} \sum_{j\in\mathcal{J}_s^k} V_{tjp}^{sk} x_{ij}^{sk} \right\}
\end{equation}
where set $\mathcal{X}_s^k$ is defined as follows:
\begin{equation}
    \mathcal{X}_s^k = \left\{ x_{ij}^{sk} \in \{0,1\} : \sum_{i\in\mathcal{J}_{s-1}^k} \sum_{j\in\mathcal{J}_s^k} x_{ij}^{sk}=1,i\in\mathcal{J}_{s-1}^k, j\in\mathcal{J}_s^k \right\}
\end{equation}

The upper bound optimization problem in Eq.~\eqref{eq:ubopt} can be solved in parallel by decomposing it into $K$ subproblems, each with the satellite index $k=1,\ldots,K$ as a parameter:
\begin{equation}
    \hat{z}_s^k=\max_{x_{ij}^{sk}\in\mathcal{X}_s^k} \left\{ \sum_{t\in\mathcal{T}_s} \sum_{p\in\mathcal{P}} \frac{\pi_{tp}^s}{r_{tp}^s} \sum_{i\in\mathcal{J}_{s-1}^k} \sum_{j\in\mathcal{J}_s^k} V_{tjp}^{sk} x_{ij}^{sk} \right\}
\end{equation}
such that
\begin{equation}
    \label{eq:z_s_hat}
    \hat{z}_s = \sum_{k\in\mathcal{K}}\hat{z}_s^k
\end{equation}

Finally, we obtain an upper bound for \MCRP by aggregating $\hat{z}_s$ [Eq.~\eqref{eq:z_s_hat}] for all stages:
\begin{equation}
    \label{eq:ub}
    \hat{z} = \sum_{s\in\mathcal{S}\setminus\{0\}}\hat{z}_s
\end{equation}

Consequently, we have the following inequalities: $z_{\text{alg}}\le z \le \hat{z}$. Algorithm~\ref{alg:ub} outlines the process of computing the proposed upper bound metric.

\begin{algorithm}[htbp]
\DontPrintSemicolon
\caption{MCRP upper bound}
\label{alg:ub}
\For{$s=1,\ldots,N$}{
Compute and store $\hat{z}_s$ as in Eq.~\eqref{eq:ubopt} \;
}
Compute $\hat{z}$ using Eq.~\eqref{eq:ub} \;
\Return $\hat{z}$ \;
\end{algorithm}

\section{Computational Experiments} \label{sec:computational_experiments}
In Sec.~\ref{sec:comparative_analysis}, we present numerical experiments to evaluate the performance and computational efficiency of the two sequential decision-making methods outlined in Sec.~\ref{sec:solution}. We test these methods across a spectrum of mission scenarios, contrasting them with the baseline \MCRP to highlight their relative merits. Subsequently, in Sec.~\ref{sec:case_study}, we undertake a case study using real-life historical natural disaster data. This case study is designed to demonstrate the practical applicability and advantages of multi-stage constellation reconfiguration facilitated by the proposed \MCRP framework, emphasizing its relevance and effectiveness in real-world situations, particularly in tracking fast-paced dynamic events.

\subsection{Comparative Analysis} \label{sec:comparative_analysis}
\subsubsection{Design of Experiments}
We set up twenty-four unique, randomly generated test instances, categorized into two distinct sets: the first representing static targets and the second, dynamic targets. Each set contains instances varying in size from medium to large, drawing combinations of parameters from each of the following sets: $N\in\{3, 4, 5\}$, $J\in\{50, 75\}$, and $K\in\{3, 5\}$. The smallest and largest instances have at most \num{1.95e15} and \num{7.52e46} potentially feasible reconfiguration plans, respectively. Each test instance is unique, ensuring that no two instances with identical problem dimensions share the same parameters.

Both static and dynamic test sets utilize a group of $K$ homogeneous satellites in inclined circular orbits, following the Walker-delta constellation pattern rule of $\SI{80}{deg}:K/K/0$. This indicates that there are $K$ orbital planes, each inclined at 80 degrees, containing one satellite and the relative phasing between satellites in adjacent planes is zero (\textit{i.e.}, the satellites' arguments of latitude are identical). The altitude of the constellation is uniquely selected for each instance, ranging between \SI{700}{km} and \SI{2000}{km}. Additionally, both sets utilize a minimum elevation angle of $\varepsilon_{\min}=\SI{5}{deg}$ for all targets, adopted from Refs.~\cite{ulybyshev2008satellite,lee2020satellite}. Each instance has a unique set of $P=10$ spot targets randomly generated within the latitude interval $[\SI{-80}{deg},\SI{80}{deg}]$ and the longitude interval $[\SI{-180}{deg},\SI{180}{deg}]$. Varying the target point set in addition to the satellite altitude allows the comparison of the framework with respect to a spectrum of conditions, effectively not limiting the comparison to a specified set of conditions external to that of stages, slots, and number of satellites. We assume that coverage by at least one satellite is required for each target to yield rewards (\textit{i.e.}, $r^s_{tp}=1, \forall s\in\mathcal{S}\setminus\{0\}, \forall t\in\mathcal{T}_s, \forall p\in\mathcal{P}$). The visibility between targets and satellites is computed using the \texttt{access} function from the MATLAB Aerospace Toolbox~\cite{MATLAB}. We set the resource availability to $c_{\max}^k=\SI{600}{m/s}$ for all satellites. The mission planning horizon is set for \SI{5}{days} and is discretized with a time step size of $t_{\text{step}}=\SI{100}{s}$, resulting in $T=\num{4320}$. It is important to note that, without compromising the generality of the proposed \MCRP framework and the presented solution methods, the parameters adopted in this comparative analysis are chosen to illustrate to the reader the overarching impact of the proposed work on the performance and computational efficiency across a wide spectrum of mission scenarios; they are not intended to narrow the scope of the applicability of the \MCRP framework. For instance, parameters such as mission planning horizon lengths, constellation orbital specifications, minimum elevation angles, and time step sizes can be altered to model suitable remote sensing applications and to increase the fidelity of parameter generation, such as visibility profiles.

The epoch is designated as March 1, 2023, at 00:00 Coordinated Universal Time (UTC). We use the Simplified General Perturbations-4 model for the propagation of the states of satellites. The propagated positions of the satellites are then used to evaluate visibility against a set of targets.

The first set of test instances considers randomly distributed spot targets present throughout the entire mission planning horizon. We define the rewards for this set such that $\pi^s_{tp}=1, \forall s\in\mathcal{S}\setminus\{0\}, \forall t\in\mathcal{T}_s, \forall p\in\mathcal{P}$, allowing every target to have the same constant unit rewards at all times. This target set includes stationary targets relevant to remote sensing applications, such as environmental monitoring, urban development tracking, and agricultural land use analysis.

The second set of test instances considers randomly distributed spot targets with time-dependent observation rewards. We assume a set of targets partitioned based on stages $\mathcal{P}=\{\mathcal{P}_1,\ldots,\mathcal{P}_N\}$, where each target set $\mathcal{P}_s$ is associated with stage $s$ and has a non-zero observation reward only during $\mathcal{T}_s$. Given target $p\in\mathcal{P}_s$, the associated observation reward $\pi_{tp}^s$ is defined as follows:
\begin{equation*}
    \pi_{tp}^s = \begin{cases}
    1, &\text{if } t \in\mathcal{T}_s\\
    0, &\text{otherwise}
    \end{cases}
\end{equation*}
We set no observation rewards for targets outside the periods of their assigned stages, intending to simulate dynamically changing environments such that a constellation configuration optimized for one stage would be drastically unfit for another.

For the test instances, we allow satellites to perform only phasing maneuvers. Orbital slots are generated such that each within $\mathcal{J}_s^k=\{1,\ldots,J\}$ possesses identical orbital elements but distinct arguments of latitude, which are uniformly spaced between zero and 360 degrees within the orbital plane of satellite $k$. We approximate the cost $c_{ij}^{sk}$ of transferring satellite $k$ from orbital slot $i \in \mathcal{J}_{s-1}^k$ to orbital slot $j \in \mathcal{J}_s^k$ by considering these two orbital slots as the boundary conditions of a circular, coplanar phasing problem, as outlined in Ref.~\cite{vallado2013fundamentals}.

We utilize the commercial software package Gurobi Optimizer (version 11.0.0) for solving \MCRP and the subproblems of MP and RHP. All computational experiments are coded and conducted on a platform equipped with an Intel Core i9-13900 2.00 GHz (base frequency) CPU processor (24 cores and 32 threads) and 32 GB of RAM. In all cases, we allow the Gurobi Optimizer to utilize all available threads. We use the default settings of the Gurobi Optimizer, except for imposing a runtime limit of \SI{10800}{s} for \MCRP and \SI{3600}{s} for the subproblems of MP and RHP. The Gurobi Optimizer returns the incumbent best solution at the runtime limit if found thus far. For \RHP, we set $L=1$. The preliminary results indicate that with $L>1$, each \RHP subproblem becomes large-scale, making the entire RHP algorithm unsuitable for solving them in series.

We define the \textit{relative performance} (RP) metric as $(z_{\text{alg}}-z)/z_{\text{alg}}$, a signed measure used to evaluate the quality of algorithmic solutions relative to the \MCRP solution obtained by the Gurobi Optimizer; $z_{\text{alg}}$ represents the objective function value of a generic algorithmic solution (derived from MP or RHP). A positive RP indicates that a solution method outperforms the Gurobi Optimizer for \MCRP. When computing the optimal solutions for \MCRP is computationally prohibitive, we infer the quality of algorithmic solutions by computing the duality gap, which bounds the optimal solution. The duality gap is computed as $|\hat{z}-z_{\text{alg}}|/|z_{\text{alg}}|$. Although the LP relaxation solution of \MCRP can serve as $\hat{z}$, its quantification may be computationally challenging. Therefore, we report the duality gaps of the algorithmic solutions using the upper bound metric derived in Sec.~\ref{sec:upper_bound}. For \MCRP, we report the duality gaps returned by the Gurobi Optimizer, which employs the same equation but uses the dual objective bound found internally by its MILP algorithm.

\subsubsection{Numerical Results}
\begin{figure}[htbp]
     \centering
     \begin{subfigure}[b]{0.435\linewidth}
        \centering
        \includegraphics[width=\linewidth]{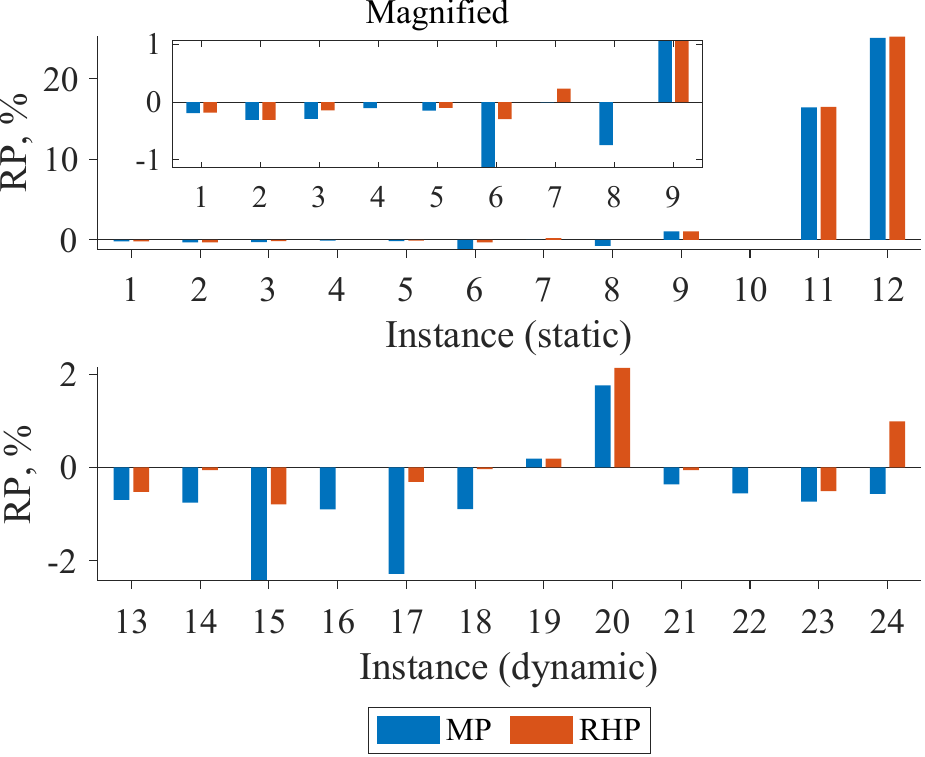}
        \caption{Relative performance versus instance.}
        \label{fig:ca:RP_instance}
     \end{subfigure}
     \hspace{1em}
     \begin{subfigure}[b]{0.435\linewidth}
        \centering
        \includegraphics[width=\linewidth]{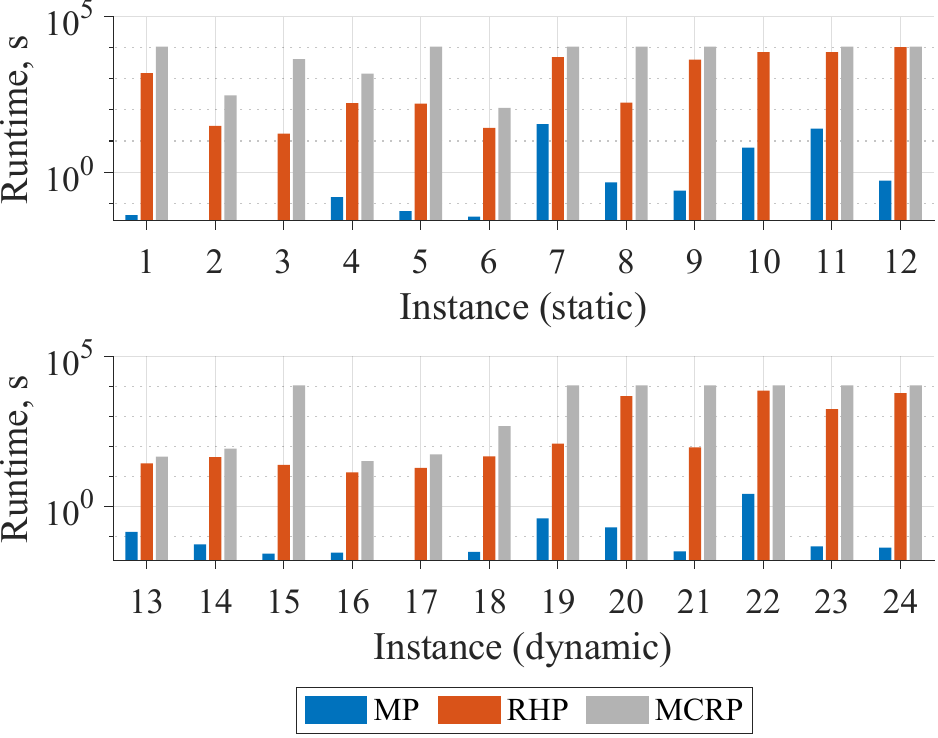}
        \caption{Runtime (log-scale) versus instance.}
        \label{fig:ca:runtime_instance}
     \end{subfigure}
    \caption{Relative performance of the algorithmic solutions with respect to MCRP and runtimes.}
    \label{fig:ca}
\end{figure}

The upper half of Table~\ref{tab:computational_results}, comprising instances 1--12, presents the results of computational experiments on the first set of test instances. Specifically, the upper halves of Figs.~\ref{fig:ca:RP_instance} and \ref{fig:ca:runtime_instance} visualize the relative performance with respect to \MCRP and computational runtime, respectively. These instances are characterized by static targets present throughout the entire mission planning horizon. Overall, when solving these instances with \MCRP, MP, and RHP, each method obtains seven, one, and seven solutions, respectively, that outperform the other methods, referred to as the best solutions (see boldface entries in the table). Out of all twelve instances, \MCRP reaches the runtime limit of \SI{10800}{s} in seven instances (indicated by hyphens), is terminated due to out of memory in one instance (indicated by Not Available, N/A), and optimally solves four instances (indicated by the zero duality gaps). MP outperforms \MCRP in four instances, while RHP does so in seven instances. The worst underperformance of MP relative to \MCRP is \SI{1.14}{\%} on instance~6. MP is computationally efficient in all twelve instances, with the maximum runtime being \SI{34.68}{s}. The most significant underperformance of RHP compared to \MCRP is \SI{0.31}{\%} on instance~2. Although RHP finds more best solutions than \MCRP and MP, it requires a more computational runtime. In five instances (\textit{i.e.}, instances 7, 9, 10, 11, and 12), RHP reaches the runtime limit of \SI{3600}{s} on some of the subproblems; these are also the instances where \MCRP fails to converge by the runtime limit of \SI{10800}{s}. However, even with RHP subproblems becoming large-scale, RHP retrieves better solutions than \MCRP in a faster overall computational runtime.

The lower part of Table~\ref{tab:computational_results}, comprising instances 13--24, reports the results of the computational experiments on the second set of test instances. Specifically, the lower halves of Figs.~\ref{fig:ca:RP_instance} and \ref{fig:ca:runtime_instance} visualize the relative performance with respect to \MCRP and computational runtime, respectively. This set is characterized by dynamic targets with time-dependent observation rewards. Overall, \MCRP, MP, and RHP obtain nine, one, and five best solutions, respectively. \MCRP optimally solves five instances, all well within the runtime limit. In instances 19, 20, and 24, \MCRP triggers the runtime limit and performs worse than the other methods. While the ability to optimize the entire mission planning horizon can be advantageous, if the scale of the instances is large, \MCRP struggles to find optimal solutions, which aligns with our original motivation to develop MP and RHP. MP outperforms \MCRP in two instances. The worst underperformance of MP relative to \MCRP is \SI{2.42}{\%}. Out of all twelve instances, MP solves all instances within \SI{2.61}{s}. MP is an attractive option that can return a high-quality solution quickly. RHP outperforms or retrieves identical solutions to \MCRP in five instances, four of which are those that \MCRP could not converge within the runtime limit, indicating their inherently large problem scale and difficulty in obtaining an optimal solution. Overall, for dynamic targets, having the ability to optimize all or at least some portion of the mission planning horizon concurrently offers better solutions. However, even with the capability to foresee the future, both \MCRP and RHP still require significant runtime when dealing with large-scale problems, which results in suboptimal solutions by the runtime limit.

Lastly, Table~\ref{tab:computational_results} also reports the baseline results without reconfiguration, $z_{\text{b}}$, and the percentage improvement of the best solutions over the baseline. In the static test set, the average improvement is approximately \SI{38.42}{\%}, ranging from a minimum of \SI{9.06}{\%} to a maximum of \SI{100.82}{\%}, which is double-fold. For the dynamic test set, the average improvement is approximately \SI{35.57}{\%}, ranging from a minimum of \SI{16.98}{\%} to a maximum of \SI{72.64}{\%}. The numerical experiments attest that having multiple opportunities for reconfiguration greatly increases the total observation rewards compared to the baseline.

It is important to note that the twenty-four test instances presented are not representative of the entire spectrum of mission scenarios. Therefore, the statistical figures of merit observed in the comparative analysis---particularly, the maximum under-/out-performance of the algorithmic solutions and the reconfiguration improvements over the baseline---should not be extrapolated to characterize the relative performances of MP and RHP compared to \MCRP in different mission scenarios. The primary objectives are (1) to validate MP and RHP for handling large-scale \MCRP instances, and (2) to illustrate the advantages of multistage reconfigurations.

In Appendix~A, we report the \MCRP runs of the test instances with a shorter runtime limit of \SI{3600}{s}. While a longer runtime limit improves the optimum and efficiency, the added benefit over the increased runtime can be limited due to the inherently large scale of the problems. Furthermore, the Gurobi Optimizer throws an out-of-memory issue due to the increased memory usage incurred during the optimization process as can be seen in instance~10. Finally, in Appendix~B, we validate the applicability and computational efficiency of the proposed upper bound metric relative to the LP relaxation bound, using the dataset from the comparative analysis.

\begin{landscape}
\begin{table}[htbp]
    \fontsize{9}{10}\selectfont
    \renewcommand{\arraystretch}{1.2}
	\caption{Results from the computational experiments. The boldfaced entries indicate the maxima of $z$, $z_{\text{mph}}$, and $z_{\text{rhp}}$.}
	\centering
	\begin{threeparttable}
	\setlength\tabcolsep{2.5pt}
	\begin{tabular}{r r r r r r r r r r r r r r r r r r}
		\hline
            \hline
		\multicolumn{4}{c}{Instance} & \multicolumn{3}{c}{\textsf{MCRP}} & UB for MP and RHP & \multicolumn{4}{c}{MP} & \multicolumn{4}{c}{RHP with $L=1$} & \multicolumn{2}{c}{Baseline} \\
		\cmidrule(lr){1-4} \cmidrule(lr){5-7} \cmidrule(lr){8-8} \cmidrule(lr){9-12} \cmidrule(lr){13-16} \cmidrule(lr){17-18}
		ID & $J$ & $N$ & $K$ & $z$ & Runtime\tnote{*}, s & DG\tnote{\textdagger}, \% & $\hat{z}$ & $z_{\text{mp}}$ & Runtime\tnote{\textdaggerdbl}, s & DG, \% & RP, \% & $z_{\text{rhp}}$ & Runtime\tnote{\textdaggerdbl}, s & DG, \% & RP, \% & $z_\text{b}$ & Improvement, \% \\
		\hline
1  & 50 & 3 & 3 & \textbf{12,066} & -        & 0.67  & 12,400 & 12,043          & 0.04  & 2.96 & -0.19                         & 12,044          & 1,524.51  & 2.96 & -0.18                         & 9,825  & 22.81  \\
2  & 50 & 4 & 3 & \textbf{8,055}  & 291.92   & 0     & 8,159  & 8,030           & 0.03  & 1.61 & -0.31                         & 8,030           & 30.45     & 1.61 & -0.31 & 6,994  & 15.17  \\
3  & 50 & 5 & 3 & \textbf{8,294}  & 4,282.03 & 0     & 8,380  & 8,270           & 0.03  & 1.33 & -0.29                         & 8,282           & 17.13     & 1.18 & -0.14                         & 7,346  & 12.90  \\
4  & 75 & 3 & 3 & \textbf{11,418} & 1,482.70 & 0     & 11,491 & 11,406          & 0.16  & 0.75 & -0.11                         & \textbf{11,418} & 164.78    & 0.64 & 0     & 7,668  & 48.90  \\
5  & 75 & 4 & 3 & \textbf{9,960}  & -        & 0.19  & 10,066 & 9,945           & 0.06  & 1.22 & -0.15                         & 9,950           & 160.75    & 1.17 & -0.10                         & 9,133  & 9.06   \\
6  & 75 & 5 & 3 & \textbf{5,746}  & 115.85   & 0     & 5,782  & 5,681           & 0.04  & 1.78 & -1.14 & 5,729           & 26.59     & 0.93 & -0.30                         & 5,202  & 10.46  \\
7  & 50 & 3 & 5 & 17,358          & -        & 2.56  & 17,918 & 17,357          & 34.68 & 3.23 & -0.01                         & \textbf{17,399} & 4,951.47  & 2.98 & 0.24  & 11,656 & 49.27  \\
8  & 50 & 4 & 5 & \textbf{10,480} & -        & 0.35  & 10,642 & 10,402          & 0.47  & 2.31 & -0.75                         & \textbf{10,480} & 171.43    & 1.55 & 0     & 7,482  & 40.07  \\
9  & 50 & 5 & 5 & 15,229          & -        & 2.76  & 15,914 & \textbf{15,393} & 0.25  & 3.38 & 1.07  & \textbf{15,393} & 4,164.89  & 3.38 & 1.07  & 10,670 & 44.26  \\
10 & 75 & 3 & 5 & N/A\tnote{\S}             & N/A      & N/A   & 17,772 & 17,314          & 6.09  & 2.65 & N/A   & \textbf{17,342} & 7,203.96  & 2.48 & N/A   & 12,569 & 37.97  \\
11 & 75 & 4 & 5 & 17,860          & -        & 23.54 & 22,254 & 21,372          & 24.96 & 4.13 & 16.43 & \textbf{21,386} & 7,247.86  & 4.06 & 16.49 & 12,632 & 69.30  \\
12 & 75 & 5 & 5 & 19,604          & -        & 35.23 & 26,672 & 26,162          & 0.53  & 1.95 & 25.07 & \textbf{26,211} & 10,575.81 & 1.76 & 25.21 & 13,052 & 100.82 \\
\hline
13 & 50 & 3 & 3 & \textbf{4,613}  & 46.03    & 0     & 4,657  & 4,581           & 0.14  & 1.66 & -0.70                         & 4,589           & 27.38     & 1.48 & -0.52                         & 3,674  & 25.56  \\
14 & 50 & 4 & 3 & \textbf{3,324}  & 84.47    & 0     & 3,364  & 3,299           & 0.05  & 1.97 & -0.76                         & 3,322           & 43.87     & 1.26 & -0.06                         & 2,541  & 30.81  \\
15 & 50 & 5 & 3 & \textbf{2,791}  & -        & 0.32  & 2,878  & 2,725           & 0.03  & 5.61 & -2.42 & 2,769           & 24.33     & 3.94 & -0.79 & 2,279  & 22.47  \\
16 & 75 & 3 & 3 & \textbf{2,356}  & 32.53    & 0     & 2,369  & 2,335           & 0.03  & 1.46 & -0.90                         & \textbf{2,356}  & 13.62     & 0.55 & 0     & 2,014  & 16.98  \\
17 & 75 & 4 & 3 & \textbf{1,925}  & 54.64    & 0     & 1,945  & 1,882           & 0.02  & 3.35 & -2.28                         & 1,919           & 19.33     & 1.35 & -0.31                         & 1,590  & 21.07  \\
18 & 75 & 5 & 3 & \textbf{2,598}  & 471.85   & 0     & 2,664  & 2,575           & 0.03  & 3.46 & -0.89                         & 2,597           & 46.14     & 2.58 & -0.04                         & 1,984  & 30.95  \\
19 & 50 & 3 & 5 & 5,369           & -        & 0.88  & 5,479  & \textbf{5,379}  & 0.40  & 1.86 & 0.19  & \textbf{5,379}  & 124.60    & 1.86 & 0.19  & 3,692  & 45.69  \\
20 & 50 & 4 & 5 & 5,138           & -        & 6.38  & 5,584  & 5,230           & 0.20  & 6.77 & 1.76  & \textbf{5,250}  & 4,673.27  & 6.36 & 2.13  & 3,508  & 49.66  \\
21 & 50 & 5 & 5 & \textbf{3,290}  & -        & 0.79  & 3,391  & 3,278           & 0.03  & 3.45 & -0.37                         & 3,288           & 93.39     & 3.13 & -0.06                         & 2,082  & 58.02  \\
22 & 75 & 3 & 5 & \textbf{5,972}  & -        & 2.03  & 6,166  & 5,939           & 2.61  & 3.82 & -0.56                         & \textbf{5,972}  & 7,203.93  & 3.25 & 0     & 4,817  & 23.98  \\
23 & 75 & 4 & 5 & \textbf{4,523}  & -        & 3.05  & 4,752  & 4,490           & 0.05  & 5.84 & -0.73                         & 4,500           & 1,732.46  & 5.60 & -0.51                         & 3,506  & 29.01  \\
24 & 75 & 5 & 5 & 4,610           & -        & 4.82  & 4,933  & 4,584           & 0.04  & 7.61 & -0.57                         & \textbf{4,656}  & 5,942.12  & 5.95 & 0.99  & 2,697  & 72.64 \\
		\hline
            \hline
	\end{tabular}
	\begin{tablenotes}
	\item[*] A hyphen (-) indicates the trigger of the runtime limit of \SI{10800}{s}.
	\item[\textdagger] The duality gap (MIPGap) is computed internally by the Gurobi Optimizer; the default optimality tolerance is \SI{0.01}{\%}.
	\item[\textdaggerdbl] Cumulative runtime of all subproblems, each imposing a runtime limit of \SI{3600}{s}.
 	\item[\S] Gurobi Optimizer runs out of memory and terminates (Gurobi error 10001: Out of memory).
	\end{tablenotes}
	\end{threeparttable}
	\label{tab:computational_results}
\end{table}
\end{landscape}

\subsection{Case Study: Tracking Hurricane Harvey} \label{sec:case_study}

In addition to the comparative analysis described in Section~\ref{sec:comparative_analysis}, a case study using historical natural disaster data is conducted. This case study serves to demonstrate the general utility of multi-stage reconfiguration, enabled by the \MCRP framework's capability to manage highly dynamic targets, thereby showcasing its practical effectiveness and relevance in real-world situations.

\begin{figure}[htbp]
    \centering
    \includegraphics[width = 0.35\linewidth]{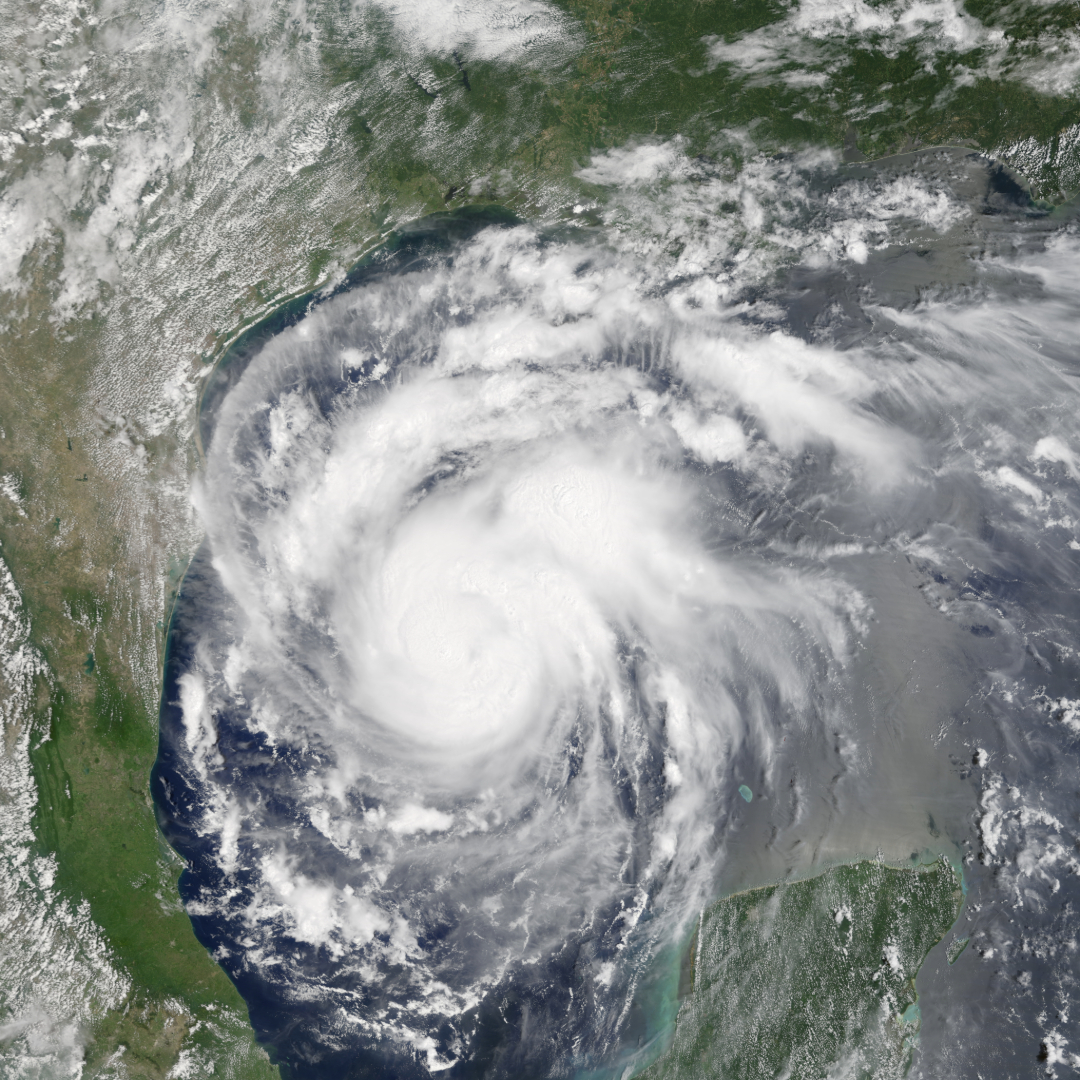}
    \caption{Hurricane Harvey, August 24, 2017 \cite{harveyimage}.}
    \label{fig:harveyimage}
\end{figure}

We selected Hurricane Harvey, a Category 4 Major Hurricane, for the case study due to its unique trajectory, extensive impact, and wealth of available data for analysis. Hurricane Harvey formed on August 16, 2017, and dissipated on September 2 of the same year. It rapidly escalated in intensity, classified as a tropical storm on August 24 and reaching Category 4 status as it made landfall on August 25 \cite{HarveyInfo}. Harvey caused widespread destruction, primarily in Texas and Louisiana, resulting in significant damages estimated at \$156.3 billion \cite{HarveyCost}. The hurricane's unprecedented rainfall led to catastrophic flooding, marking it as one of the most severe weather events in U.S. history. Figure~\ref{fig:harveyimage} shows imagery of Hurricane Harvey in the Gulf of Mexico captured by NASA satellite Terra \cite{harveyimage}.

Satellites responding to Harvey provided crucial data and imagery in order to evaluate growth and direction over time through onboard scientific instruments. This includes the NASA Global Precipitation Mission, a constellation gathering rain and snow observations, the NOAA Advanced Microwave Sounding Unit, gathering upper atmosphere temperature profiles, and the ESA Advanced Scatterometer, collecting wind speed and direction through radar measurements of electromagnetic backscatter \cite{Blake2017}. Similarly, data acquired by the Visible Infrared Imaging Radiometer Suite aboard the Suomi National Polar-Orbiting Partnership Satellite were paramount in monitoring flood water levels caused by Hurricane Harvey \cite{Li2022HarveyFlooding}. Additional instruments aboard similar sun-synchronous EO satellites that prove useful in the specific case of hurricane observation include the Temporal Experiment for Storms and Tropical Systems Technology (TEMPEST), Moderate Resolution Imaging Spectroradiometer (MODIS), and Compact Ocean Wind Vector Radiometer (COWVR). The TEMPEST systems gather ice accumulation and precipitation data within clouds \cite{TEMPEST-8517330}, the MODIS captures visual data \cite{Wilson2018MODIS}, and the COWVR gathers ocean vector winds \cite{COWVR-7943884}. Wind direction and speed data allow a path to be predicted for landfall preparation and evacuation planning by national response agencies, while rain and temperature observations provide details about development and severity over time.

\subsubsection{MCRP Modeling}

The objectives of the case study are two-fold: (1) to compare the concept of reconfiguration with a baseline case of no reconfiguration, and (2) to compare the value provided by multi-stage reconfiguration applications over the state-of-the-art single-stage reconfiguration \cite{lee2023regional}. To achieve this, multiple cases are tested by varying a single parameter of interest: the number of stages, from baseline (no reconfiguration) to single-stage (state-of-the-art), and to multiple stages (the proposed \MCRP framework). The numbers of stages considered are $N=1,2,3,4,6,$ and $8$; this configuration allows us to directly deduce the impact of multi-stage reconfiguration on the observational throughput of the system. Note that the $N=5$ and $N=7$ cases are excluded from this analysis because their $T_s = T/N$ values are not integers. In this case study, we report and compare the optimal solutions obtained by \MCRP (\textit{i.e.}, the DG of 0.01\% or less). The use of MP and RHP may sacrifice the optimality of their solutions, which is not suitable for the objectives of this case study aimed at accurately assessing the impact of the number of stages on the system's observational throughput.

\begin{figure}[htbp]
    \centering
    \includegraphics[width = 0.53\linewidth]{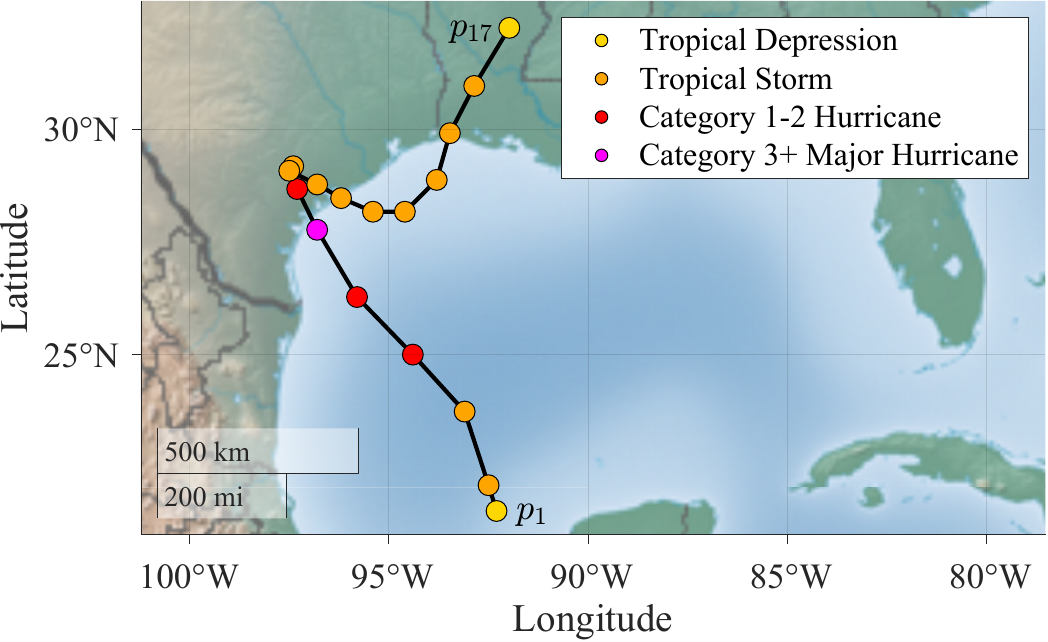}
    \caption{Sampled historical path of Hurricane Harvey; data obtained from Ref.~\cite{Blake2017}.}
    \label{fig:Harvey}
\end{figure}

We sample the approximate time when Hurricane Harvey reached its maximum intensity as a Category 4 Major Hurricane, including the periods when it was classified as a tropical depression before and after reaching this peak intensity. The duration of the selected range of Hurricane Harvey was approximately $\SI{8.5}{days}$ or $\SI{734400}{s}$ from the first to the final observation, which is discretized using $t_{\text{step}}=\SI{100}{s}$, resulting in $T=\num{7344}$. We retrieved the historical path of Hurricane Harvey from Ref.~\cite{Blake2017}. This modeling provides us with the starting point and the epoch of the simulation on August 23, 2017, at 12:00 UTC. Starting at the epoch, we sampled the position of Hurricane Harvey every twelve hours, resulting in a total of $P=17$ spot targets of interest, as shown in Table~\ref{tab:Harvey_path}. The sampled path of Hurricane Harvey is visualized in Fig.~\ref{fig:Harvey}, where the path tracks the ``best track,'' a smoothed representation of the overall hurricane location \cite{NOAAGlossary}. In the figure, $p_1$ and $p_{17}$ each represent the positions of the storm at the epoch and the last position, respectively.

\begin{table}[htbp]
    \fontsize{9}{10}\selectfont
    \renewcommand{\arraystretch}{1.2}
    \caption{Sampled historical path of Hurricane Harvey.}
    \centering
    \begin{tabular}{r r r r l r}
	\hline
        \hline
        Target & Time, UTC & Latitude, \degree N & Longitude, \degree W & Category & Reward, $R_{\text{cat}}$ \\
        \hline
        $p_1$ & 08/23/2017 12:00 & 21.40 & 92.30 & Tropical Depression & 1 \\
        $p_2$ & 08/24/2017 00:00 & 22.00 & 92.50 & Tropical Storm & 2 \\
        $p_3$ & 08/24/2017 12:00 & 23.70 & 93.10 & Tropical Storm & 2 \\
        $p_4$ & 08/25/2017 00:00 & 25.00 & 94.40 & Category 1-2 Hurricane & 3 \\
        $p_5$ & 08/25/2017 12:00 & 26.30 & 95.80 & Category 1-2 Hurricane & 3 \\
        $p_6$ & 08/26/2017 00:00 & 27.80 & 96.80 & Category 3+ Major Hurricane & 4 \\
        $p_7$ & 08/26/2017 12:00 & 28.70 & 97.30 & Category 1-2 Hurricane & 3 \\
        $p_8$ & 08/27/2017 00:00 & 29.20 & 97.40 & Tropical Storm & 2 \\
        $p_9$ & 08/27/2017 12:00 & 29.10 & 97.50 & Tropical Storm & 2 \\
        $p_{10}$ & 08/28/2017 00:00 & 28.80 & 96.80 & Tropical Storm & 2 \\
        $p_{11}$ & 08/28/2017 12:00 & 28.50 & 96.20 & Tropical Storm & 2 \\
        $p_{12}$ & 08/29/2017 00:00 & 28.20 & 95.40 & Tropical Storm & 2 \\
        $p_{13}$ & 08/29/2017 12:00 & 28.20 & 94.60 & Tropical Storm & 2 \\
        $p_{14}$ & 08/30/2017 00:00 & 28.90 & 93.80 & Tropical Storm & 2 \\
        $p_{15}$ & 08/30/2017 12:00 & 29.92 & 93.47 & Tropical Storm & 2 \\
        $p_{16}$ & 08/31/2017 00:00 & 30.94 & 92.86 & Tropical Storm & 2 \\
        $p_{17}$ & 08/31/2017 12:00 & 32.18 & 91.98 & Tropical Depression & 1 \\
        \hline
        \hline
    \end{tabular}
    \label{tab:Harvey_path}
\end{table}

By sampling the path of Hurricane Harvey, we assume that each sampled spot point remains stationary at its position for a duration of twelve hours. In addition, we set the minimum elevation angle to \SI{10}{deg} for all spot targets. Since these spot points appear sequentially over time, we vary the observation rewards to model Hurricane Harvey's motion within the \MCRP framework. Additionally, we assign more weight to observing targets when the intensity of the tropical cyclone system, measured based on sustained surface wind speed, is greater. There are four different categories of the system for the sampled period: tropical depression (38 mph or less), tropical storm (39-73 mph), hurricane (74 mph or greater), and major hurricane (111 mph or greater) \cite{NOAAGlossary}. Each category is assigned a numerical value $R_{\text{cat}}$ of one, two, three, and four as summarized in Table~\ref{tab:Harvey_path}, respectively, to model their varying importance. Assigning a higher weighted reward value during periods of increased storm intensity motivates the reconfiguration process to prioritize observations in these intervals. This approach potentially provides valuable information to first responders for disaster management and helps predict its future path. Mathematically, we impose the following observation rewards on the sampled spot targets:
\begin{equation*}
    \pi_{tp}=\begin{cases}
        R_{\text{cat}}, &\text{if } t\in\mathcal{T}_p \\
        0, &\text{otherwise}
    \end{cases}
\end{equation*}
where $\mathcal{T}_p := \left\{t:t_p\le t \le t_{p+1},t\in\mathcal{T}\right\}$ represents the time interval at which spot target $p$ has non-zero observation rewards, and $t_p := (p-1)T/P+1$. For example, $p_1$ has its observation reward equal to one in the interval $t\in[1, T/P]$ and zero elsewhere, while $p_{10}$ has its observation reward equal to two during $t\in[9T/P+1, 10T/P]$, and zero elsewhere before and after.

The following parameters are chosen for the case study. We consider a constellation of four satellites in circular orbits, with their states defined at the epoch as shown in Table~\ref{tab:sat_specs}. The satellites' states, including the orbital parameters and $c_{\max}^k$, are randomly generated with altitudes ranging from \SI{700}{km} to \SI{1200}{km}, inclinations between \SI{60}{deg} and \SI{110}{deg}, right ascension of the ascending node (RAAN), and argument of latitude from \SI{0}{deg} to \SI{360}{deg}, and $c_{\max}^k$ values from \SI{0.8}{km/s} to \SI{1.2}{km/s}. Unlike the homogeneous constellation configurations from the test instances of the comparative analysis in Sec.~\ref{sec:comparative_analysis}, the configuration here represents a federated system of heterogeneous satellites, each placed in a unique orbit and with a unique budget.

\begin{table}[htbp]
    \fontsize{9}{10}\selectfont
    \renewcommand{\arraystretch}{1.2}
    \caption{Key satellite specification parameters defined at the epoch.}
    \centering
    \begin{tabular}{r r r r r r}
    \hline
    \hline
    Satellite, $k$ & Altitude, km & Inclination, deg. & RAAN, deg. & Arg. of Lat., deg. & $c_{\max}^k$, km/s  \\
    \hline
    1 & 1,135.06 & 80.56 & 200.24 & 160.93 & 0.930 \\
    2 & 991.14 & 65.87 & 132.15 & 210.76 & 1.079 \\
    3 & 839.42 & 94.25 & 144.85 & 58.31 & 0.947 \\
    4 & 792.96 & 81.88 & 40.69 & 187.46 & 1.135 \\
    \hline
    \hline
    \end{tabular}
    \label{tab:sat_specs}
\end{table}

In this case study, we define the orbital slots to allow for changes in inclination, RAAN, and argument of latitude, with the capability to change any combination of options from one stage to the next. Each orbital slot in inclination and RAAN is evenly distributed in the positive and negative directions of the initial slot such that there are eight additional inclination and RAAN slots, with a shared option at the initial slot, resulting in 17 unique plane slots along an axis of inclination and RAAN. Similarly, each argument of latitude option is evenly distributed among the interval $[\SI{0}{deg}, \SI{360}{deg})$ such that there are 24 unique phase slots. As a result of the plane and phase slot combinations, the total number of slots is $J=408$. The equations regarding the degree of separation between plane slots and an example of the cost computation between two sets of orbital elements are given in Appendix~C.

All other parameters utilized in the case study are identical to those in the two sets of test instances shown in Sec.~\ref{sec:comparative_analysis}. 

\subsubsection{Results Analysis}
Using the previously described parameters, we run a total of seven different cases: the baseline, single-stage, and multi-stage with $N=2,3,4,6,$ and $8$. In this section, we analyze the obtained results.

\begin{figure}[htbp]
    \centering
    \includegraphics[width=0.5\textwidth]{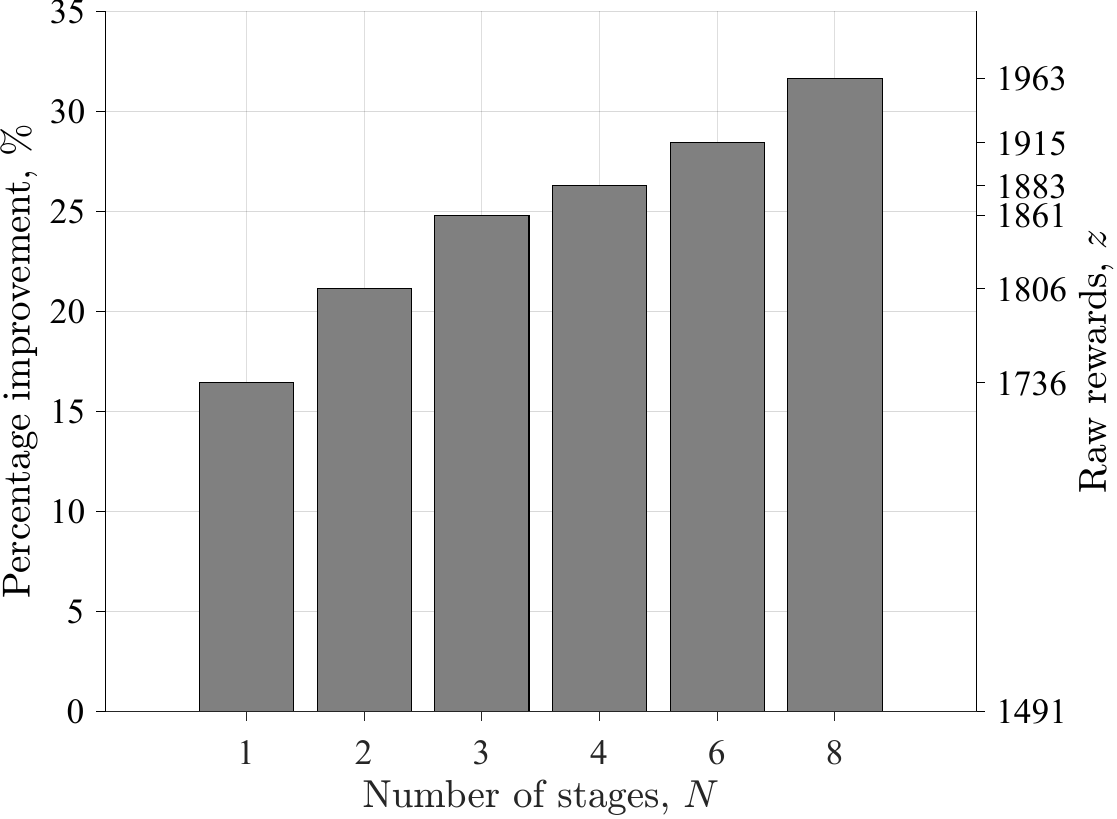}
    \caption{Percentage improvement over the baseline, reported per number of stages.}
    \label{fig:HarveyResults}
\end{figure}

Figure~\ref{fig:HarveyResults} depicts the percentage improvement over the baseline scenario across various reconfiguration cases, shown on the left y-axis, alongside the raw rewards $z$ on the right y-axis. The baseline scenario, without any reconfiguration, achieved a raw reward of $z=1491$. In contrast, the single-stage reconfiguration yielded a higher raw reward of $z=1736$, marking a \SI{16.43}{\%} improvement. The results also demonstrate that increasing the number of stages in the reconfiguration process enhances overall performance. This is exemplified by the eight-stage case, which achieved the highest improvement, with a reward of $z=1963$, an increase of \SI{31.66}{\%} over the baseline. These findings underscore the significant benefits of reconfiguration, evident even in the single-stage case, and further suggest that more stages lead to greater improvements in observational throughput.

\begin{figure}[htbp]
    \centering
    \includegraphics[width=0.48\textwidth]{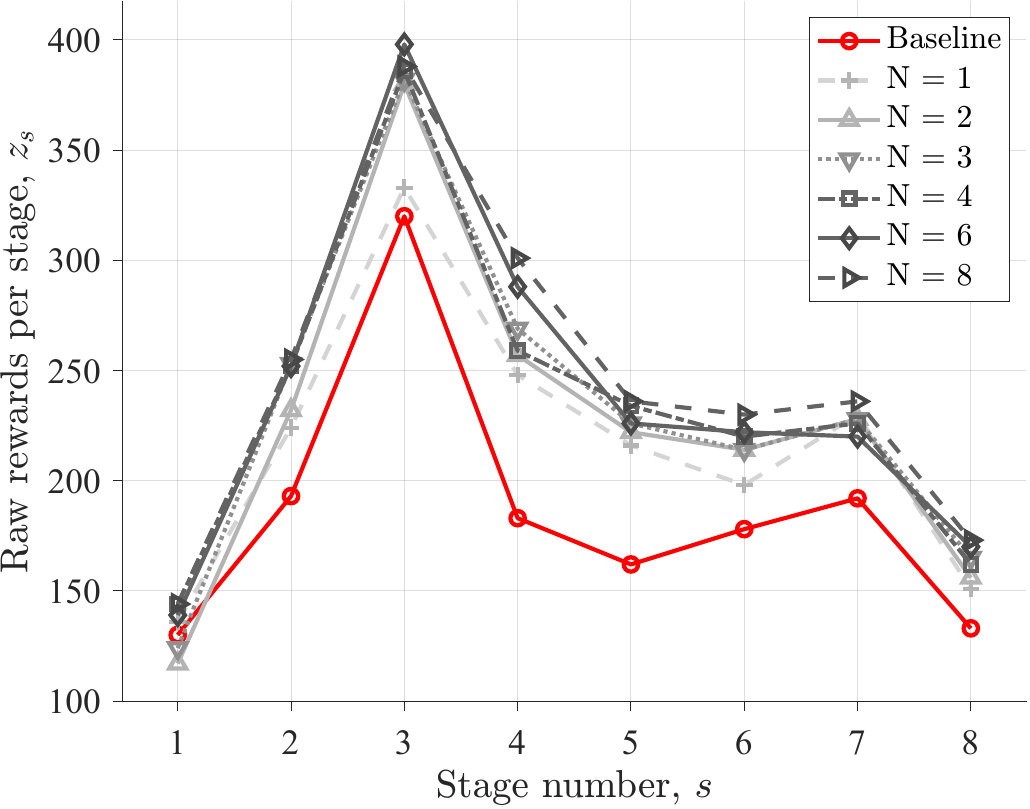}
    \caption{Raw rewards over time per case.}
    \label{fig:VTWRewards}
\end{figure}

Figure~\ref{fig:VTWRewards} presents the raw reward, $z$, for each discretized time interval relative to the stage time interval $\mathcal{T}_s$ in the eight-stage scenario. Here, rewards are discretized by time intervals rather than by stage. Doing so enables us to perform a direct and detailed comparison of rewards across identical time intervals for each case. As a result, the baseline case, despite having no stages, is represented over eight discrete time intervals aligned with the eight-stage scenario; this is similarly applied to reconfiguration cases that have fewer than eight stages. The available rewards for each time interval resulting from Table~\ref{tab:Harvey_path} are \num{1404}, \num{2376}, \num{3186}, \num{2106}, \num{1836}, \num{1836}, \num{1836}, and \num{1404} for the first through eighth time interval, respectively. The figure demonstrates the benefits of increasing the number of stages, showing a marked overall improvement over the baseline case, especially notable at the conclusion of the third time interval. The third time interval displays a large spike in rewards in comparison to the other time intervals as a result of the associated rewards reported in Table~\ref{tab:Harvey_path}. The third time interval contains the largest attainable score of all eight time intervals, and all six cases dedicate the most $\Delta v$ consumption to maneuvers performed prior to the time interval conclusion. For the third time interval, the amount of $\Delta v$ consumed by satellites is \SI{3.05}{km/s}, \SI{2.89}{km/s}, \SI{2.22}{km/s}, \SI{3.29}{km/s}, \SI{2.94}{km/s}, and \SI{2.48}{km/s} for each case. This amounts to \SI{74.49}{\%}, \SI{70.76}{\%}, \SI{54.22}{\%}, \SI{80.31}{\%}, \SI{71.82}{\%}, and \SI{60.73}{\%} of the total budget for each case, demonstrating a significant commitment in $\Delta v$ expenditure to secure the higher-valued rewards.

Next in the results analysis, we select the eight-stage case for an in-depth analysis due to the high level of rewards gathered and the high level of complexity allowing demonstration of the \MCRP formulations capabilities.

In the analysis of the eight-stage case, key information can be gained through analysis of the type of transfer performed at each stage. Table~\ref{tab:HarveySlotProgression} depicts the flow of orbit transfer from one stage to the next, the corresponding orbital elements of the destination slot upon arrival at $t_s$, and the type of transfer performed. The majority of plane changes occur in the early stages, prior to the most reward-populated stage, with only one plane change occurring later on, demonstrating that the constellation is aligning itself with the most optimal plane when the majority of the rewards are available. Emphasis on phasing maneuvers extends to 23 out of the 32 transfers being selected as phasing only, and all maneuvers performed having phasing involved, suggesting a large importance of the location in the argument of latitude. This is well illustrated by the configurations of the satellite constellation throughout the mission horizon, as depicted in Appendix~D.

\begin{table}[htb]
    \fontsize{9}{10}\selectfont
    \renewcommand{\arraystretch}{1.2}
    \caption{Flow of optimal slot transfer, $N=8$.}
    \centering
    \resizebox{\textwidth}{!}{
    \begin{tabular}{r r r r r r r r}
    \hline
    \hline
    \multirow{2}{*}{Satellite, $k$} & \multirow{2}{*}{Stage, $s$} & \multirow{2}{*}{Consumed $\Delta v$, km/s} & \multirow{2}{*}{Transfer type} & \multicolumn{4}{c}{Destination slot}\\
		\cmidrule(lr){5-8}
  & & & & Index, $j$ & Inclination, deg. & RAAN, deg. & Arg. of Lat. at $t_s$, deg. \\
\hline
    1 & 1 & 0.568 & RAAN/Phase      & 260 & 80.55 & 195.79 & 85.93 \\
      & 2 & 0.043 & Phase           & 255 & 80.55 & 195.79 & 70.07 \\
      & 3 & 0.085 & Phase           & 249 & 80.55 & 195.79 & 39.20 \\
      & 4 & 0.043 & Phase           & 244 & 80.55 & 195.79 & 23.34 \\
      & 5 & 0.043 & Phase           & 263 & 80.55 & 195.79 & 7.47  \\
      & 6 & 0.002 & Phase           & 259 & 80.55 & 195.79 & 6.61  \\
      & 7 & 0.085 & Phase           & 253 & 80.55 & 195.79 & 335.75\\
      & 8 & 0     & None            & 253 & 80.55 & 195.79 & 34.88 \\
\hline
    2 & 1 & 0.650 & Inc/Phase       & 47  & 60.82 & 132.15 & 180.76\\
      & 2 & 0     & None            & 47  & 60.82 & 132.15 & 30.00 \\
      & 3 & 0.255 & Inc/Phase       & 20  & 59.13 & 132.15 & 194.25\\
      & 4 & 0.118 & Phase           & 9   & 59.13 & 132.15 & 238.49\\
      & 5 & 0.043 & Phase           & 18  & 59.13 & 132.15 & 222.74\\
      & 6 & 0     & None            & 18  & 59.13 & 132.15 & 71.98 \\
      & 7 & 0     & None            & 18  & 59.13 & 132.15 & 281.22\\
      & 8 & 0.002 & Phase           & 4   & 59.13 & 132.15 & 280.47\\
\hline
    3 & 1 & 0.226 & RAAN/Phase      & 292 & 94.25 & 143.39 & 103.31\\
      & 2 & 0     & None            & 292 & 94.25 & 143.39 & 118.94\\
      & 3 & 0.422 & RAAN/Phase      & 264 & 94.25 & 140.46 & 74.57 \\
      & 4 & 0.125 & Phase           & 260 & 94.25 & 140.46 & 30.19 \\
      & 5 & 0.125 & Phase           & 256 & 94.25 & 140.46 & 345.82\\
      & 6 & 0     & None            & 256 & 94.25 & 140.46 & 1.44  \\
      & 7 & 0.040 & Phase           & 254 & 94.25 & 140.46 & 347.07\\
      & 8 & 0.002 & Phase           & 253 & 94.25 & 140.46 & 347.69\\
\hline
    4 & 1 & 0.199 & Phase           & 318 & 81.88 & 40.69  & 262.46\\
      & 2 & 0.018 & Phase           & 313 & 81.88 & 40.69  & 255.80\\
      & 3 & 0.018 & Phase           & 332 & 81.88 & 40.69  & 249.15\\
      & 4 & 0.103 & Phase           & 325 & 81.88 & 40.69  & 212.49\\
      & 5 & 0.023 & Phase           & 105 & 81.88 & 40.69  & 220.84\\
      & 6 & 0.060 & Phase           & 99  & 81.88 & 40.69  & 199.18\\
      & 7 & 0.381 & RAAN/Phase      & 355 & 81.88 & 42.46  & 147.53\\
      & 8 & 0.060 & Phase           & 349 & 81.88 & 42.46  & 125.88\\
\hline
\hline
    \end{tabular}
    }
    \label{tab:HarveySlotProgression}
\end{table}

Overall, the Hurricane Harvey case study highlights the significance of constellation reconfiguration in general, demonstrating a performance level higher than the baseline. Our findings reveal that all reconfiguration cases, including the state-of-the-art single-stage reconfiguration \cite{lee2023regional}, outperform the baseline by providing highly flexible orbital slots as options. Multi-stage reconfiguration, facilitated by the proposed \MCRP framework, has empirically demonstrated further enhancement in performance by providing additional reconfiguration opportunities, even under identical parameters such as budget.

\section{Conclusions} \label{sec:conclusions}
This paper addresses the MCRP, a problem of optimizing the sequence of orbital maneuvers for satellites to maximize the total observation rewards obtained by covering a set of targets of interest. To model stage transitions and fuel consumption by satellites, we adopt the concept of TEGs. This involves expanding the vertices (the orbital slots) forward in time and constructing directed edges (the orbital transfers) between the vertices of any two adjacent stages. Based on this model, we propose a novel ILP formulation of the MCRP, which allows for the use of commercial MILP solvers, facilitating convenient handling and provably optimal solutions.

To address the issue of computational intractability in solving large-scale \MCRP instances, we propose two sequential decision-making approaches: MP and RHP. Through computational experiments, we empirically demonstrate that \MCRP performs well for small-scale instances. However, for large-scale instances, whether the mission scenarios involve static or dynamic targets, both MP and RHP provide high-quality solutions, outperforming \MCRP in numerous test instances. Generally, RHP outperforms MP by making informed decisions that exploit the deterministic forecast of the impact of current-stage decisions on subsequent stages. Nevertheless, as observed in several large-scale test instances, RHP, similar to \MCRP, is subject to an exponential expansion of the solution space, even with a single-stage lookahead policy. This challenge renders MP a high-performing and computationally efficient alternative.

Both the comparative analysis and the case study of tracking Hurricane Harvey confirm the findings in the existing literature regarding the general utility of constellation reconfiguration in EO applications as a flexible, system-level operational approach, effective regardless of the number of stages employed. Furthermore, our experiment, employing our proposed \MCRP framework, uniquely demonstrates that increased opportunities for reconfiguration (\textit{i.e.}, more stages), even when constrained by the same budget, lead to a higher total observational throughput than that of single-stage reconfiguration, which is considered the state of the art in the literature.

Some parameters and parameter generation methods used in the computational experiments are subject to improvement for a more rigorous assessment of the values provided by the presented \MCRP framework across a wider spectrum of mission scenarios. A more diverse set of parameters, such as mission planning lengths, constellation orbital elements, and the number and spread of target points, could be explored to better characterize the sensitivity of the model in representing real-world remote sensing applications. Additionally, the components of the framework and the algorithms used to generate parameters, such as the reward mechanism, visibility computation, and trajectory optimization techniques, could be enhanced to explore the full potential of the proposed \MCRP framework and enable higher-fidelity results.

There are several fruitful directions for future research. The first is to enhance the applicability of the proposed \MCRP framework to real-world satellite systems operations. This paper considers a case where orbital transfers are the sole decision variables. To accurately assess the impact of the concept of multi-stage reconfiguration in Earth observation systems, the proposed problem, \MCRP, should be integrated with a scheduler framework known as the Earth observation satellite scheduling problem, which aims to optimize the schedules of satellites' actions under various operational constraints. Such a scheduling framework may include tasks of observation and data downlink and/or constraints of onboard data and battery capacity. The second direction is to address computational intractability. Various algorithmic approaches and simplifications can be explored, such as relaxing constraints~\eqref{mcrp:cov} for non-cooperative satellite observational maneuver planning or using approximate dynamic programming, to reduce the time complexity for \MCRP and RHP. Lastly, in this paper, MCRP is explored in a deterministic setting. Therefore, an interesting follow-up research topic would be the investigation of a stochastic variant, incorporating the trajectory of a dynamic target under some probability distribution.

\appendix

\section*{Appendix A: MCRP with a Shorter Runtime Limit}
We report \MCRP runs of the twenty-four test instances conducted in Sec.~\ref{sec:comparative_analysis} with a runtime limit of \SI{3600}{s} and juxtapose the results with a runtime limit of \SI{10800}{s} in Table~\ref{tab:1hr_vs_3hr}.

\begin{table}[htbp]
    \fontsize{9}{10}\selectfont
    \renewcommand{\arraystretch}{1.2}
	\caption{Comparison of MCRP results with the one-hour and three-hour runtime limits.}
	\centering
        \begin{threeparttable}
	\setlength\tabcolsep{2.5pt}
	\begin{tabular}{r r r r r r r r}
		\hline
            \hline
		\multicolumn{4}{c}{Instance} & \multicolumn{2}{c}{\textsf{MCRP (1-hr limit)}} & \multicolumn{2}{c}{\textsf{MCRP (3-hr limit)}} \\
		\cmidrule(lr){1-4} \cmidrule(lr){5-6} \cmidrule(lr){7-8}
		ID & $J$ & $N$ & $K$ & $z$ & DG\tnote{\textdagger}, \% & $z$ & DG\tnote{\textdagger}, \% \\
		\hline
1 & 50 & 3 & 3 & 12,055 & 1.09 & 12,066 & 0.67 \\
2 & 50 & 4 & 3 & 8,055 & 0 & 8,055 & 0 \\
3 & 50 & 5 & 3 & 8,294 & 0.12 & 8,294 & 0 \\
4 & 75 & 3 & 3 & 11,418 & 0 & 11,418 & 0 \\
5 & 75 & 4 & 3 & 9,959 & 0.38 & 9,960 & 0.19 \\
6 & 75 & 5 & 3 & 5,746 & 0 & 5,746 & 0 \\
7 & 50 & 3 & 5 & 15,535 & 14.62 & 17,358 & 2.56 \\
8 & 50 & 4 & 5 & 10,480 & 0.60 & 10,480 & 0.35 \\
9 & 50 & 5 & 5 & 14,506 & 7.89 & 15,229 & 2.76 \\
10 & 75 & 3 & 5 & 16,039 & 9.81 & N/A\tnote{\S} & N/A\tnote{\S} \\
11 & 75 & 4 & 5 & 17,716 & 24.57 & 17,860 & 23.54 \\
12 & 75 & 5 & 5 & 18,484 & 43.42 & 19,604 & 35.23 \\
\hline
13 & 50 & 3 & 3 & 4,613 & 0 & 4,613 & 0 \\
14 & 50 & 4 & 3 & 3,324 & 0 & 3,324 & 0 \\
15 & 50 & 5 & 3 & 2,791 & 0.68 & 2,791 & 0.32 \\
16 & 75 & 3 & 3 & 2,356 & 0 & 2,356 & 0 \\
17 & 75 & 4 & 3 & 1,925 & 0 & 1,925 & 0 \\
18 & 75 & 5 & 3 & 2,598 & 0 & 2,598 & 0 \\
19 & 50 & 3 & 5 & 5,369 & 0.91 & 5,369 & 0.88 \\
20 & 50 & 4 & 5 & 5,091 & 7.39 & 5,138 & 6.38 \\
21 & 50 & 5 & 5 & 3,284 & 1.10 & 3,290 & 0.79 \\
22 & 75 & 3 & 5 & 5,969 & 2.08 & 5,972 & 2.03 \\
23 & 75 & 4 & 5 & 4,504 & 3.49 & 4,523 & 3.05 \\
24 & 75 & 5 & 5 & 4,531 & 6.62 & 4,610 & 4.82 \\
		\hline
            \hline
	\end{tabular}
 \begin{tablenotes}
	\item[\textdagger] The duality gap (MIPGap) is computed internally by the Gurobi Optimizer; the default optimality tolerance is \SI{0.01}{\%}.
        \item[\S] Gurobi Optimizer runs out of memory and terminates (Gurobi error 10001: Out of memory).
	\end{tablenotes}
	\end{threeparttable}
	\label{tab:1hr_vs_3hr}
\end{table}

\section*{Appendix B: Comparison of Upper Bound Metrics}
In this appendix, we compare the performance and computational efficiency of the upper bound metric proposed in Sec.~\ref{sec:upper_bound} with the LP relaxation method, using the dataset from the comparative analysis conducted in Sec.~\ref{sec:comparative_analysis}. Although the proposed upper bound metric can be straightforwardly computed using a maximum operator given a set of feasible sets, we formulated the upper bound metric as an optimization problem and solved it using the Gurobi Optimizer to encourage a fair comparison with the LP relaxation.

The results are reported in Table~\ref{tab:upper_bound_results}. In all test instances, the proposed upper bound metric $\hat{z}$ was able to find a value almost immediately. The metric simply involves finding a maximum value given a feasible set; the Gurobi Optimizer easily found the solution during the presolve stage. In all instances, the LP relaxation found tighter upper bounds, $\hat{z}_\text{LP}$, than those proposed. However, as expected, we can observe that for certain large-scale instances, particularly those from the static target set, the LP remains difficult to solve, with the maximum runtime of \SI{1807.67}{s}. Here, we define the relative performance metric as follows: $(\hat{z}-\hat{z}_\text{LP})/\hat{z}$. The proposed upper bound metric is able to find values with at most \SI{2.01}{\%} difference from the LP relaxation upper bound. Both the proposed upper bound metric and the LP relaxation bound can be used as proper upper bounds, but the LP relaxation comes with additional computational runtime.

\begin{table}[htb]
    \fontsize{9}{10}\selectfont
    \renewcommand{\arraystretch}{1.2}
	\caption{Comparison of the upper bound metrics for MCRP.}
	\centering
	\begin{threeparttable}
	\setlength\tabcolsep{2.5pt}
	\begin{tabular}{r r r r r r r r r}
		\hline
            \hline
		\multicolumn{4}{c}{Instance} & \multicolumn{3}{c}{Upper bound} & \multicolumn{2}{c}{LP relaxation}  \\
		\cmidrule(lr){1-4} \cmidrule(lr){5-7} \cmidrule(lr){8-9}
		ID & $J$ & $N$ & $K$ & $\hat{z}$ & Runtime\tnote{*}, s & RP, \% & $\hat{z}_\text{LP}$ & Runtime, s \\
		\hline
1  & 50 & 3 & 3 & 12,400 & \textless{}0.01 & -1.42 & 12,223.71 & 228.02   \\
2  & 50 & 4 & 3 & 8,159  & \textless{}0.01 & -0.98 & 8,079.33  & 39.44    \\
3  & 50 & 5 & 3 & 8,380  & \textless{}0.01 & -0.38 & 8,348.04  & 35.98    \\
4  & 75 & 3 & 3 & 11,491 & \textless{}0.01 & -0.30 & 11,456.40 & 362.70   \\
5  & 75 & 4 & 3 & 10,066 & \textless{}0.01 & -0.55 & 10,010.72 & 655.75   \\
6  & 75 & 5 & 3 & 5,782  & \textless{}0.01 & -0.45 & 5,756.03  & 43.63    \\
7  & 50 & 3 & 5 & 17,918 & 0.02            & -0.53 & 17,822.87 & 1,485.73 \\
8  & 50 & 4 & 5 & 10,642 & \textless{}0.01 & -0.70 & 10,567.30 & 213.57   \\
9  & 50 & 5 & 5 & 15,914 & 0.02            & -1.53 & 15,670.86 & 908.93   \\
10 & 75 & 3 & 5 & 17,772 & 0.02            & -0.83 & 17,625.27 & 1,807.67 \\
11 & 75 & 4 & 5 & 22,254 & \textless{}0.01 & -0.83 & 22,069.43 & 1,541.30 \\
12 & 75 & 5 & 5 & 26,672 & \textless{}0.01 & -0.60 & 26,511.29 & 1,743.13 \\
\hline
13 & 50 & 3 & 3 & 4,657  & \textless{}0.01 & -0.85 & 4,617.42  & 41.84    \\
14 & 50 & 4 & 3 & 3,364  & \textless{}0.01 & -0.90 & 3,333.67  & 13.75    \\
15 & 50 & 5 & 3 & 2,878  & \textless{}0.01 & -1.81 & 2,825.78  & 7.85     \\
16 & 75 & 3 & 3 & 2,369  & \textless{}0.01 & -0.40 & 2,359.59  & 11.45    \\
17 & 75 & 4 & 3 & 1,945  & 0.02            & -0.81 & 1,929.26  & 13.97    \\
18 & 75 & 5 & 3 & 2,664  & \textless{}0.01 & -1.98 & 2,611.22  & 20.52    \\
19 & 50 & 3 & 5 & 5,479  & \textless{}0.01 & -0.80 & 5,435.22  & 63.72    \\
20 & 50 & 4 & 5 & 5,584  & \textless{}0.01 & -2.00 & 5,472.06  & 38.49    \\
21 & 50 & 5 & 5 & 3,391  & \textless{}0.01 & -1.84 & 3,328.52  & 11.10    \\
22 & 75 & 3 & 5 & 6,166  & \textless{}0.01 & -0.99 & 6,104.88  & 147.84   \\
23 & 75 & 4 & 5 & 4,752  & \textless{}0.01 & -1.77 & 4,667.91  & 51.85    \\
24 & 75 & 5 & 5 & 4,933  & \textless{}0.01 & -2.01 & 4,833.72  & 45.36   \\
		\hline
            \hline
	\end{tabular}
	\begin{tablenotes}
	\item[*] \textless{}0.01 indicates the runtime of less than \SI{0.01}{s}.
	\end{tablenotes}
	\end{threeparttable}
	\label{tab:upper_bound_results}
\end{table}

\section*{Appendix C: Orbital Slot Bounds and Example Cost Computation}

We utilize orbital boundary value problems found in Chapter Six of Ref.~\cite{vallado2013fundamentals} to compute the bounds of plane change slots. These bounds are computed as if the entire budget, $c^k_{\max}$, is utilized in a single maneuver and reduced by a scaling factor, $\eta$.  The equation used to calculate the inclination bounds, the maximum difference in inclination $\Delta i$ from the initial condition, is shown in Eq.~\eqref{eq:delta_i}:
\begin{equation}
\delta = 2 \eta \arcsin\left( \frac{c^k_{\max}}{2\sqrt{(\mu/a_k)}} \right) \label{eq:delta_i}
\end{equation}
where we let $\Delta i=\delta$ for the inclination change only. Similarly, the equation used to calculate the RAAN bounds, the maximum difference in RAAN $\Delta \Omega$ from the initial condition, is shown in Eq.~\eqref{eq:delta_RAAN}:
\begin{equation}
    \Delta \Omega = \eta \arccos\left( \frac{\cos \delta - \cos^2 i_k}{\sin^2 i_k} \right) \label{eq:delta_RAAN}
\end{equation}
where $\delta$ is computed from Eq.~\eqref{eq:delta_i}. In Eqs.~\eqref{eq:delta_i} and \eqref{eq:delta_RAAN}, $a_k$ is the semi-major axis of satellite $k$, $i_k$ is the inclination of satellite $k$, and $\mu$ is the gravitational parameter of Earth. In this paper, the scaling factor is set to \num{0.8}.

Many potential transfer costs are possible given the diversity of the orbital slots available to each satellite. These include changes in inclination, RAAN, simultaneous changes in inclination and RAAN, and changes in the argument of latitude. If both inclination and RAAN differ between slots, we change them simultaneously to minimize the cost. If both plane and argument of latitude changes are required, we first perform the plane change. Additionally, phasing considers a maximum of five revolutions of the target and transfer orbit regarding the circular, coplanar phasing problem in Chapter Six of Ref.~\cite{vallado2013fundamentals}. As a result of this diversity, some examples are prudent to provide information for reference to the costs provided as the result of the reconfiguration process. For the purpose of these examples, a satellite on a circular orbit at an altitude of \SI{1000}{km} and with a \SI{45}{deg} inclination, RAAN, and initial argument of latitude $u$ will be used. Some prominent examples are provided in Table~\ref{tab: example costs}, including various combinations of possible changes and directions for the change.

\begin{table}[htbp]
    \centering
    \caption{Transfer costs for various cases.}
    \begin{tabular}{l l r}
        \hline
        \hline
        Transfer type & Transfer case & Cost, km/s \\
        \hline
        Inclination-only & $\Delta i = \pm \SI{5}{deg}$                             & $0.641$ \\ 
        RAAN-only & $\Delta \Omega = \pm \SI{5}{deg}$                         & $0.453$ \\
        Phase & $\Delta u = \pm \SI{5}{deg}$   & $0.014$ \\
        Phase & $\Delta u = +\SI{120}{deg}$      & $0.306$ \\
        Phase & $\Delta u = -\SI{120}{deg}$      & $0.350$ \\
        Combined Inc/RAAN & $\Delta i = \Delta \Omega = +\SI{5}{deg}$                  & $0.796$ \\
        Combined Inc/RAAN & $\Delta i = \Delta \Omega = -\SI{5}{deg}$                  & $0.773$ \\
        Combined Inc/RAAN/Phase & $\Delta i = \Delta \Omega = +\SI{5}{deg},~\Delta u = +\SI{120}{deg}$ & $1.127$ \\
        Combined Inc/RAAN/Phase & $\Delta i = \Delta \Omega = -\SI{5}{deg},~\Delta u = -\SI{120}{deg}$ & $1.129$ \\
        \hline
        \hline
    \end{tabular}
    \label{tab: example costs}
\end{table}

\section*{Appendix D: Constellation Formation Over Eight Stages in Response to Hurricane Harvey}
\begin{figure}[!h]
            \centering
            \renewcommand{\arraystretch}{0.8}
            \begin{tabular}{cccc}
                        \includegraphics[height=3.5cm]{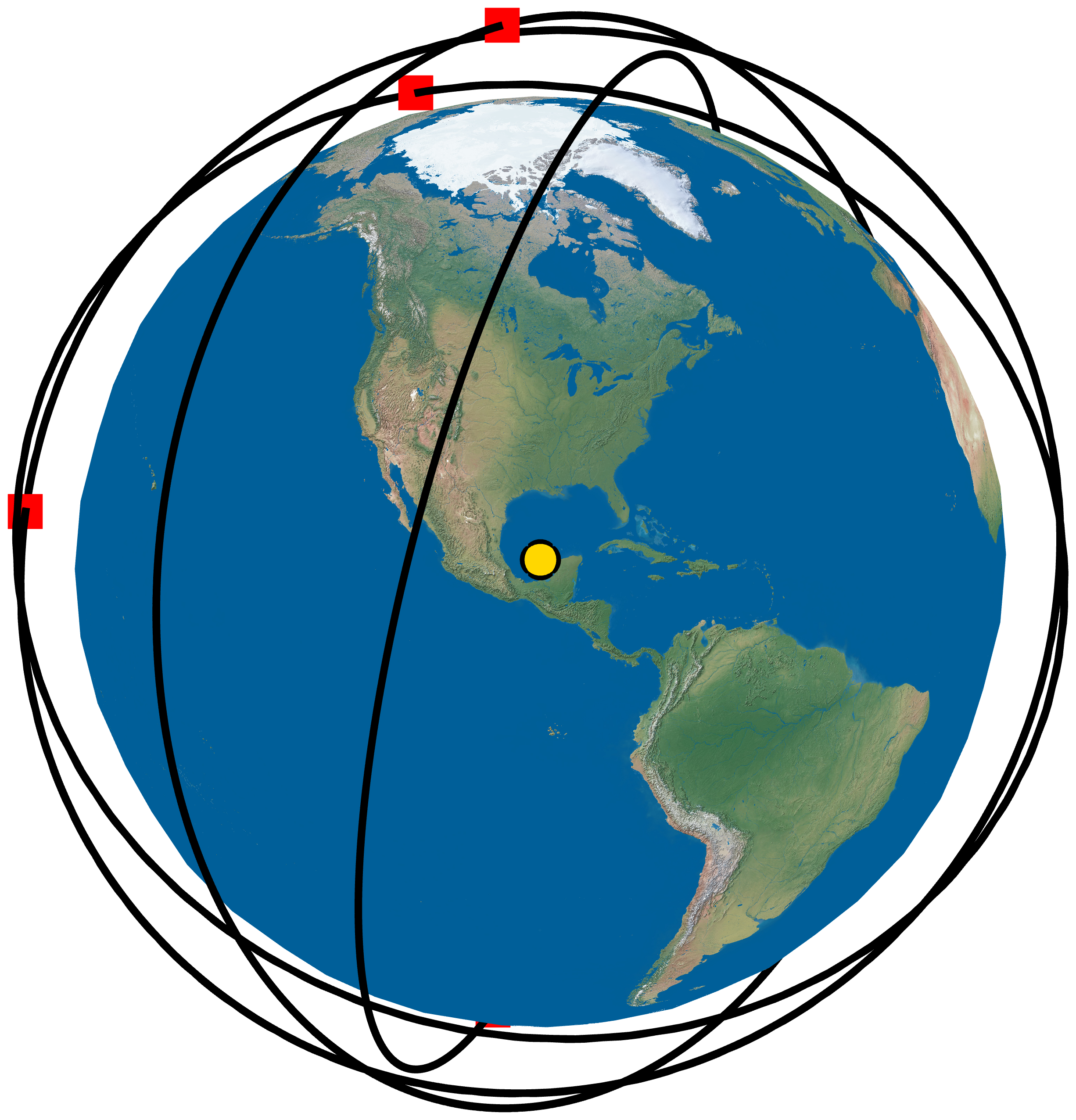} &
                        \includegraphics[height=3.5cm]{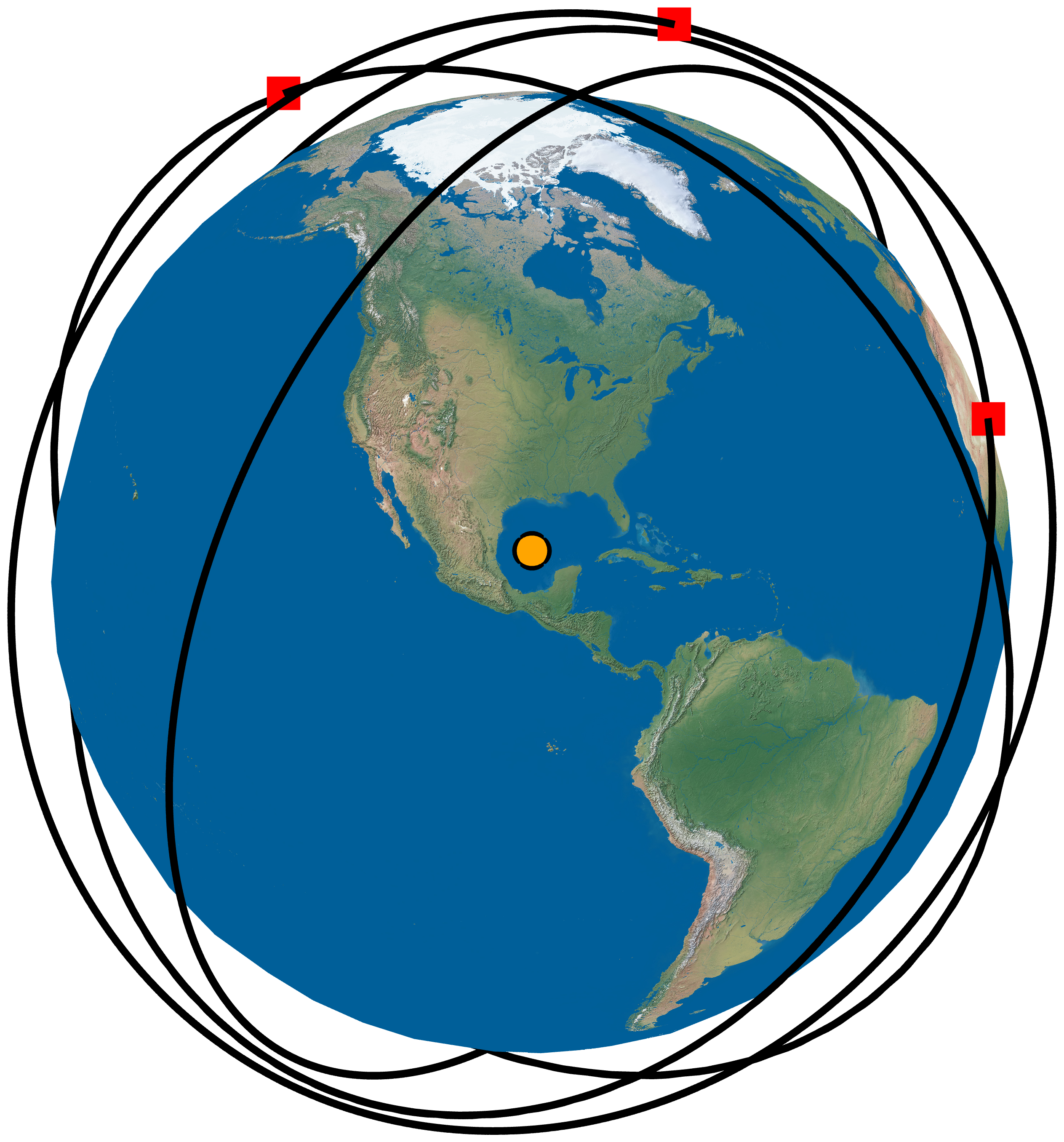} &
                        \includegraphics[height=3.5cm]{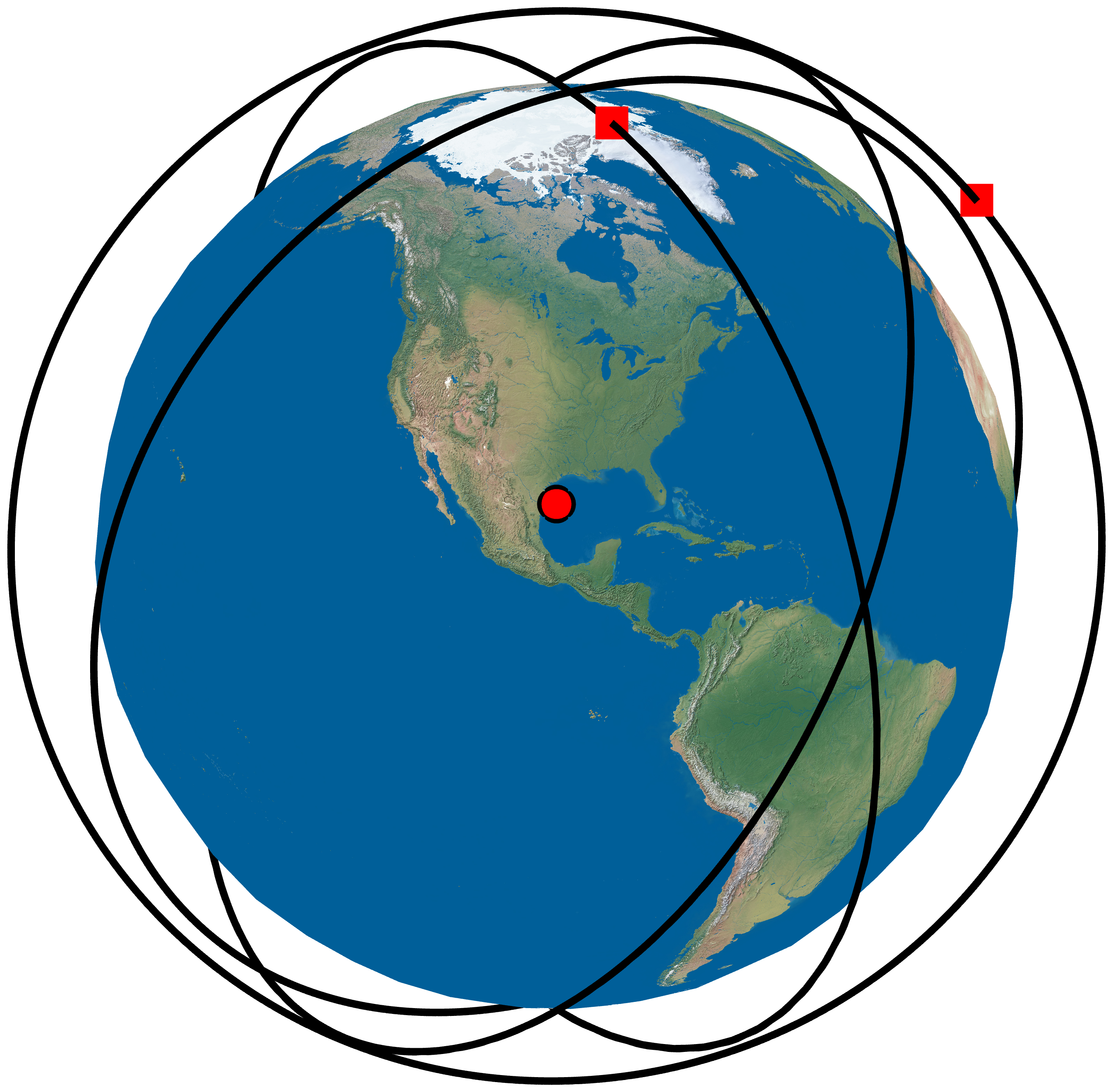} &
                        \includegraphics[height=3.5cm]{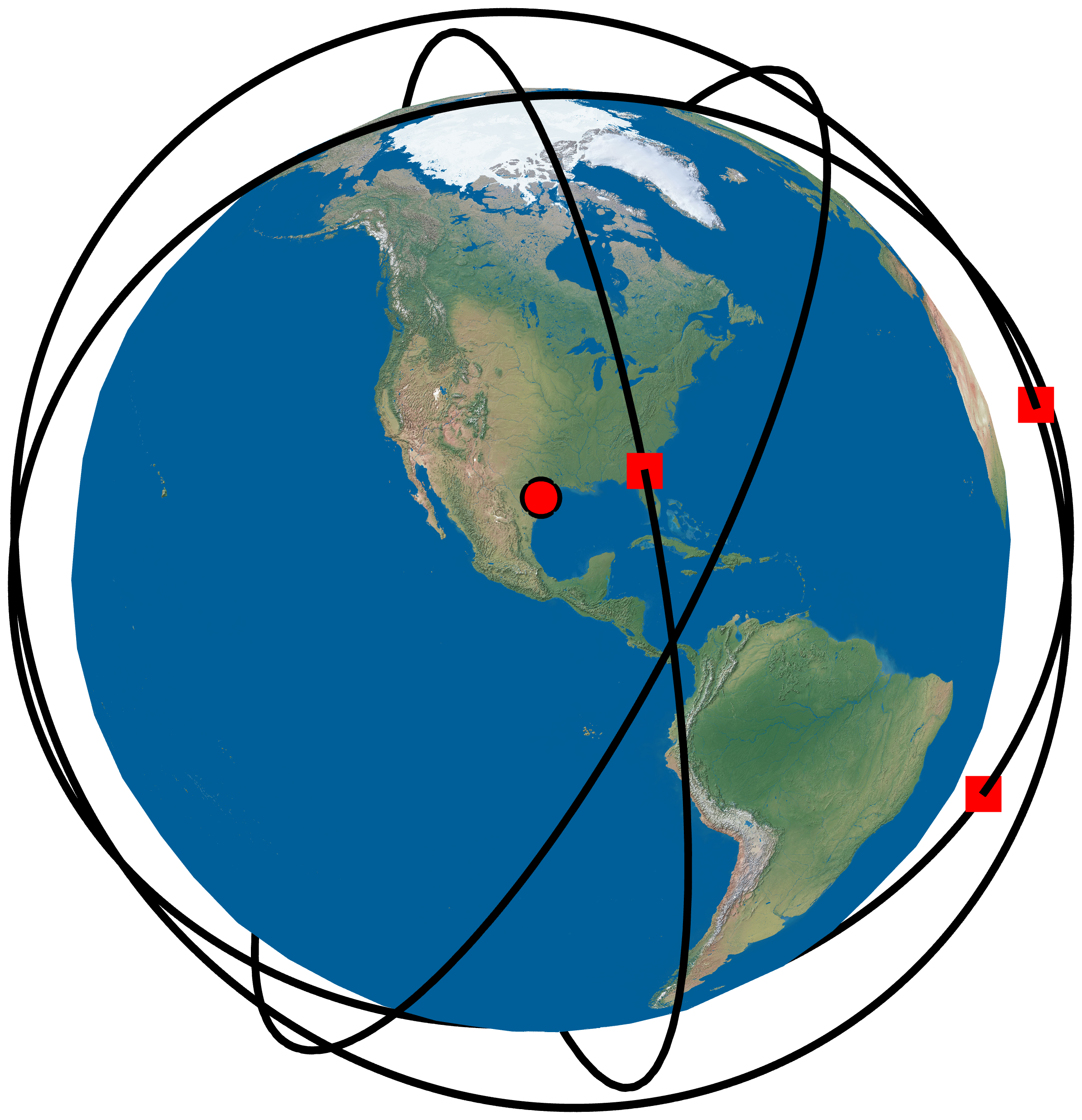} \\
                        \textbf{(a)} Stage 1; target $p_1$ & \textbf{(b)} Stage 2; target $p_3$& \textbf{(c)} Stage 3; target $p_5$ & \textbf{(d)} Stage 4; target $p_7$\\
            \end{tabular}
            \begin{tabular}{cccc}
                        \includegraphics[height=3.5cm]{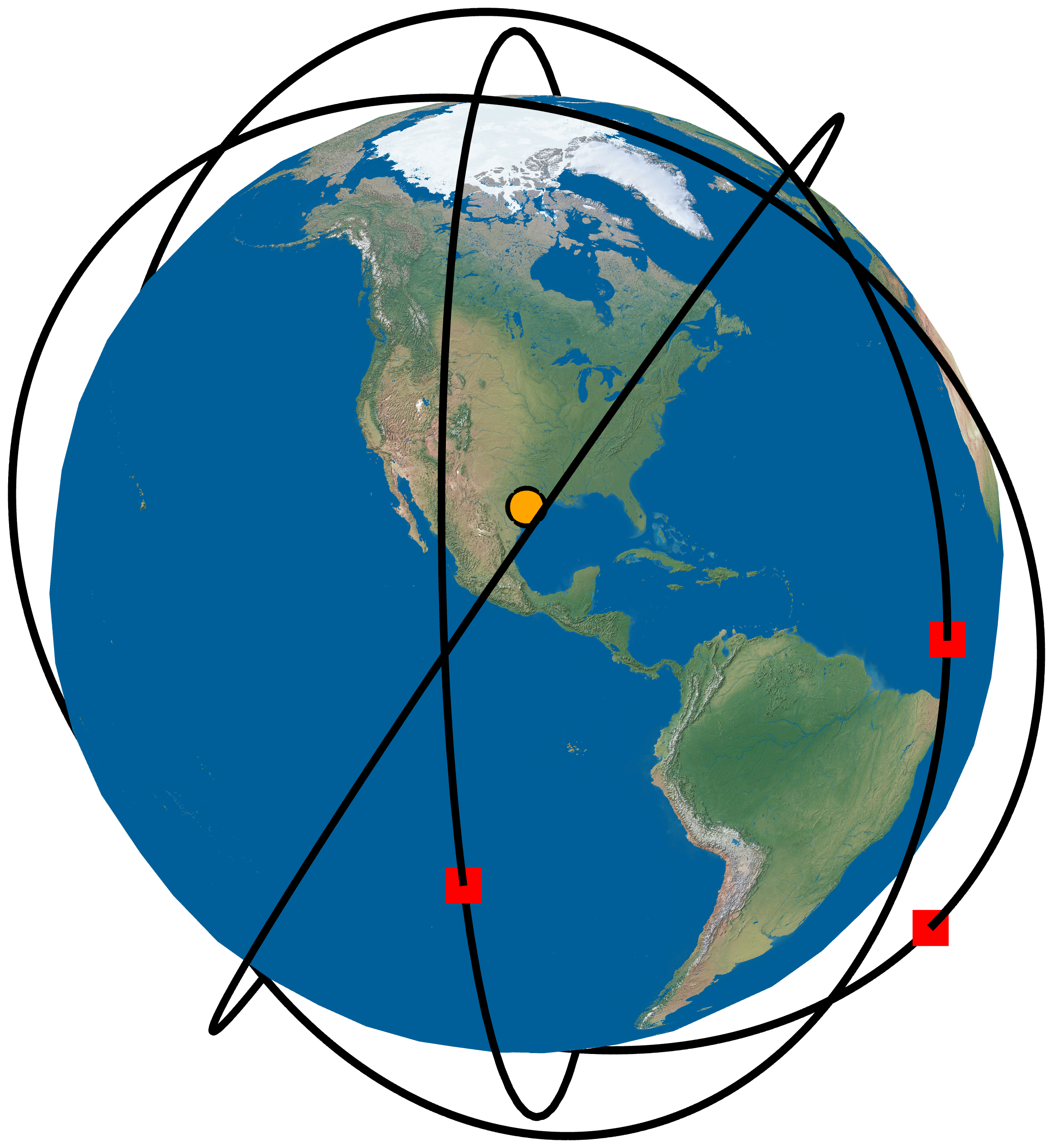} &
                        \includegraphics[height=3.5cm]{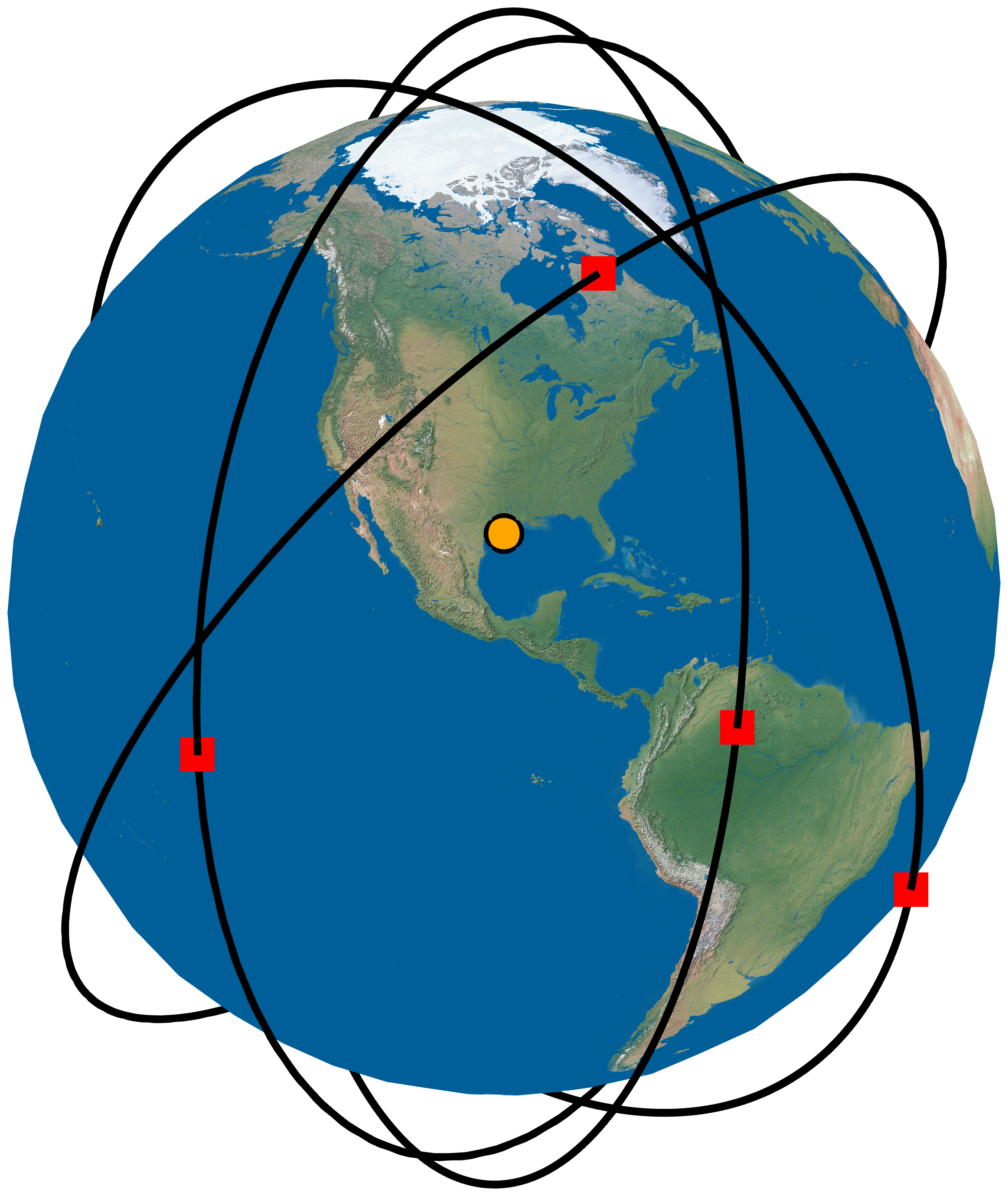} &
                        \includegraphics[height=3.5cm]{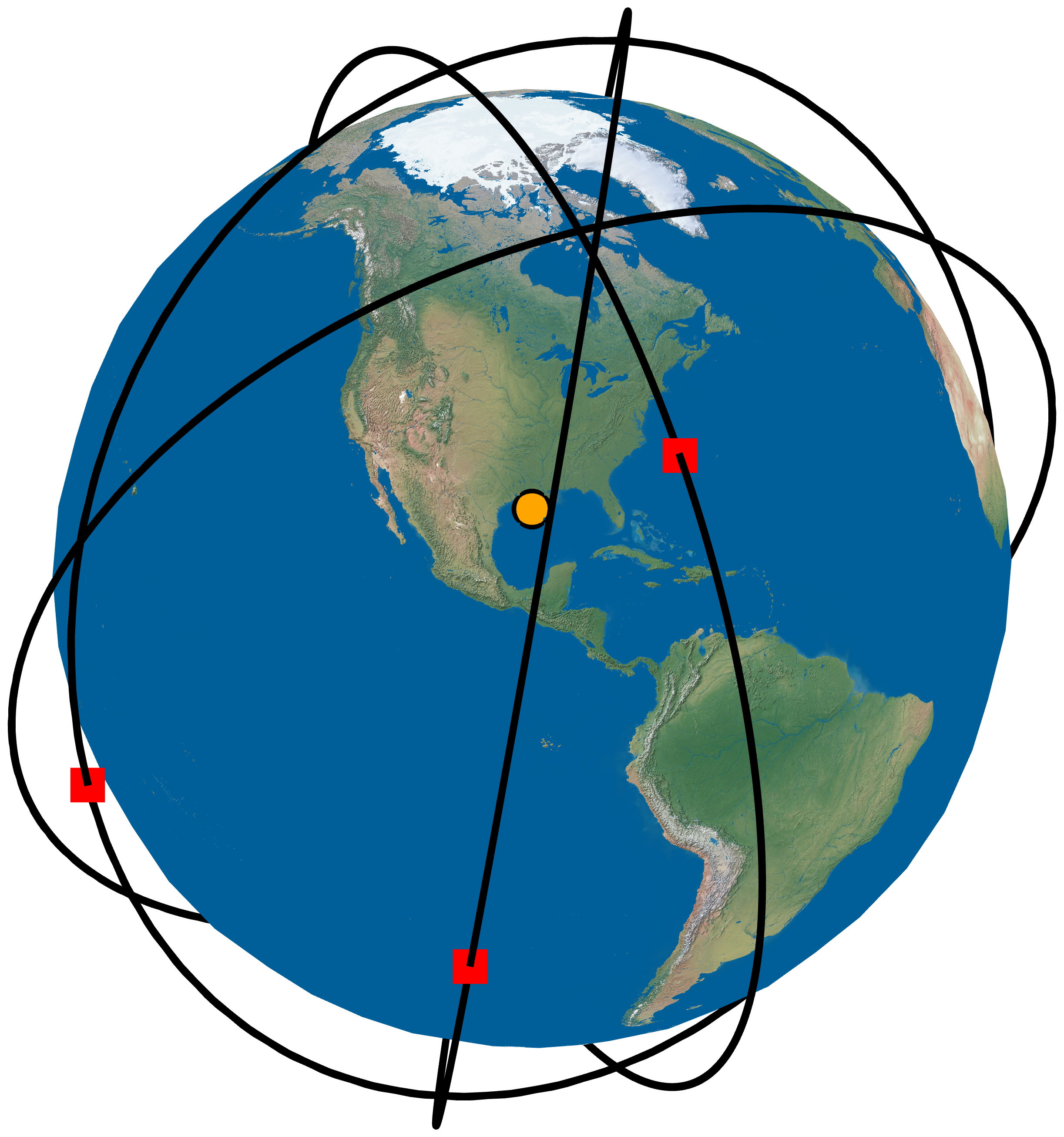} &
                        \includegraphics[height=3.5cm]{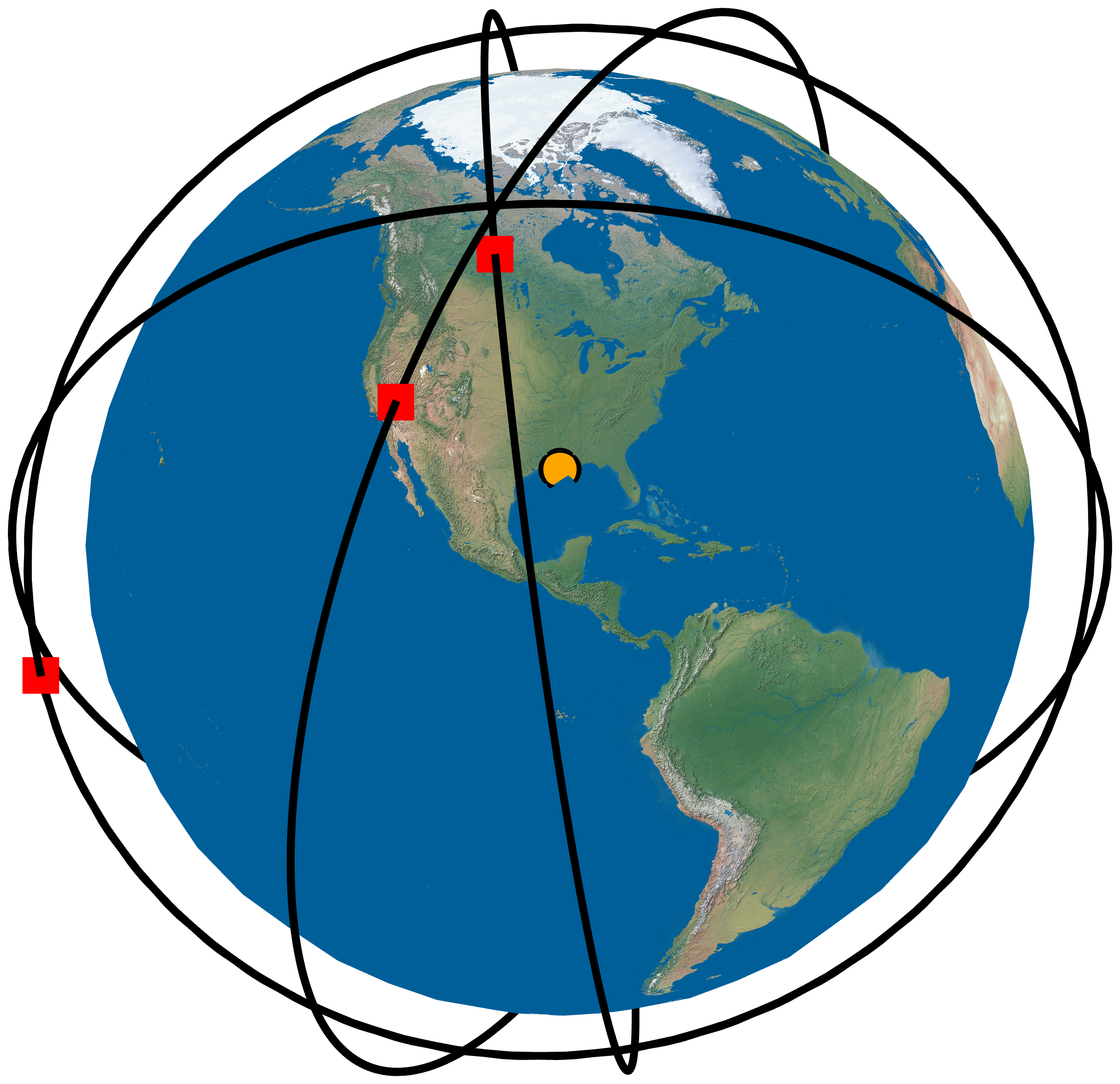} \\
                        \textbf{(e)} Stage 5; target $p_9$ &\textbf{(f)} Stage 6; target $p_{11}$& \textbf{(g)} Stage 7; target $p_{13}$ & \textbf{(h)} Stage 8; target $p_{15}$\\
                        \\
            \end{tabular}
            \caption{Constellation reconfiguration and Hurricane Harvey's evolution. Red squares represent satellites, and circles indicate Hurricane Harvey; their colors indicate their intensity.}
            \label{fig:3dbyStage}
\end{figure}

\begin{figure}[!h]
            \centering
            \renewcommand{\arraystretch}{0.8}
            \begin{tabular}{cccc}
                        \includegraphics[height=3.5cm]{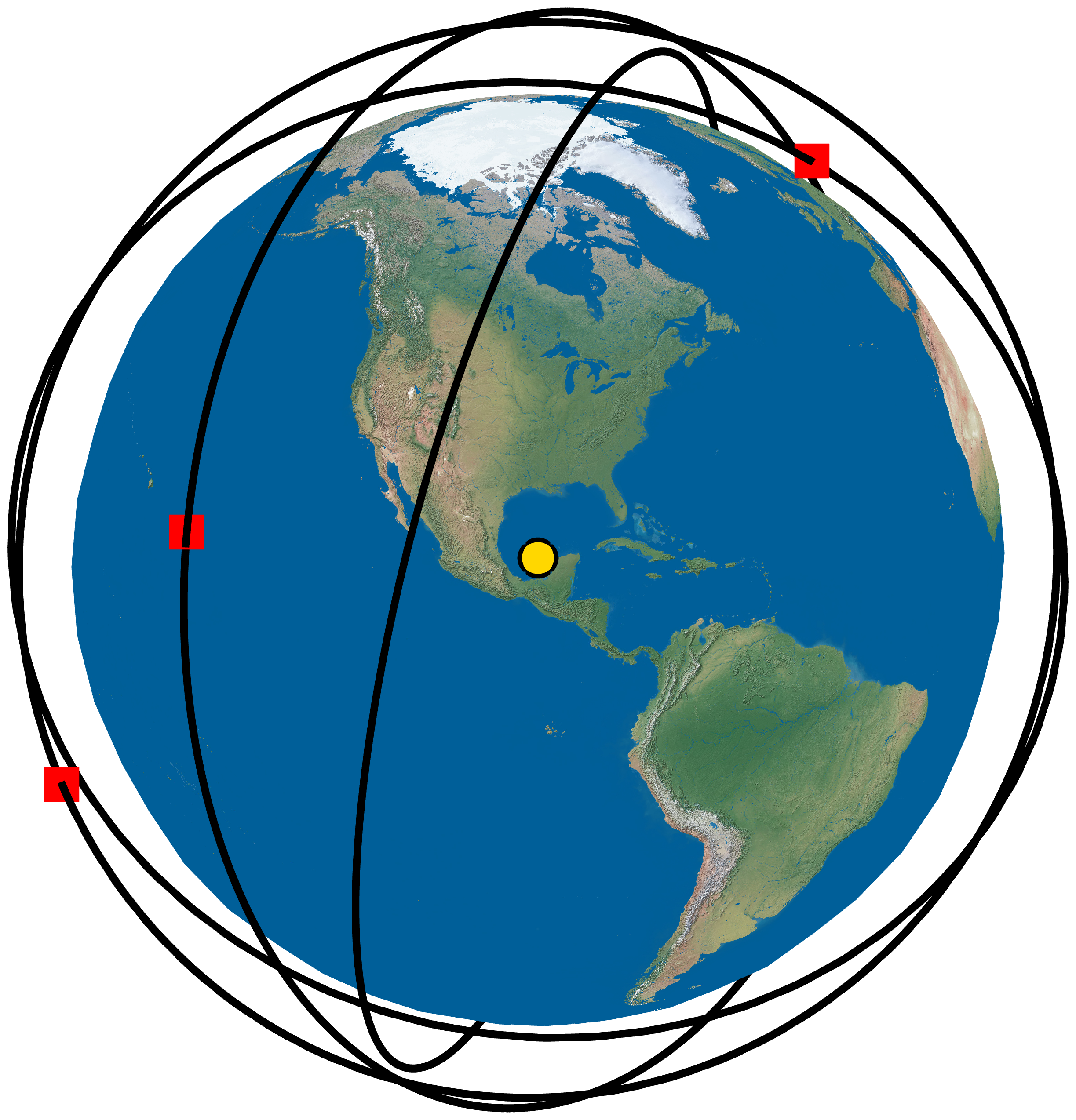} &
                        \includegraphics[height=3.5cm]{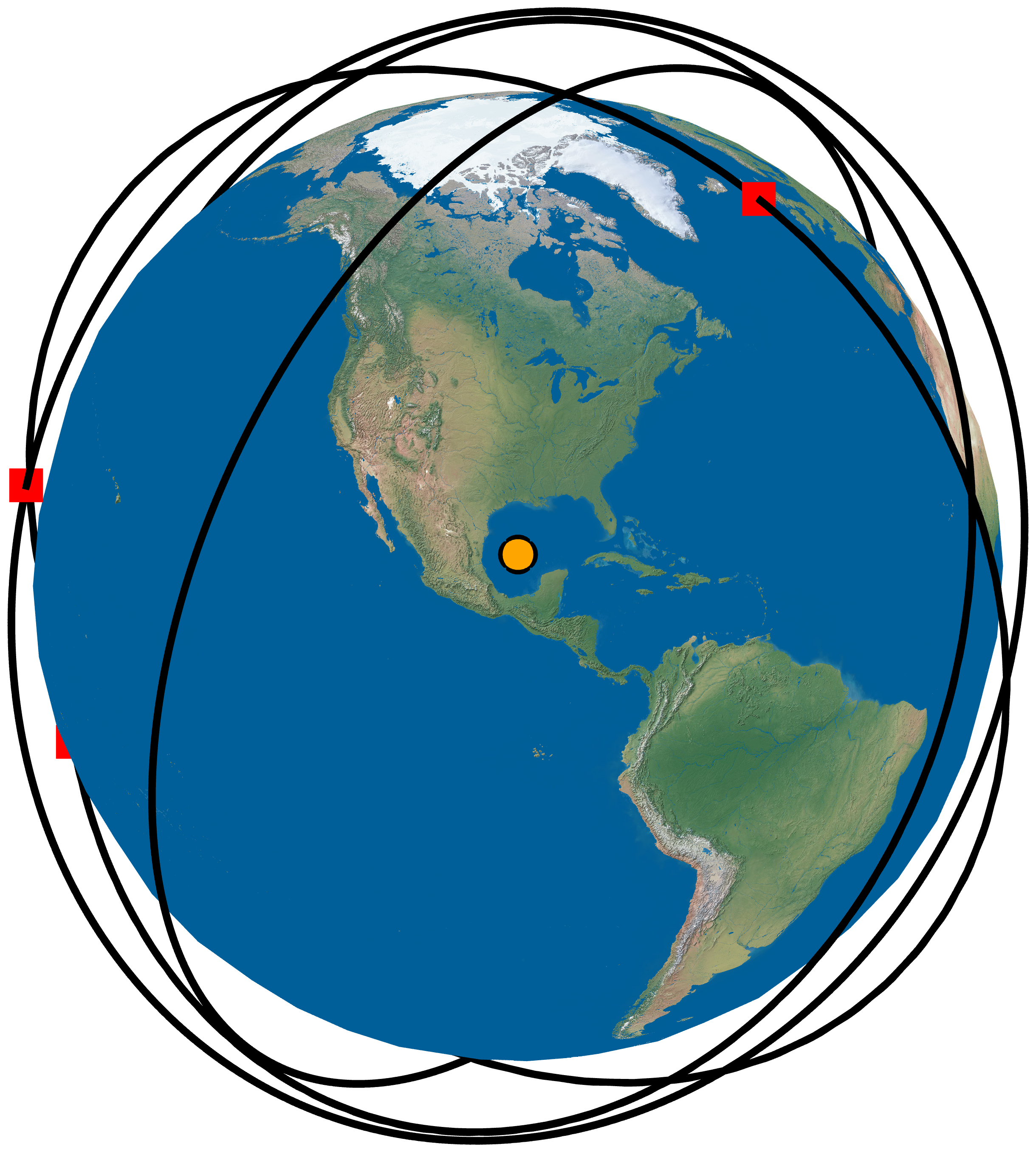} &
                        \includegraphics[height=3.5cm]{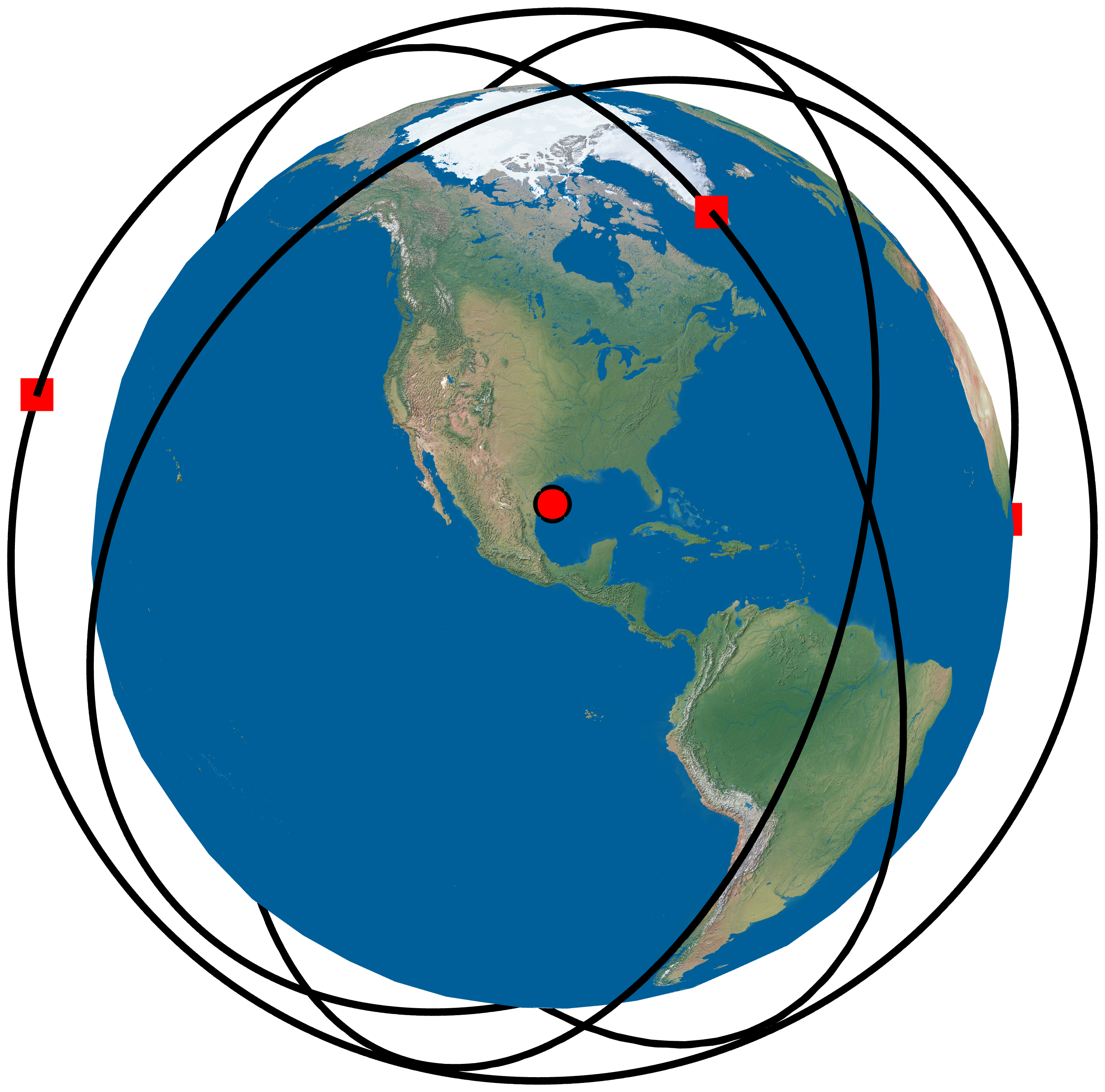} &
                        \includegraphics[height=3.5cm]{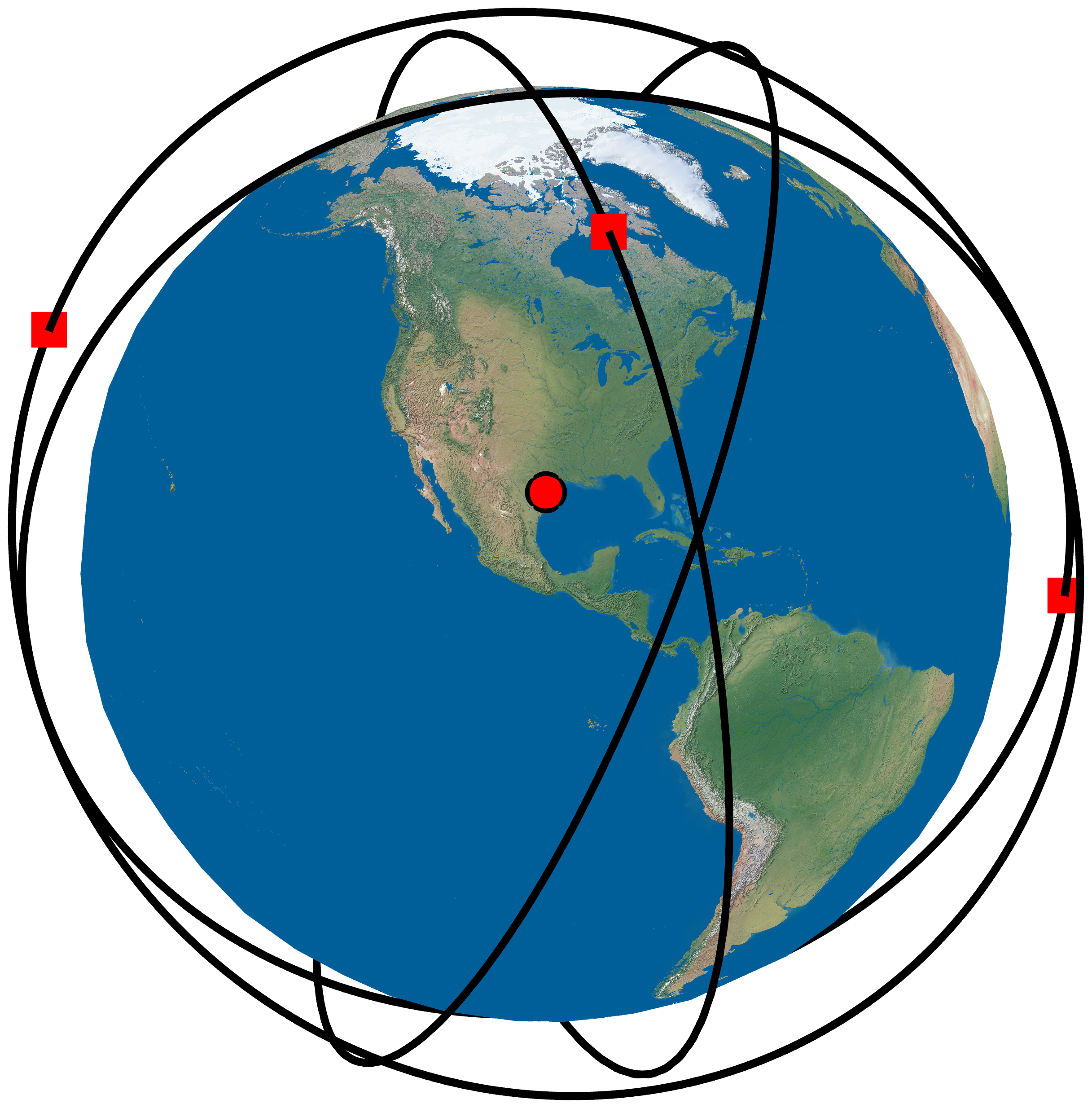} \\
                        \textbf{(a)} Stage 1; target $p_1$ & \textbf{(b)} Stage 2; target $p_3$& \textbf{(c)} Stage 3; target $p_5$ & \textbf{(d)} Stage 4; target $p_7$\\
            \end{tabular}
            \begin{tabular}{cccc}
                        \includegraphics[height=3.5cm]{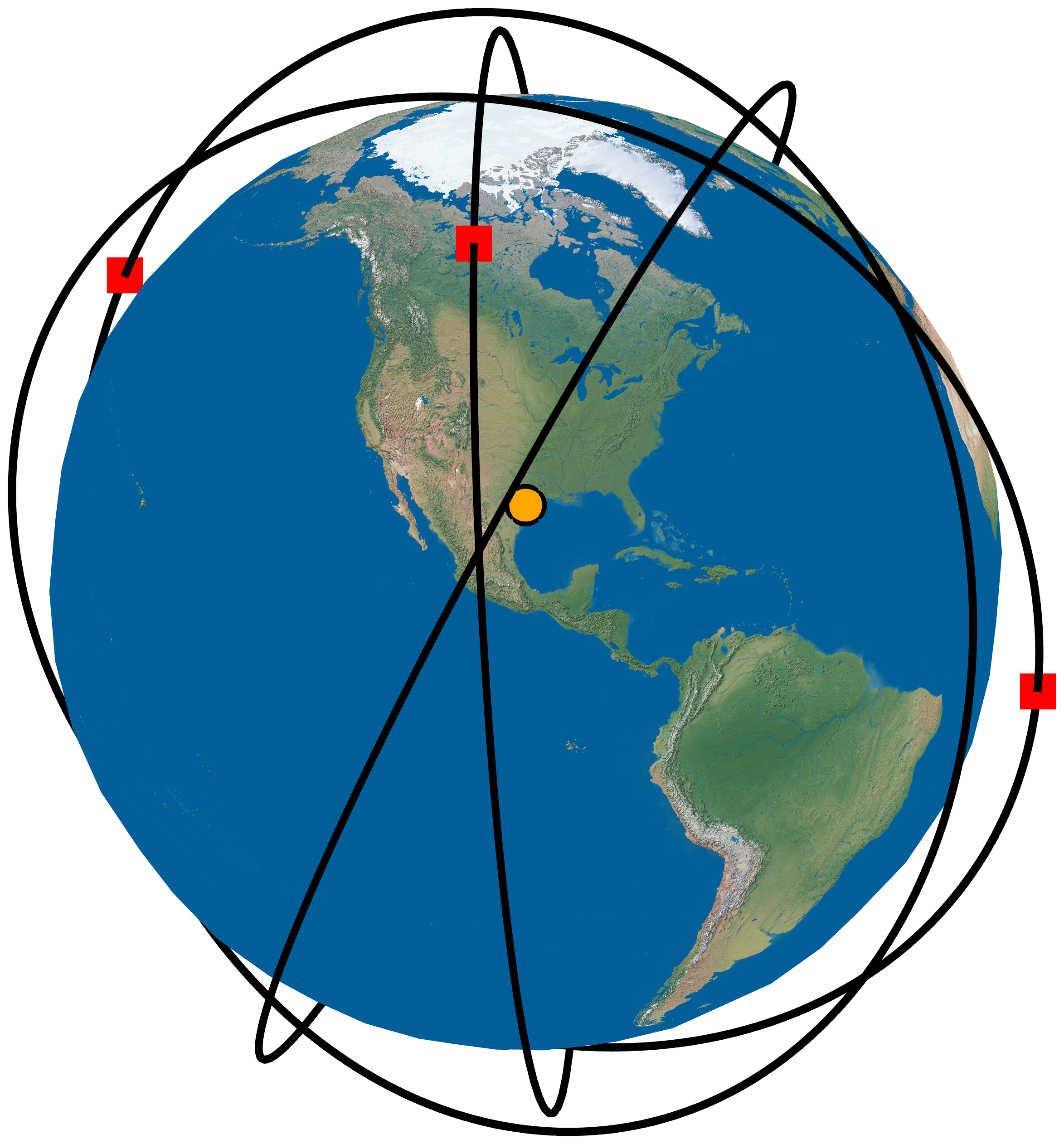} &
                        \includegraphics[height=3.5cm]{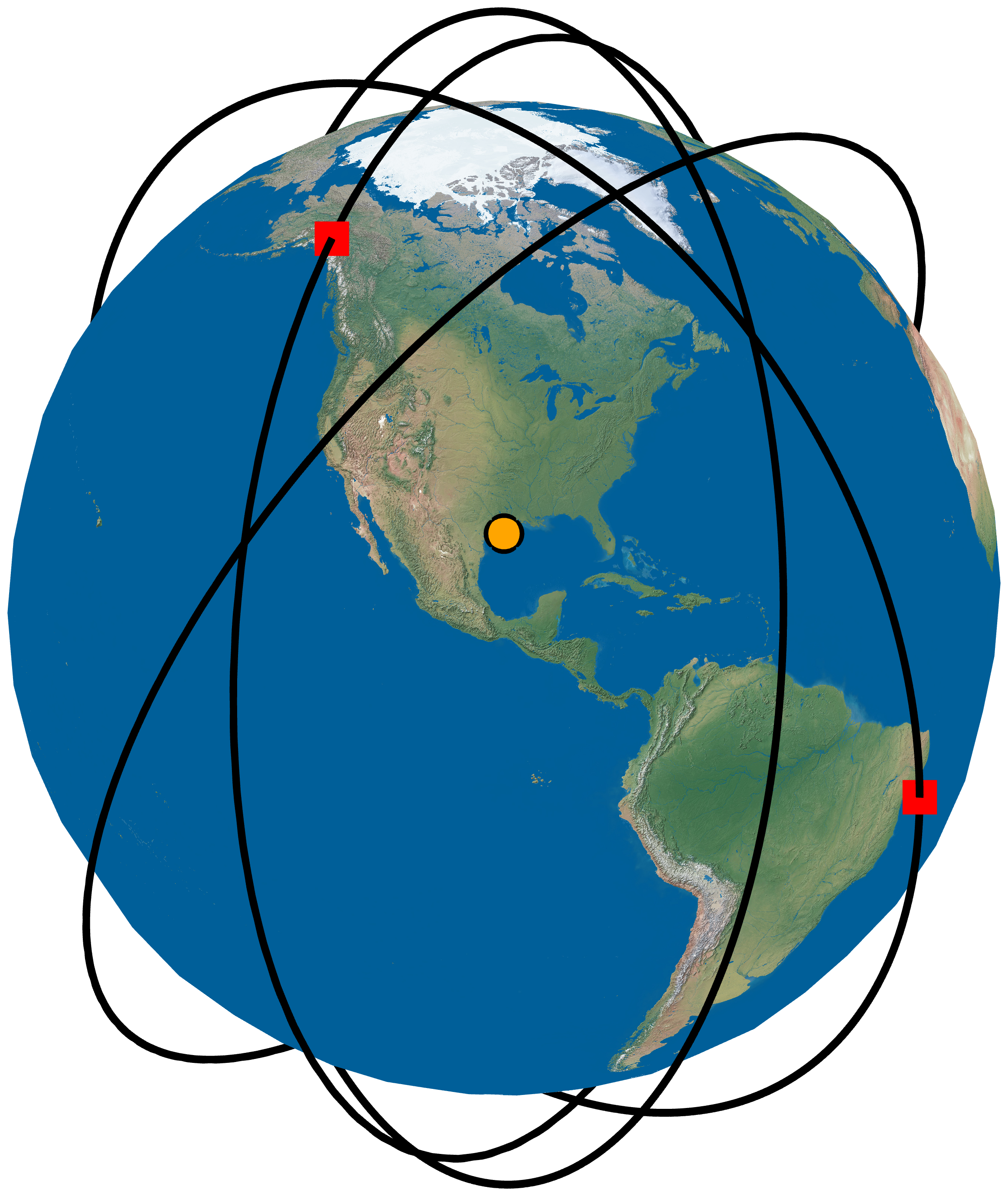} &
                        \includegraphics[height=3.5cm]{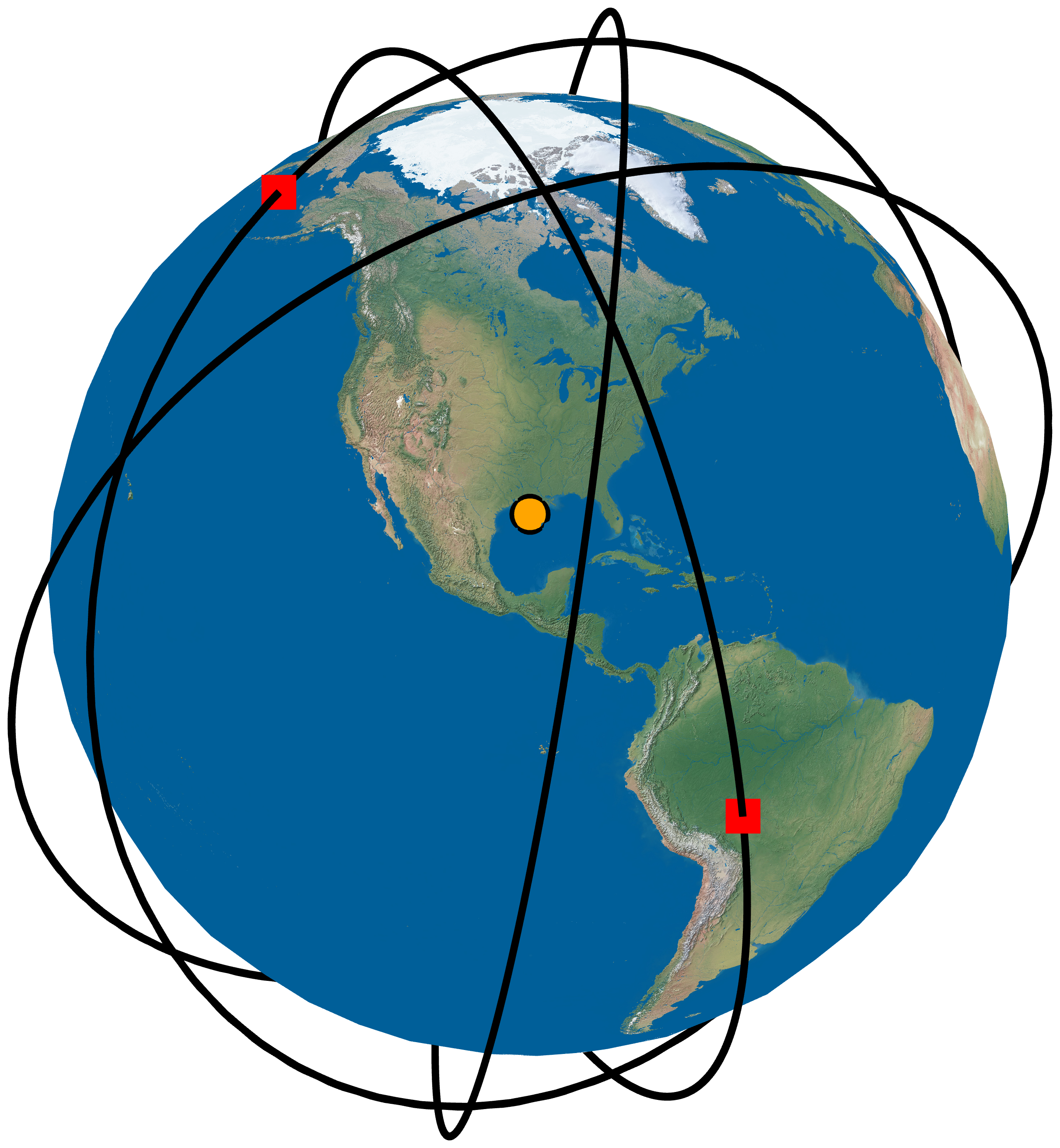} &
                        \includegraphics[height=3.5cm]{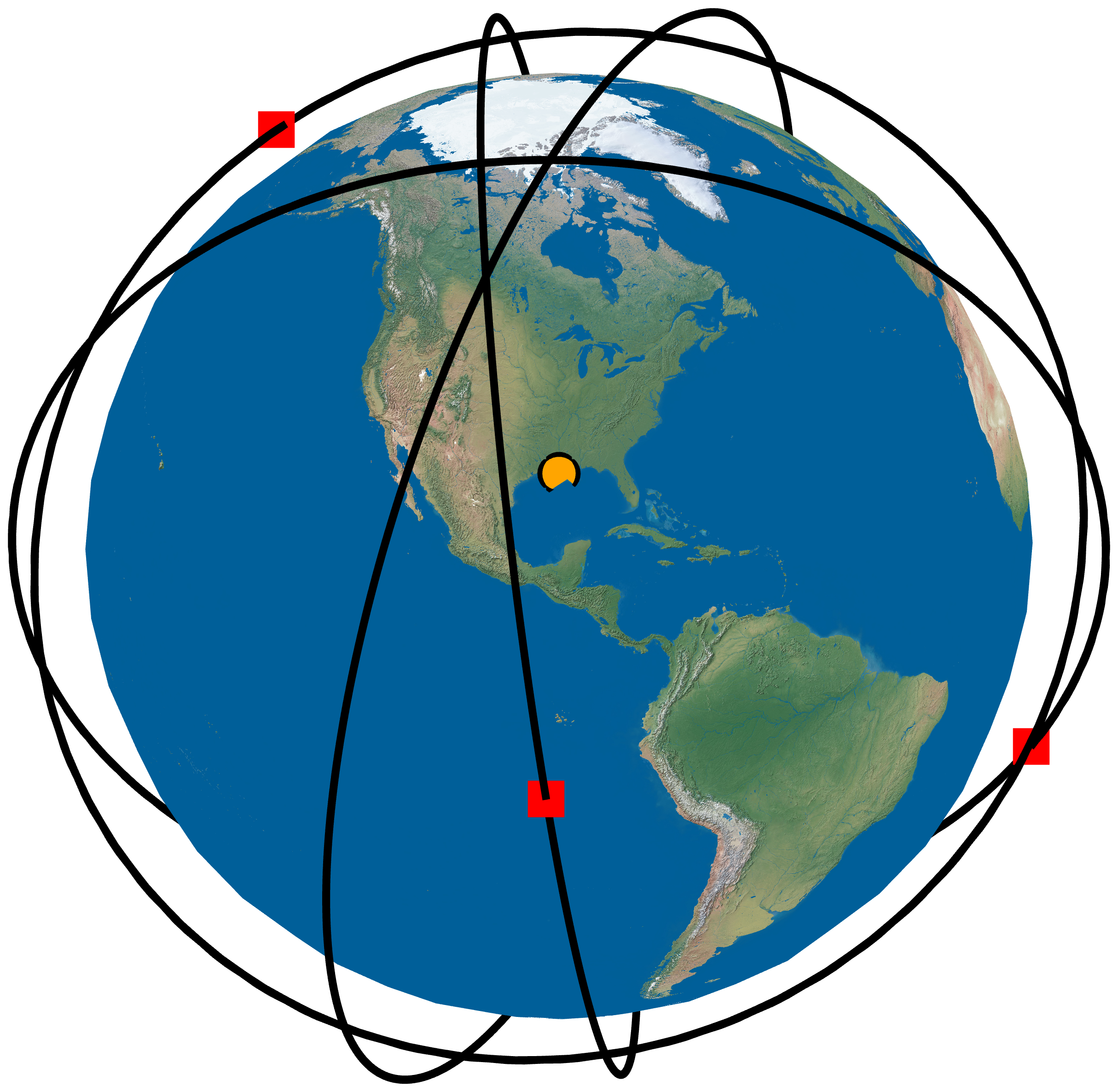} \\
                        \textbf{(e)} Stage 5; target $p_9$ &\textbf{(f)} Stage 6; target $p_{11}$& \textbf{(g)} Stage 7; target $p_{13}$ & \textbf{(h)} Stage 8; target $p_{15}$\\
                        \\
            \end{tabular}
            \caption{Baseline constellation and Hurricane Harvey's evolution. Red squares represent satellites, and circles indicate Hurricane Harvey; their colors indicate their intensity.}
            \label{fig:3dbyStage_Baseline}
\end{figure}

Figure~\ref{fig:3dbyStage} depicts the configurations of the satellite constellation over the course of the mission horizon for Hurricane Harvey with the consideration of eight stages. Each stage is shown at its associated time and the spot target of interest. In addition, Fig.~\ref{fig:3dbyStage_Baseline} depicts the configurations of the baseline constellation at the same times, allowing a direct comparison between the eight-stage reconfiguration case and the fixed configuration of the baseline constellation. The difference in the figures further reflects the observations made in Sec.~\ref{sec:case_study}, showing that even with limited changes in orbital planes, simple phasing maneuvers combined with multiple opportunities for reconfiguration significantly improve the observations realized by the constellation.


\section*{Acknowledgment}
This material is based upon work supported partially by TelePIX and partially by the National Science Foundation Graduate Research Fellowship Program under Grant No. DGE-2039655. Any opinions, findings, and conclusions or recommendations expressed in this material are those of the author(s) and do not necessarily reflect the views of the National Science Foundation. The authors thank the anonymous reviewers for their insightful suggestions, which have greatly contributed to improving the overall quality of the paper.

\bibliography{references}

\begin{thebibliography}{47}
\newcommand{\enquote}[1]{``#1''}
\providecommand{\natexlab}[1]{#1}
\providecommand{\url}[1]{\texttt{#1}}
\providecommand{\urlprefix}{URL }
\expandafter\ifx\csname urlstyle\endcsname\relax
  \providecommand{\doi}[1]{\discretionary{}{}{}https://doi.org/#1}\else
  \providecommand{\doi}[1]{\discretionary{}{}{}\urlstyle{rm}\url{https://doi.org/#1}}\fi

\bibitem[{Chen et~al.(2021)Chen, Wang, and Chu}]{Chen2021Ocean}
Chen, G., Wang, Q., and Chu, X., \enquote{Accelerated spread of Fukushima’s
  waste water by ocean circulation,} \emph{The Innovation}, Vol.~2, 2021, p.
  100119.
\newblock \doi{10.1016/j.xinn.2021.100119}.

\bibitem[{Irrgang et~al.(2019)Irrgang, Saynisch, and Thomas}]{Irrgang2019Ocean}
Irrgang, C., Saynisch, J., and Thomas, M., \enquote{Estimating global ocean
  heat content from tidal magnetic satellite observations,} \emph{Scientific
  Reports}, Vol.~9, No.~1, 2019.
\newblock \doi{10.1038/s41598-019-44397-8}.

\bibitem[{Rogier~de Jong and Dent(2011)}]{Jong2011LandDeg}
Rogier~de Jong, M.~S., Sytze de~Bruin, and Dent, D., \enquote{Quantitative
  mapping of global land degradation using Earth observations,}
  \emph{International Journal of Remote Sensing}, Vol.~32, No.~21, 2011, pp.
  6823--6853.
\newblock \doi{10.1080/01431161.2010.512946}.

\bibitem[{Guo(2010)}]{Guo2010Disaster}
Guo, H., \enquote{Understanding global natural disasters and the role of earth
  observation,} \emph{International Journal of Digital Earth}, Vol.~3, No.~3,
  2010, pp. 221--230.
\newblock \doi{10.1080/17538947.2010.499662}.

\bibitem[{Politi et~al.(2020)Politi, Paterson, Scarrott, Tuohy, O’Mahony, and
  Cámaro-García}]{Politi2020Coasts}
Politi, E., Paterson, S.~K., Scarrott, R., Tuohy, E., O’Mahony, C., and
  Cámaro-García, W.~C., \enquote{Earth observation applications for coastal
  sustainability: potential and challenges for implementation,} \emph{Crafting
  Options, Approaches, and Solutions Towards Sustainability (COASTS) for
  Coastal Regions of the World}, Vol.~1, No.~1, 2020, pp. 306--329.
\newblock \doi{10.1139/anc-2018-0015}.

\bibitem[{Duncan et~al.(2003)Duncan, Martin, Staudt, Yevich, and
  Logan}]{Duncan2003Biomass}
Duncan, B.~N., Martin, R.~V., Staudt, A.~C., Yevich, R., and Logan, J.~A.,
  \enquote{Interannual and seasonal variability of biomass burning emissions
  constrained by satellite observations,} \emph{Journal of Geophysical
  Research: Atmospheres}, Vol. 108, No.~D2, 2003, pp. ACH 1--1--ACH 1--22.
\newblock \doi{https://doi.org/10.1029/2002JD002378}.

\bibitem[{Srivastava et~al.(2018)Srivastava, Piles, and
  Pearson}]{Srivastava2018Drought}
Srivastava, P., Piles, M., and Pearson, S., \enquote{Earth Observation-Based
  Operational Estimation of Soil Moisture and Evapotranspiration for
  Agricultural Crops in Support of Sustainable Water Management,}
  \emph{Sustainability}, Vol.~10, 2018, p. 181.
\newblock \doi{10.3390/su10010181}.

\bibitem[{{van der Meer} et~al.(2012){van der Meer}, {van der Werff}, {van
  Ruitenbeek}, Hecker, Bakker, Noomen, {van der Meijde}, Carranza, de~Smeth,
  and Woldai}]{VANDERMEER2012Geology}
{van der Meer}, F.~D., {van der Werff}, H.~M., {van Ruitenbeek}, F.~J., Hecker,
  C.~A., Bakker, W.~H., Noomen, M.~F., {van der Meijde}, M., Carranza, E.
  J.~M., de~Smeth, J.~B., and Woldai, T., \enquote{Multi- and hyperspectral
  geologic remote sensing: A review,} \emph{International Journal of Applied
  Earth Observation and Geoinformation}, Vol.~14, No.~1, 2012, pp. 112--128.
\newblock \doi{https://doi.org/10.1016/j.jag.2011.08.002}.

\bibitem[{Sogno et~al.(2020)Sogno, Traidl-Hoffmann, and
  Kuenzer}]{Sogno2020Disease}
Sogno, P., Traidl-Hoffmann, C., and Kuenzer, C., \enquote{Earth Observation
  Data Supporting Non-Communicable Disease Research: A Review,} \emph{Remote
  Sensing}, Vol.~12, 2020, p.~34.
\newblock \doi{10.3390/rs12162541}.

\bibitem[{Hansen and Loveland(2012)}]{HANSEN2012LandCover}
Hansen, M.~C., and Loveland, T.~R., \enquote{A review of large area monitoring
  of land cover change using Landsat data,} \emph{Remote Sensing of
  Environment}, Vol. 122, 2012, pp. 66--74.
\newblock \doi{https://doi.org/10.1016/j.rse.2011.08.024}, landsat Legacy
  Special Issue.

\bibitem[{Kelly et~al.(2009)Kelly, Loverro, Case, Quéruel, Maréchal, and
  Barroso}]{Kelly2009}
Kelly, A.~C., Loverro, A., Case, W.~F., Quéruel, N., Maréchal, C., and
  Barroso, T., \enquote{Small Earth observing satellites flying with large
  satellites in the A-train,} \emph{Small Satellite Missions for Earth
  Observation}, 2009, p. 19–28.
\newblock \doi{10.1007/978-3-642-03501-2_2}.

\bibitem[{Gierach and Subrahmanyam(2007)}]{Gierach2007}
Gierach, M.~M., and Subrahmanyam, B., \enquote{Satellite Data Analysis of the
  Upper Ocean Response to Hurricanes Katrina and Rita (2005) in the Gulf of
  Mexico,} \emph{IEEE Geoscience and Remote Sensing Letters}, Vol.~4, No.~1,
  2007, pp. 132--136.
\newblock \doi{10.1109/LGRS.2006.887145}.

\bibitem[{Stephens et~al.(2003)Stephens, Cooksley, Da~Silva~Curiel, Boland,
  Jason, Northham, Brewer, Anzalchi, Newell, Underwood, Machin, Sun, and
  Sweeting}]{Stephens2003DMC}
Stephens, P., Cooksley, J., Da~Silva~Curiel, A., Boland, L., Jason, S.,
  Northham, J., Brewer, A., Anzalchi, J., Newell, H., Underwood, C., Machin,
  S., Sun, W., and Sweeting, S., \enquote{Launch of the international Disaster
  Monitoring Constellation; the development of a novel international
  partnership in space,} \emph{International Conference on Recent Advances in
  Space Technologies (RAST)}, 2003, pp. 525 -- 535.
\newblock \doi{10.1109/RAST.2003.1303972}.

\bibitem[{{de Weck} et~al.(2008){de Weck}, Scialom, and
  Siddiqi}]{deweck2008optimal}
{de Weck}, O.~L., Scialom, U., and Siddiqi, A., \enquote{Optimal
  reconfiguration of satellite constellations with the auction algorithm,}
  \emph{Acta Astronautica}, Vol.~62, No.~2, 2008, pp. 112--130.
\newblock \doi{10.1016/j.actaastro.2007.02.008}.

\bibitem[{He et~al.(2020)He, Li, Yang, and Zhao}]{he2020reconfigurable}
He, X., Li, H., Yang, L., and Zhao, J., \enquote{{Reconfigurable Satellite
  Constellation Design for Disaster Monitoring Using Physical Programming},}
  \emph{International Journal of Aerospace Engineering}, Vol. 2020, 2020, p.
  8813685.
\newblock \doi{10.1155/2020/8813685}.

\bibitem[{Xiaoyu et~al.(2022)Xiaoyu, Bai, Xu, Li, Zhou, Yu, and
  Zhang}]{Zuo2022Surrogate}
Xiaoyu, Z., Bai, X., Xu, M., Li, M., Zhou, J., Yu, L., and Zhang, J.,
  \enquote{Satellite Constellation Reconfiguration Using Surrogate-Based
  Optimization,} \emph{Journal of Aerospace Engineering}, Vol.~35, 2022.
\newblock \doi{10.1061/(ASCE)AS.1943-5525.0001438}.

\bibitem[{Lee and Ho(2023)}]{lee2023regional}
Lee, H., and Ho, K., \enquote{Regional Constellation Reconfiguration Problem:
  Integer Linear Programming Formulation and Lagrangian Heuristic Method,}
  \emph{Journal of Spacecraft and Rockets}, Vol.~60, No.~6, 2023, pp.
  1828--1845.
\newblock \doi{10.2514/1.A35685}.

\bibitem[{Chen et~al.(2015)Chen, Mahalec, Chen, Liu, He, and
  Sun}]{chen2015reconfiguration}
Chen, Y., Mahalec, V., Chen, Y., Liu, X., He, R., and Sun, K.,
  \enquote{Reconfiguration of satellite orbit for cooperative observation using
  variable-size multi-objective differential evolution,} \emph{European Journal
  of Operational Research}, Vol. 242, No.~1, 2015, pp. 10--20.
\newblock \doi{10.1016/j.ejor.2014.09.025}.

\bibitem[{Paek et~al.(2019)Paek, Kim, and de~Weck}]{paek2019optimization}
Paek, S.~W., Kim, S., and de~Weck, O., \enquote{Optimization of Reconfigurable
  Satellite Constellations Using Simulated Annealing and Genetic Algorithm,}
  \emph{Sensors}, Vol.~19, No.~4, 2019.
\newblock \doi{10.3390/s19040765}.

\bibitem[{Appel et~al.(2014)Appel, Guelman, and Mishne}]{Appel2014Optimization}
Appel, L., Guelman, M., and Mishne, D., \enquote{Optimization of satellite
  constellation reconfiguration maneuvers,} \emph{Acta Astronautica}, Vol.~99,
  2014, p. 166–174.
\newblock \doi{10.1016/j.actaastro.2014.02.016}.

\bibitem[{Anderson et~al.(2022)Anderson, Cardin, and
  Grogan}]{Anderson2022megaconstellation}
Anderson, J.~F., Cardin, M.-A., and Grogan, P.~T., \enquote{Design and analysis
  of flexible multi-layer staged deployment for satellite mega-constellations
  under demand uncertainty,} \emph{Acta Astronautica}, Vol. 198, 2022, p.
  179–193.
\newblock \doi{10.1016/j.actaastro.2022.05.022}.

\bibitem[{{Ferringer} et~al.(2009){Ferringer}, {Spencer}, and
  {Reed}}]{ferringer2009many}
{Ferringer}, M.~P., {Spencer}, D.~B., and {Reed}, P., \enquote{Many-objective
  reconfiguration of operational satellite constellations with the
  Large-Cluster Epsilon Non-dominated Sorting Genetic Algorithm-II,} \emph{2009
  IEEE Congress on Evolutionary Computation}, 2009, pp. 340--349.
\newblock \doi{10.1109/CEC.2009.4982967}.

\bibitem[{Morgan et~al.(2023)Morgan, McGrath, and de~Weck}]{Morgan2023}
Morgan, S.~J., McGrath, C.~N., and de~Weck, O.~L., \enquote{Optimization of
  multispacecraft maneuvers for mobile target tracking from low Earth orbit,}
  \emph{Journal of Spacecraft and Rockets}, Vol.~60, No.~2, 2023, p. 581–590.
\newblock \doi{10.2514/1.a35457}.

\bibitem[{McGrath and Macdonald(2019)}]{McGrath2019General}
McGrath, C.~N., and Macdonald, M., \enquote{General Perturbation Method for
  Satellite Constellation Reconfiguration Using Low-Thrust Maneuvers,}
  \emph{Journal of Guidance, Control, and Dynamics}, Vol.~42, No.~8, 2019, pp.
  1676--1692.
\newblock \doi{10.2514/1.G003739}.

\bibitem[{Jiaxin et~al.(2021)Jiaxin, Leping, Huan, and
  Yanwei}]{Jiaxin2021Evolution}
Jiaxin, H., Leping, Y., Huan, H., and Yanwei, Z., \enquote{Optimal
  reconfiguration of constellation using adaptive innovation driven
  multiobjective evolutionary algorithm,} \emph{Journal of Systems Engineering
  and Electronics}, Vol.~32, 2021, pp. 1527--1538.
\newblock \doi{10.23919/JSEE.2021.000128}.

\bibitem[{Lee and Ho(2020)}]{Lee2020binary}
Lee, H., and Ho, K., \enquote{Binary Integer Linear Programming Formulation for
  Optimal Satellite Constellation Reconfiguration,} \emph{AAS/AIAA
  Astrodynamics Specialist Conference}, 2020.

\bibitem[{Lee and Ho(2021)}]{Lee2021lagrangian}
Lee, H., and Ho, K., \enquote{A Lagrangian Relaxation-Based Heuristic Approach
  to Regional Constellation Reconfiguration Problem,} \emph{AAS/AIAA
  Astrodynamics Specialist Conference}, 2021.

\bibitem[{Lee et~al.(2022)Lee, Chen, and Ho}]{Lee2022maximizing}
Lee, H., Chen, H., and Ho, K., \enquote{Maximizing Observation Throughput via
  Multi-Stage Satellite Constellation Reconfiguration,} \emph{AAS/AIAA
  Astrodynamics Specialist Conference}, 2022.

\bibitem[{Denis et~al.(2016)Denis, {de Boissezon}, Hosford, Pasco, Montfort,
  and Ranera}]{denis2016}
Denis, G., {de Boissezon}, H., Hosford, S., Pasco, X., Montfort, B., and
  Ranera, F., \enquote{The evolution of Earth Observation satellites in Europe
  and its impact on the performance of emergency response services,} \emph{Acta
  Astronautica}, Vol. 127, 2016, pp. 619--633.
\newblock \doi{10.1016/j.actaastro.2016.06.012}.

\bibitem[{Voigt et~al.(2016)Voigt, Giulio-Tonolo, Lyons, Kučera, Jones,
  Schneiderhan, Platzeck, Kaku, Hazarika, Czaran, Li, Pedersen, James, Proy,
  Muthike, Bequignon, and Guha-Sapir}]{voigt2016global}
Voigt, S., Giulio-Tonolo, F., Lyons, J., Kučera, J., Jones, B., Schneiderhan,
  T., Platzeck, G., Kaku, K., Hazarika, M.~K., Czaran, L., Li, S., Pedersen,
  W., James, G.~K., Proy, C., Muthike, D.~M., Bequignon, J., and Guha-Sapir,
  D., \enquote{Global trends in satellite-based emergency mapping,}
  \emph{Science}, Vol. 353, No. 6296, 2016, pp. 247--252.
\newblock \doi{10.1126/science.aad8728}.

\bibitem[{Lee et~al.(2020)Lee, Shimizu, Yoshikawa, and Ho}]{lee2020satellite}
Lee, H., Shimizu, S., Yoshikawa, S., and Ho, K., \enquote{Satellite
  Constellation Pattern Optimization for Complex Regional Coverage,}
  \emph{Journal of Spacecraft and Rockets}, Vol.~57, No.~6, 2020, pp.
  1309--1327.
\newblock \doi{10.2514/1.A34657}.

\bibitem[{{Williams Rogers} et~al.(2023){Williams Rogers}, Kim, Lee, Kim, and
  Lee}]{david2023designing}
{Williams Rogers}, D.~O., Kim, S.-W., Lee, M., Kim, Y.-H., and Lee, H.,
  \enquote{Designing Optimal Satellite Constellation Patterns with Facility
  Location Problem Models and Mixed Integer Linear Programming,} \emph{ASCEND
  2023}, 2023.
\newblock \doi{10.2514/6.2023-4658},
  \urlprefix\url{https://arc.aiaa.org/doi/abs/10.2514/6.2023-4658}.

\bibitem[{Powell(2009)}]{powell2009}
Powell, W.~B., \enquote{What you should know about approximate dynamic
  programming,} \emph{Naval Research Logistics (NRL)}, Vol.~56, No.~3, 2009,
  pp. 239--249.
\newblock \doi{10.1002/nav.20347}.

\bibitem[{Sethi and Sorger(1991)}]{Sethi1991theory}
Sethi, S., and Sorger, G., \enquote{{A theory of rolling horizon decision
  making},} \emph{Annals of Operations Research}, Vol.~29, No.~1, 1991, pp.
  387--415.
\newblock \doi{10.1007/BF02283607}.

\bibitem[{Megiddo et~al.(1983)Megiddo, Zemel, and Hakimi}]{megiddo1983}
Megiddo, N., Zemel, E., and Hakimi, S.~L., \enquote{The Maximum Coverage
  Location Problem,} \emph{SIAM Journal on Algebraic Discrete Methods}, Vol.~4,
  No.~2, 1983, pp. 253--261.
\newblock \doi{10.1137/0604028}.

\bibitem[{Ulybyshev(2008)}]{ulybyshev2008satellite}
Ulybyshev, Y., \enquote{Satellite constellation design for complex coverage,}
  \emph{Journal of Spacecraft and Rockets}, Vol.~45, No.~4, 2008, pp. 843--849.

\bibitem[{MATLAB(2023)}]{MATLAB}
MATLAB, \emph{23.2.0.2459199 (R2023b) Update 5}, The MathWorks Inc., Natick,
  Massachusetts, 2023.

\bibitem[{Vallado(2013)}]{vallado2013fundamentals}
Vallado, D., \emph{Fundamentals of Astrodynamics and Applications}, Space
  technology library, Microcosm Press, 2013.

\bibitem[{{NOAA National Weather Service}(2017)}]{HarveyInfo}
{NOAA National Weather Service}, \enquote{Major Hurricane Harvey,} , 2017.
\newblock \urlprefix\url{https://www.weather.gov/crp/hurricane_harvey}, last
  Accessed January 17, 2024.

\bibitem[{{NOAA's National Centers for Environmental Information, The National
  Hurricane Center}(2023)}]{HarveyCost}
{NOAA's National Centers for Environmental Information, The National Hurricane
  Center}, \enquote{Costliest U.S. Tropical Cyclones,} , 12 2023.
\newblock \urlprefix\url{https://www.ncei.noaa.gov/access/billions/dcmi.pdf},
  last Accessed January 17, 2024.

\bibitem[{{NASA Earth Observatory}(2017)}]{harveyimage}
{NASA Earth Observatory}, \enquote{Hurricane Harvey Stirs Up the Gulf of
  Mexico,} , 2017.
\newblock
  \urlprefix\url{https://earthobservatory.nasa.gov/images/90818/hurricane-harvey-stirs-up-the-gulf-of-mexico},
  last Accessed January 17, 2024.

\bibitem[{Blake and Zelinsky(2018)}]{Blake2017}
Blake, E.~S., and Zelinsky, D.~A., \enquote{National Hurricane Center Tropical
  Cyclone Report Hurricane Harvey (AL092017) 17 August 2017 - 1 September
  2017,} \emph{National Hurricane Center}, 2018.

\bibitem[{Li et~al.(2022)Li, Goldberg, Kalluri, Lindsey, Sjoberg, Zhou,
  Helfrich, Green, Borges, Yang, and Sun}]{Li2022HarveyFlooding}
Li, S., Goldberg, M., Kalluri, S., Lindsey, D.~T., Sjoberg, B., Zhou, L.,
  Helfrich, S., Green, D., Borges, D., Yang, T., and Sun, D., \enquote{High
  Resolution 3D Mapping of Hurricane Flooding from Moderate-Resolution
  Operational Satellites,} \emph{Remote Sensing}, Vol.~14, No.~21, 2022.
\newblock \doi{10.3390/rs14215445}.

\bibitem[{Reising et~al.(2018)Reising, Gaier, Padmanabhan, Lim, Heneghan,
  Kummerow, Berg, Chandrasekar, Radhakrishnan, Brown, Carvo, and
  Pallas}]{TEMPEST-8517330}
Reising, S.~C., Gaier, T.~C., Padmanabhan, S., Lim, B.~H., Heneghan, C.,
  Kummerow, C.~D., Berg, W., Chandrasekar, V., Radhakrishnan, C., Brown, S.~T.,
  Carvo, J., and Pallas, M., \enquote{An Earth Venture In-Space Technology
  Demonstration Mission for Temporal Experiment for Storms and Tropical Systems
  (Tempest),} \emph{IGARSS 2018 - 2018 IEEE International Geoscience and Remote
  Sensing Symposium}, 2018, pp. 6301--6303.
\newblock \doi{10.1109/IGARSS.2018.8517330}.

\bibitem[{Wilson et~al.(2018)Wilson, Angal, and Xiong}]{Wilson2018MODIS}
Wilson, T.~M., Angal, A., and Xiong, X., \enquote{Sensor Performance Assessment
  for Terra and Aqua Modis using unscheduled lunar observations,}
  \emph{Sensors, Systems, and Next-Generation Satellites XXII}, 2018.
\newblock \doi{10.1117/12.2324873}.

\bibitem[{Brown et~al.(2017)Brown, Focardi, Kitiyakara, Maiwald, Milligan,
  Montes, Padmanabhan, Redick, Russel, Bach, and Walkemeyer}]{COWVR-7943884}
Brown, S., Focardi, P., Kitiyakara, A., Maiwald, F., Milligan, L., Montes, O.,
  Padmanabhan, S., Redick, R., Russel, D., Bach, V., and Walkemeyer, P.,
  \enquote{The COWVR Mission: Demonstrating the capability of a new generation
  of small satellite weather sensors,} \emph{2017 IEEE Aerospace Conference},
  2017, pp. 1--7.
\newblock \doi{10.1109/AERO.2017.7943884}.

\bibitem[{{National Oceanic and Atmospheric
  Administration}(2009)}]{NOAAGlossary}
{National Oceanic and Atmospheric Administration}, \enquote{Glossary of NHC
  Terms,} , 2009.
\newblock \urlprefix\url{https://www.nhc.noaa.gov/aboutgloss.shtml}, last
  Accessed January 17, 2024.

\end{thebibliography}

\end{document}